\documentclass[10pt]{article}
\pdfoutput=1

\usepackage[activate={true, nocompatibility}, final, tracking=true, kerning=true, spacing=true, factor=1100, stretch=10, shrink=10]{microtype}
\microtypecontext{spacing=nonfrench}

\usepackage[T1]{fontenc}

\usepackage{authblk}

\usepackage{amsmath,amsfonts,amssymb,mathtools,amsthm}
\usepackage{empheq}
\usepackage{xfrac}
\usepackage{stmaryrd}
\usepackage{accents}
\usepackage{bm}

\usepackage{graphicx}
\usepackage{siunitx}
\usepackage{booktabs,multirow,multicol}
\usepackage{caption}
\usepackage{subcaption}
\usepackage[section]{placeins}

\usepackage{tikz}
\usetikzlibrary{matrix, shapes, arrows.meta, positioning, fit, backgrounds}
\usepackage{pgfplots}
\pgfplotsset{compat=newest}

\usepackage[frozencache]{minted}

\usepackage{natbib}
\usepackage[nottoc]{tocbibind}

\usepackage[linesnumbered, vlined, ruled]{algorithm2e}

\usepackage{url}
\usepackage[colorlinks=false]{hyperref}
\usepackage[noabbrev]{cleveref}
\hypersetup{pdfprintscaling=None}

\pgfmathsetlengthmacro{\myaxiswidth}{0.47\textwidth}
\pgfmathsetlengthmacro{\myaxisheight}{0.3\textwidth}
\pgfplotsset{
every axis/.append style={const plot, no markers, 
width=\myaxiswidth, 
height=\myaxisheight,
enlargelimits=false, scale only axis,
ymajorticks=false,
xmin=1, ymin=0, ymax=1,
xticklabel style={anchor=south east},
},
every axis plot/.append style={mark=none, line width=2, draw opacity=0.5},
}

\newcommand{\ubar}[1]{\underaccent{\bar}{#1}}
\newcommand{\utilde}[1]{\underaccent{\tilde}{#1}}

\crefname{equation}{}{}

\SetKwComment{Comment}{$\triangleright$\ }{}
\SetCommentSty{}
\SetKw{Continue}{continue}
\SetKw{Break}{break}

\makeatletter
\newenvironment{model}[1]{
    \subequations
    \def\@currentlabel{#1}
}{
    \endsubequations
}
\makeatother

\makeatletter
\g@addto@macro\@floatboxreset\centering
\makeatother

\let\emptyset\varnothing

\DeclareMathOperator{\cl}{cl}

\DeclareMathOperator{\rank}{rank}

\DeclareMathOperator{\cone}{cone}
\DeclareMathOperator{\dom}{dom}

\DeclareMathOperator{\smat}{smat}
\DeclareMathOperator{\svec}{svec}

\newcommand{\bbR}[0]{\mathbb{R}}
\newcommand{\bbZ}[0]{\mathbb{Z}}
\newcommand{\bbS}[0]{\mathbb{S}}

\newcommand{\iI}[0]{\llbracket I \rrbracket}
\newcommand{\iL}[0]{\llbracket L \rrbracket}
\newcommand{\iK}[0]{\llbracket K \rrbracket}

\newcommand{\iJk}[0]{\llbracket J_k \rrbracket}
\newcommand{\iNs}[0]{\llbracket n \rrbracket}

\newcommand{\cX}[0]{\mathcal{X}}
\newcommand{\cN}[0]{\mathcal{N}}
\newcommand{\cZ}[0]{\mathcal{Z}}

\newcommand{\K}[0]{\mathcal{K}}
\newcommand{\KZ}[0]{\{0\}}
\newcommand{\KL}[0]{\mathcal{L}}
\newcommand{\KR}[0]{\mathcal{V}}
\newcommand{\KP}[0]{\mathcal{S}}
\newcommand{\KE}[0]{\mathcal{E}}

\newcommand{\be}[0]{\bm{e}}
\newcommand{\bt}[0]{\bm{t}}
\newcommand{\bT}[0]{\bm{T}}
\newcommand{\btau}[0]{\bm{\tau}}
\newcommand{\bw}[0]{\bm{w}}
\newcommand{\bom}[0]{\bm{\omega}}
\newcommand{\bW}[0]{\bm{W}}
\newcommand{\bx}[0]{\bm{x}}
\newcommand{\by}[0]{\bm{y}}
\newcommand{\bz}[0]{\bm{z}}
\newcommand{\bl}[0]{\bm{l}}
\newcommand{\bu}[0]{\bm{u}}
\newcommand{\til}[0]{\utilde{l}}
\newcommand{\tiu}[0]{\utilde{u}}
\newcommand{\tibl}[0]{\utilde{\bm{l}}}
\newcommand{\tibu}[0]{\utilde{\bm{u}}}
\newcommand{\bpi}[0]{\bm{\pi}}
\newcommand{\bc}[0]{\bm{c}}
\newcommand{\bb}[0]{\bm{b}}
\newcommand{\bA}[0]{\bm{A}}

\newcommand{\bmu}[0]{\bm{\mu}}
\newcommand{\bnu}[0]{\bm{\nu}}

\newcommand{\micp}[0]{\mathfrak{M}}
\newcommand{\conic}[0]{\mathfrak{C}}
\newcommand{\poly}[0]{\mathfrak{P}}

\begin{document}

\title{Outer Approximation With Conic Certificates For Mixed-Integer Convex Problems%
\thanks{M. Lubin is now at Google. The authors thank MOSEK for access to alpha releases of MOSEK 9 and Russell Bent and Emre Yamangil for support for and early contributions to Pajarito. This work has been partially funded by the National Science Foundation under grant CMMI-1351619, the Department of Energy Computational Science Graduate Fellowship under grant DE-FG02-97ER25308, and the Office of Naval Research under grant N00014-18-1-2079.}}

\author[$\S$]{Chris Coey}
\author[$\S$]{Miles Lubin}
\author[$\ddag$]{Juan Pablo Vielma}
\affil[$\S$]{MIT Operations Research Center}
\affil[$\ddag$]{MIT Sloan School of Management}
\affil[ ]{\texttt{coey@mit.edu, miles.lubin@gmail.com, jvielma@mit.edu}}

\maketitle

\begin{abstract}
A mixed-integer convex (MI-convex) optimization problem is one that becomes convex when all integrality constraints are relaxed. We present a branch-and-bound LP outer approximation algorithm for an MI-convex problem transformed to \textit{MI-conic} form. 
The polyhedral relaxations are refined with \textit{$\K^*$ cuts} derived from \textit{conic certificates} for continuous primal-dual conic subproblems. 
Under the assumption that all subproblems are \textit{well-posed}, the algorithm detects infeasibility or unboundedness or returns an optimal solution in finite time. 
Using properties of the conic certificates, we show that the $\K^*$ cuts imply certain practically-relevant guarantees about the quality of the polyhedral relaxations, and demonstrate how to maintain helpful guarantees when the LP solver uses a positive feasibility tolerance. 
We discuss how to \textit{disaggregate} $\K^*$ cuts in order to tighten the polyhedral relaxations and thereby improve the speed of convergence, and propose fast heuristic methods of obtaining useful $\K^*$ cuts.
Our new open source MI-conic solver \textit{Pajarito} (\href{http://github.com/JuliaOpt/Pajarito.jl}{github.com/JuliaOpt/Pajarito.jl}) uses an external mixed-integer linear (MILP) solver to manage the search tree and an external continuous conic solver for subproblems. 
Benchmarking on a library of mixed-integer second-order cone (MISOCP) problems, we find that Pajarito greatly outperforms Bonmin (the leading open source alternative) and is competitive with CPLEX's specialized MISOCP algorithm. 
We demonstrate the robustness of Pajarito by solving diverse MI-conic problems involving mixtures of positive semidefinite, second-order, and exponential cones, and provide evidence for the practical value of our analyses and enhancements of $\K^*$ cuts.
\end{abstract}

\newpage
\setcounter{tocdepth}{2}
\tableofcontents


\newpage
\section{Mixed-Integer Convex Optimization}
\label{sec:intro}

A \textit{mixed-integer convex} (MI-convex) problem is a finite-dimensional optimization problem that minimizes a convex objective function over convex constraints and integrality restrictions on a subset of the variables. 
\citet{Belotti2013a} and \citet{Bonami2012} review MI-convex applications and \citet{Lubin2016b} characterize which nonconvex feasible regions are MI-convex-representable. 
Since an MI-convex problem without integrality restrictions is just a convex problem, MI-convex optimization generalizes both mixed-integer linear optimization (MILP) and convex optimization.
This structure also leads to effective branch-and-bound (B\&B) algorithms, which recursively partition the possible values of the integer variables in a search tree and obtain objective bounds and feasible solutions from efficiently-solvable subproblems.

\subsection{Branch-And-Bound Algorithms}
\label{sec:intro:alg}

A \textit{nonlinear B\&B} (B\&B-NL) algorithm for a MI-convex problem solves a nonlinear subproblem that includes all of the convex constraints at every node of the search tree.
The Bonmin solver package \citep{Bonami2008} implements a B\&B-NL variant by calling the derivative-oracle-based nonlinear programming (NLP) solver Ipopt to solve the subproblems. 
The relatively new SCIP-SDP \citep{Gally2016} B\&B-NL implementation for mixed-integer semidefinite (MISDP) problems uses a primal-dual conic interior-point solver for the SDP subproblems. 

Typically, B\&B-NL methods need to solve a large number of very similar nonlinear subproblems to near-global optimality and feasibility in order to obtain accurate objective bounds.
Linear optimization (LP) solvers based on the Simplex algorithm are able to rapidly reoptimize after variable bounds are changed or linear cuts are added, thus typically benefiting from warm-starting much more so than state-of-the-art NLP or conic solvers. 
\textit{B\&B LP outer approximation} (B\&B-OA) algorithms take advantage of LP warm-starting by solving a polyhedral relaxation of the nonlinear subproblem at every node. Implementations often take advantage of the speed and stability of advanced MILP branch-and-cut solvers.

B\&B-OA algorithms differ in how they refine the polyhedral relaxations of the nonlinear constraints and how they obtain feasible solutions. 
In a \textit{separation-based} algorithm, no nonlinear solver is used. At each node the LP optimal point is first checked for feasibility for the convex constraints; if the violation exceeds a positive tolerance, valid cuts separating the point are added to the LP, otherwise the point may be accepted as a new incumbent if it is integral.%
\footnote{Commercial mixed-integer second-order cone optimization (MISOCP) solvers use separation-based algorithms, but also occasionally solve SOCP subproblems to obtain feasible solutions and fathom nodes. SCIP-SDP offers both B\&B-NL and separation-based B\&B-OA methods for mixed-integer semidefinite problems.}
\citet{Quesada1992} and \citet{Leyffer1993} describe \textit{subproblem-based} B\&B-OA algorithms that solve smooth subproblems at a subset of the nodes. The subproblems provide points at which cuts based on gradient inequalities can be derived.%
\footnote{For a convex function $f : \bbR^n \to \bbR$, the set $\cX = \{\bx \in \bbR^n : f(\bx) \leq 0\}$ is convex. If $f$ is smooth, then given a point $\bar{\bx} \in \bbR^n$, the following \textit{gradient cut} yields a polyhedral relaxation of $\cX$:
\begin{equation}
f (\bar{\bx}) + (\nabla f (\bar{\bx}))^T (\bx - \bar{\bx}) \leq 0
\label{con:gradcut}
.
\end{equation}}

\citet{Bonami2008} found that Bonmin's B\&B-OA method generally outperforms its B\&B-NL method. 
Since both of these methods rely on NLP subproblems, they frequently fail in the presence of nonsmoothness. 
Continuous conic solvers are more numerically robust than derivative-oracle-based NLP solvers on nonsmooth problems (such as SOCPs and SDPs). 
For the special case of MISOCP, \citet{Drewes2012} propose a conic subproblem-based B\&B-OA algorithm that derives cuts from subgradients satisfying subproblem KKT optimality conditions, and hence does not require smoothness assumptions.

Another advantage of conic solvers is that they return simple \textit{certificates} proving primal or dual infeasibility or optimality of a primal-dual solution pair \citep{Friberg2015}.
Using the theory of conic duality, it is possible to describe an elegant OA algorithm for generic MI-conic problems that uses conic certificates returned by primal-dual conic solvers, with no need to examine KKT conditions or solve a second modified subproblem in the case of infeasibility (as in the algorithm by \citet{Drewes2012}).
\citet{Lubin2016} propose this idea in an iterative OA algorithm. 
However a B\&B algorithm using a single search tree, instead of solving a sequence of MILP instances each with their own search tree, is more flexible and likely to be significantly faster in practice. 
We fill this gap with the first conic-certificate-based B\&B-OA algorithm.

\subsection{Mixed-Integer Conic Form}
\label{sec:intro:conic}

We use the following general form for a mixed-integer conic (MI-conic) problem:
\begin{model}{$\micp$}
\label{mod:micp}
\begin{empheq}[left=\micp\empheqlbrace\quad]{alignat=3}
\inf_{\mathclap{\bx \in \bbR^N}} \quad \bc^T \bx & : && \qquad
\label{obj:m} 
\\
\bb - \bA \bx & \in \K \subset \bbR^M
\label{con:mk} 
\\
x_i & \in \bbZ &&& \forall i & \in \iI
,
\label{con:mi}
\end{empheq}
\end{model}
where $\K$ is a \textit{closed convex cone}, i.e. a closed subset of $\bbR^M$ that contains all conic (nonnegative) combinations of its points \citep{Ben-Tal2001}:
\begin{equation}
\alpha_1 \by_1 + \alpha_2 \by_2 \in \K \qquad \forall \alpha_1, \alpha_2 \geq 0 \qquad \forall \by_1, \by_2 \in \K
\label{eq:cone}
.
\end{equation} 
The decision variables in \ref{mod:micp} are represented by the column vector $\bx \in \bbR^N$, so the objective \cref{obj:m} minimizes a linear function of $\bx$ subject to (denoted by `:') the constraints \cref{con:mk,con:mi}. 
The index set of integer decision variables is $\iI = \{1, \ldots, I\}$, so the integrality constraints \cref{con:mi} restrict only the first $I$ variables $x_1, \ldots, x_I$ to the set of integers $\bbZ$. 
The conic constraint \cref{con:mk}, which restricts the affine transformation $\bb - \bA \bx$ of $\bx$ to $\K$, is a convex constraint, so relaxing the integrality constraints \cref{con:mi} results in a convex conic optimization problem. 

Any MI-convex problem can be expressed in MI-conic form, by homogenizing the convex constraints, for example through perspective transformations \citep{Boyd2004, Lubin2016}. 
\textit{Disciplined Convex Programming} (DCP) modeling packages such as CVX \citep{Grant2014}, CVXPy \citep{Diamond2016}, and Convex.jl \citep{Udell2014} perform conic transformations automatically, conveniently enabling modelers to access powerful conic solvers such as ECOS \citep{Domahidi2013}, SCS \citep{ODonoghue2016}, MOSEK \citep{MOSEK2016}, and CSDP \citep{Borchers1999}. 

Conic solvers recognize the cone $\K$ as a Cartesian product of standard \textit{primitive} cones \citep{Friberg2016}. A primitive closed convex cone cannot be written as a Cartesian product of two or more lower-dimensional closed convex cones. 
\citet{Lubin2017} claim that the following classes of standard primitive cones are extremely versatile, encoding all of the problems in the Conic Benchmark Library (CBLIB) compiled by \citet{Friberg2016}, and all 333 MI-convex problems in MINLPLIB2 \citep{GAMS2018}.
\begin{description}
\item[Linear cones] naturally express affine constraints; any mixed-integer linear optimization (MILP) problem can be written in \ref{mod:micp} form using nonnegative, nonpositive, and zero cones.
\item[Second-order cones] (and rotated-second-order cones) are widely used to model rational powers, norms, and geometric means, as well as convex quadratic objectives and constraints \citep{Ben-Tal2001}. 
\item[Positive semidefinite cones] can model robust norms and functions of eigenvalues \citep{Ben-Tal2001}, and sum-of-squares constraints for polynomial optimization problems \citep{Parrilo2003}.
\item[Exponential cones] can model exponentials, logarithms, entropy, and powers, as well as log-sum-exp functions that arise from convex transformations of geometric programs \citep{Serrano2015}.
\end{description}
We note that conic representations are useful for constructing tight formulations for disjunctions or unions of convex sets \citep{Lubin2016b,Lubin2017a,Vielma2017}.

\subsection{Overview And Contributions}
\label{sec:intro:over}

In \cref{sec:alg}, we start by reviewing the relevant foundations of conic duality and certificates. 
We then introduce the notion of $\K^*$ cuts, and describe how to refine LP OAs of conic constraints using certificates obtained from continuous primal-dual conic solvers. 
For a MI-conic problem \ref{mod:micp}, we propose the first B\&B-OA algorithm based on conic certificates. 
We show that our algorithm detects infeasibility or unboundedness or terminates with an optimal solution in finite time under minimal assumptions.

In \cref{sec:guar}, we demonstrate that a $\K^*$ cut from a conic certificate implies useful guarantees about the infeasibility or optimal objective of an LP OA, suggesting that our algorithm can often fathom a node immediately after solving the LP rather than proceeding to the expensive conic subproblem solve.
We consider how these guarantees may be lost in the more realistic setting of an LP solver with a positive feasibility tolerance, and propose a practical methodology for scaling a certificate $\K^*$ cut to recover similar guarantees.

In \cref{sec:tight}, we describe how to strengthen the LP OAs by \textit{disaggregating} $\K^*$ cuts, and show that this methodology maintains the guarantees from \cref{sec:guar}.
We argue for initializing the LP OAs using \textit{initial fixed $\K^*$ cuts}, and offer a procedure for cheaply obtaining \textit{separation $\K^*$ cuts} to cut off an infeasible LP OA solution.
These proposed techniques require minimal modifications to our algorithm and are practical to implement.
In \cref{sec:spec}, we specialize these techniques for the second-order, positive semidefinite, and exponential cones. 

In \cref{sec:soft}, we describe the software architecture and algorithmic implementation of Pajarito, our open source MI-convex solver.%
\footnote{The new version of Pajarito that we implemented for this paper is the first conic-certificate-based OA solver. Although Pajarito solver was introduced in \citet{Lubin2016}, this early implementation used NLP solvers instead of primal-dual conic solvers for continuous subproblems, and was built to assess the value of \textit{extended formulations} by counting iterations before convergence. To avoid confusing users, we recently moved this old NLP-based functionality out of Pajarito and into \textit{Pavito} solver at \href{https://github.com/JuliaOpt/Pavito.jl}{github.com/JuliaOpt/Pavito.jl}.}
This section may be of particular interest to advanced users and developers of mathematical optimization software. 
We emphasize that our implementations diverge from the idealized algorithmic description in \cref{sec:alg}, because of our decision to leverage powerful external mixed-integer linear (MILP) solvers through limited, solver-independent interfaces.
In \cref{sec:socef}, we describe how Pajarito lifts $\K^*$ cuts for the second-order cone using an \textit{extended formulation}, resulting in tighter LP OAs.
In \cref{sec:psdsoc}, we show how Pajarito can optionally tighten OAs for PSD cone constraints by strengthening $\K^*$ cuts to rotated second-order cone constraints, which can be added to an MISOCP OA model.

In \cref{sec:exp}, we summarize computational experiments demonstrating the speed and robustness of Pajarito. 
We benchmark Pajarito and several MISOCP solver packages, and conclude that Pajarito is the fastest and most-reliable open source solver for MISOCP. 
Finally, we compare the performance of several of Pajarito's algorithmic variants on MI-conic instances involving mixtures of positive semidefinite, second-order, and exponential cones, demonstrating practical advantages of several methodological contributions from \cref{sec:guar,sec:tight} and \cref{sec:spec}. 

\section{A Branch-And-Bound LP Outer Approximation Algorithm}
\label{sec:alg}

For a MI-conic problem \ref{mod:micp}, we propose a branch-and-bound LP outer approximation (B\&B-OA) algorithm, the first such method based on conic certificates. 
In \cref{sec:alg:sub}, we describe the continuous conic subproblems that a nonlinear branch-and-bound (B\&B-NL) algorithm would solve at each node, and review the relevant foundations of conic duality. 
In \cref{sec:alg:poly}, we introduce the notion of $\K^*$ cuts and describe how to refine polyhedral relaxations of the conic subproblems using information from conic certificates. 
Finally, we outline our B\&B-OA algorithm in \cref{sec:alg:alg}, and discuss finiteness of convergence.

\subsection{Continuous Subproblems And Conic Duality}
\label{sec:alg:sub}

Recall from \ref{mod:micp} that the first $I$ variables in $\bx$ are constrained to be integer. 
Branch-and-bound algorithms recursively partition the valid integer assignments, so for convenience we assume known finite lower bounds $\bl^0 \in \bbZ^I$ and upper bounds $\bu^0 \in \bbZ^I$ on the integer variables $x_1, \ldots x_I$.
At a node of the branch-and-bound search tree with lower bounds $\bl \in \bbZ^I$ and upper bounds $\bu \in \bbZ^I$ on $x_1, \ldots, x_I$, the natural continuous conic subproblem is \ref{mod:conic}: 
\begin{model}{$\conic (\bl, \bu)$}
\label{mod:conic}
\begin{empheq}[left = {\conic (\bl, \bu)} \empheqlbrace \quad]{alignat=3}
\inf_{\mathclap{\bx}} \quad \bc^T \bx & : && \qquad
\\
\bb - \bA \bx & \in \K
\label{con:conic}
\\
l_i - x_i & \in \bbR_- &&& \forall i & \in \iI 
\label{con:low}
\\
u_i - x_i & \in \bbR_+ &&& \forall i & \in \iI 
\label{con:up}
\\
\bx & \in \bbR^N
\label{con:xR}
,
\end{empheq}
\end{model}
where the bound constraints \cref{con:low,con:up} are expressed in conic form using the nonpositive cone $\bbR_-$ (the nonpositive reals) and the nonnegative cone $\bbR_+$ (the nonnegative reals). 

There exist primal-dual conic algorithms for \ref{mod:conic} that are powerful in both theory and practice. 
The foundation for these methods and for much of this paper is the elegant theory of conic duality, described by \citet{Ben-Tal2001,Boyd2004}. 
Recall that the cone $\K$ in \ref{mod:conic} is a closed convex cone; we let $\K^*$ denote the \textit{dual cone} of $\K$, i.e. the set of points that have nonnegative inner product with all points in $\K$:
\begin{equation}
\K^* = \{ \bz \in \bbR^M : \by^T \bz \geq 0, \forall \by \in \K \}
\label{eq:kd}
.
\end{equation}
$\K^*$ is also a closed convex cone \citep{Boyd2004}. The standard conic dual of \ref{mod:conic} can be written as \ref{mod:conic*}:
\begin{model}{$\conic^* (\bl, \bu)$}
\label{mod:conic*}
\begin{empheq}[left = {\conic^* (\bl, \bu)} \empheqlbrace \quad]{alignat=3}
\sup_{\mathclap{\bz, \bmu, \bnu}} \quad - \bb^T \bz - \bl^T \bmu - \bu^T \bnu & : && \qquad
\\
\bz & \in \K^*
\label{con:dualconic}
\\
\bmu & \in \bbR_-^I
\label{con:duallow}
\\
\bnu & \in \bbR_+^I
\label{con:dualup}
\\
\bc + \bA^T \bz + \bmu' + \bnu' & \in \KZ^N
\label{con:dualxR}
,
\end{empheq}
\end{model}
where for ease of exposition we let $\bmu' = (\mu_1, \ldots, \mu_I, 0, \ldots, 0) \in \bbR^N$ and similarly for $\bnu'$. 
Note that the nonnegative and nonpositive cones are both \textit{self-dual}, i.e. $\bbR^*_- = \bbR_-$ and $\bbR^*_+ = \bbR_+$. The zero cone $\KZ$ (containing only the origin) is dual to the free cone $\bbR$. 
The variables $\bz$ in the dual constraint \cref{con:dualconic} are associated with the primal constraint \cref{con:conic}, and similarly for \cref{con:duallow} and \cref{con:low}, \cref{con:dualup} and \cref{con:up}, and \cref{con:xR} and \cref{con:dualxR}. 

If the conic primal-dual pair \ref{mod:conic}--\ref{mod:conic*} is \textit{well-posed}, then conic duality can be thought of as a simple generalization of LP duality.%
\footnote{If $\K$ is polyhedral, then $\K^*$ is polyhedral, and hence $\conic (\bl, \bu)$ and $\conic^* (\bl, \bu)$ are both LPs. All LPs are well-posed.}
In particular, the $\inf$ and $\sup$ can be replaced with $\min$ and $\max$, and the possible status combinations for \ref{mod:conic} and \ref{mod:conic*} are: both infeasible, one unbounded and the other infeasible, or both feasible and bounded with equal objective values attained by optimal solutions. 
The conditions for well-posedness in conic duality are described by \citet{Friberg2016}, and are outside the scope of this paper, so we assume that any primal-dual subproblem we encounter has the well-posed property. 

\citet{Friberg2016} discusses certificates that provide easily-verifiable proofs of unboundedness or infeasibility of the primal or dual problems or of optimality of a given pair of primal and dual points. In terms of the primal subproblem \ref{mod:conic}, the three possible mutually-exclusive cases and their interpretations are as follows.
\begin{description}
\item[A dual improving ray] certifies that \ref{mod:conic} is infeasible, via the conic generalization of Farkas' lemma. 
The improving ray $(\bar{\bz}, \bar{\bmu}, \bar{\bnu}) \in \bbR^{M+2I}$ of \ref{mod:conic*} is a feasible direction for \ref{mod:conic*} along which the objective value of any feasible point of \ref{mod:conic*} can be improved indefinitely. 
It satisfies the following conditions:
\begin{subequations}
\begin{alignat}{1}
-\bb^T \bar{\bz} - \bl^T \bar{\bmu} - \bu^T \bar{\bnu} & > 0 
\label{eq:certinf:obj}
\\
\bar{\bz} & \in \K^*
\label{eq:certinf:z}
\\
\bar{\bmu} & \leq \bm{0}
\\
\bar{\bnu} & \geq \bm{0}
\\
\bA^T \bar{\bz} + \bar{\bmu}' + \bar{\bnu}' & = \bm{0}
\label{eq:certinf:ray}
.
\end{alignat}
\end{subequations}
Clearly, if \ref{mod:conic*} itself has a feasible point, then it is unbounded, otherwise it is infeasible.%
\footnote{Conditions \crefrange{eq:certinf:z}{eq:certinf:ray} imply that $(\bar{\bz}, \bar{\bmu}, \bar{\bnu})$ is feasible for a modified \ref{mod:conic*} problem in which $\bc = \bm{0}$.}
\item[A primal improving ray and a primal feasible point] certifies that \ref{mod:conic} is unbounded, because the improving ray is a feasible direction along which the objective value of the feasible point can be improved indefinitely. 
The improving ray $\bar{\bx} \in \bbR^N$ of \ref{mod:conic} also implies infeasibility of \ref{mod:conic*} and satisfies the following conditions:
\begin{subequations}
\begin{alignat}{3}
\bc^T \bar{\bx} & < 0 && \qquad
\label{eq:certunb:start}
\\
- \bA \bar{\bx} & \in \K
\\
\bar{x}_i & = 0 &&& \forall i & \in \iI 
\label{eq:certunb:end}
,
\end{alignat}
\end{subequations}
and the feasible point $\hat{\bx} \in \bbR^N$ of \ref{mod:conic} simply satisfies the primal feasibility conditions \crefrange{con:conic}{con:xR}.
Note that if the objective coefficients of the continuous variables are all zero ($c_{I+1} = \ldots = c_N = 0$), then conditions \cref{eq:certunb:start,eq:certunb:end} can never be satisfied, so there cannot be a primal improving ray.%
\footnote{This matches intuition because if the continuous variables have zero objective coefficients and the integer variables are bounded, \ref{mod:micp} cannot be unbounded.}
\item[A complementary solution pair] certifies optimality for \ref{mod:conic} of the primal feasible point $\hat{\bx} \in \bbR^N$ in the pair $(\hat{\bx}, (\hat{\bz}, \hat{\bmu}, \hat{\bnu}))$, via conic weak duality. 
The dual feasible point $(\hat{\bz}, \hat{\bmu}, \hat{\bnu}) \in \bbR^{M+2I}$ is also optimal for \ref{mod:conic*}, and the pair have equal primal and dual objective values:
\begin{equation}
\bc^T \hat{\bx} = -\bb^T \hat{\bz} - \bl^T \hat{\bmu} - \bu^T \hat{\bnu}
\label{eq:certsd}
.
\end{equation}
\end{description}

\subsection{Dynamic Polyhedral Relaxations}
\label{sec:alg:poly}

Recall from equation \cref{eq:kd} that $\by \in \K$ if and only if $\bz^T \by \geq 0, \forall \bz \in \K^*$. 
This implies that a nonpolyhedral conic constraint $\bb - \bA \bx \in \K$ has the following equivalent semi-infinite linear representation:
\begin{equation}
\bz^T ( \bb - \bA \bx ) \geq 0 \qquad \forall \bz \in \K^*
.
\end{equation}
We refer to a point $\bz \in \K^*$ as a \textit{$\K^*$ point}, and call the corresponding linear constraint $\bz^T ( \bb - \bA \bx ) \geq 0$ a \textit{$\K^*$ cut}. 
A $\K^*$ cut cannot exclude any point $\bx$ that satisfies $\bb - \bA \bx \in \K$, so any finite set of $\K^*$ cuts defines a valid polyhedral relaxation of the conic constraint \cref{con:conic}. 

Given a finite set $\cZ \subset \K^*$ of $\K^*$ points, consider modifying the subproblem \ref{mod:conic} by relaxing the conic constraint and instead imposing the finite number of $\K^*$ cuts implied by $\cZ$. 
We refer to the resulting LP as \ref{mod:poly}, which we choose to write in inequality form rather than conic form:
\begin{model}{$\poly (\cZ, \bl, \bu)$}
\label{mod:poly}
\begin{empheq}[left = {\poly (\cZ, \bl, \bu)} \empheqlbrace \quad]{alignat=3}
\min_{\mathclap{\bx}} \quad \bc^T \bx & : && \qquad
\\
x_i & \geq l_i &&& \forall i & \in \iI
\label{con:poly:l}
\\
x_i & \leq u_i &&& \forall i & \in \iI
\label{con:poly:u}
\\
\bz^T ( \bb - \bA \bx ) & \geq 0 &&& \forall \bz & \in \cZ
\label{con:poly:k}
.
\end{empheq}
\end{model}

Since \ref{mod:poly} and \ref{mod:conic} have the same objective function, and the feasible set of \ref{mod:poly} is a polyhedral relaxation of the feasible set of \ref{mod:conic}, solving \ref{mod:poly} with an LP solver may give us useful information about \ref{mod:conic}. 
If \ref{mod:poly} is infeasible, then \ref{mod:conic} must be infeasible. If \ref{mod:poly} has an optimal objective value of $L$, then \ref{mod:conic} is either infeasible or has an optimal objective no smaller than $L$. 
In these cases, we may be able to immediately fathom the node by infeasibility or by bound, or even use a fractional optimal solution for \ref{mod:poly} to make a branching decision, without needing to solve \ref{mod:conic}. 
However, if \ref{mod:poly} is unbounded, it does not provide useful information about the status or optimal value of \ref{mod:conic}. 

An LP solver based on the Simplex algorithm is able to rapidly reoptimize \ref{mod:poly} at each new node, after the integer variable bounds are updated and any new $\K^*$ cuts are added. 
As noted by \citet{Skajaa2013}, state-of-the-art conic solvers benefit much less from warm-starting, so it may be computationally faster to sacrifice some information about the conic subproblems in order to avoid some expensive conic subproblem solves. 

In analogy to B\&B-OA algorithms based on gradient cuts, we add a new cut after every infeasible or bounded conic subproblem solve. Our key innovation, however, is to obtain this cut directly from the conic certificate found by the conic subproblem solver.
Suppose that at some node, a primal-dual conic subproblem solver yields a dual improving ray $(\bar{\bz}, \bar{\bmu}, \bar{\bnu})$: from condition \cref{eq:certinf:z}, $\bar{\bz} \in \K^*$, so $\bar{\bz}$ is a $\K^*$ point. 
Now suppose that the subproblem solver yields a complementary solution $(\hat{\bx}, (\hat{\bz}, \hat{\bmu}, \hat{\bnu}))$: by the dual feasibility condition \cref{con:dualconic}, $\hat{\bz} \in \K^*$, so $\hat{\bz}$ is a $\K^*$ point. 
In both cases, a subvector of the ray or solution for the dual subproblem \ref{mod:conic*} allows us to augment $\cZ \subset \K^*$, refining our LP OA model \ref{mod:poly}.
In \cref{sec:guar}, we use conic duality theory to show that these $\K^*$ cuts derived from certificates encode important information about conic subproblems into the subsequent polyhedral relaxations.

\subsection{The Conic-Certificate-Based Algorithm}
\label{sec:alg:alg}

Our conic-certificate-based B\&B-OA algorithm for the MI-conic problem \ref{mod:micp} is outlined in \cref{alg:oa}. 
Recall that \ref{mod:micp} is in minimization form. \Cref{alg:oa} maintains an upper bound $U$ (initially $\infty$), a corresponding best feasible solution set $\cX$ (initially empty), and a set of active nodes $\cN$ of the search tree. 
A node $(\bl, \bu, L)$ is characterized by the finite variable bound vectors $\bl$ and $\bu$ and a lower bound value $L$. 
The node's lower bound $L$ signifies that all feasible solutions for \ref{mod:micp} that satisfy the node's bounds on integer variables have an objective value of at least $L$. 
The node set $\cN$ is initialized to contain only the root node $(\bl^0, \bu^0, -\infty)$, where $\bl^0, \bu^0 \in \bbR^I$ are the finite initial global bounds on the integer variables.

\begin{algorithm}
\caption{Conic-certificate-based branch-and-bound LP outer approximation for \ref{mod:micp}.}
\label{alg:oa}
\DontPrintSemicolon
initialize incumbent solution set $\cX$ to $\emptyset$, upper bound $U$ to $\infty$\;
initialize $\K^*$ point set $\cZ$ to $\emptyset$\; \label{ln:zinit}
initialize node list $\cN$ with root node $(\bl^0, \bu^0, -\infty)$\;
\While{$\cN$ contains nodes}{
    remove a node $(\bl, \bu, L)$ from $\cN$\; \label{ln:rem}
    \If{lower bound $L \geq U$}{
        \Continue \Comment*[r]{fathomed by bound} \label{ln:fbnd}
    }
    call LP solver on \ref{mod:poly}\; \label{ln:lp}
    \uIf{get an infeasibility proof \label{ln:lpinf}}{
        \Continue \Comment*[r]{fathomed by infeasibility} \label{ln:lpfinf}
    }
    \ElseIf{get an optimal solution $\hat{\bx}$ \label{ln:lpsd}}{
        update $L$ to $\bc^T \hat{\bx}$\; \label{ln:lplb}
        \uIf{$L \geq U$}{
            \Continue \Comment*[r]{fathomed by bound} \label{ln:lpfbnd}
        }
        \ElseIf{$\hat{\bx}$ is fractional}{
            add branch nodes to $\cN$ using $\hat{\bx}$ and $L$\; \label{ln:lpbrnch0}
            \Continue \Comment*[r]{branched} \label{ln:lpbrnch}
        }
    }
    call primal-dual continuous conic solver on \ref{mod:conic}\; \label{ln:cp}
    \uIf{get a dual improving ray $(\bar{\bz}, \bar{\bmu}, \bar{\bnu})$ \label{ln:dir}}{
        add $\K^*$ point $\bar{\bz}$ to $\cZ$\; \label{ln:cutinf}
        \Continue \Comment*[r]{fathomed by infeasibility} \label{ln:finf}
    }
    \uElseIf{get a primal improving ray $\bar{\bx}$ and feasible point $\hat{\bx}$ \label{ln:pir}}{
        \If{$\hat{\bx}$ is integral \label{ln:intunb}}{
            update $U$ to $-\infty$\;
            \Break \Comment*[r]{proven unbounded} \label{ln:funb}
        }
    }
    \ElseIf{get a complementary solution $(\hat{\bx}, (\hat{\bz}, \hat{\bmu}, \hat{\bnu}))$ \label{ln:comp}}{
        add $\K^*$ point $\hat{\bz}$ to $\cZ$\; \label{ln:cutsd}
        update $L$ to $\bc^T \hat{\bx}$\; \label{ln:cplb}
        \uIf{$L \geq U$}{
            \Continue \Comment*[r]{fathomed by bound} \label{ln:fsdbd}
        }
        \ElseIf{$\hat{\bx}$ is integral \label{ln:intsd}}{
            update $\cX$ to $\{\hat{\bx}\}$ and $U$ to $\bc^T \hat{\bx}$\; \label{ln:upsdop}
            \Continue \Comment*[r]{fathomed by integrality} \label{ln:fsdop}
        }
    }
    add branch nodes to $\cN$ using $\hat{\bx}$ (fractional) and $L$\; \label{ln:brnch}
}
\Return $\cX$, $U$ \label{ln:ret}
\end{algorithm}

On \cref{ln:rem}, the main loop removes a node $(\bl, \bu, L)$ from $\cN$. If the node's lower bound $L$ is no smaller than the current global best upper bound $U$, \cref{ln:fbnd} fathoms the node by bound as it cannot yield a better incumbent. 
Otherwise, \cref{ln:lp} solves the node's LP OA model \ref{mod:poly}, taking advantage of an LP warm-start from a previous node. 

If \ref{mod:poly} is infeasible, \cref{ln:lpfinf} immediately fathoms the node by infeasibility. 
If \ref{mod:poly} has an optimal solution $\hat{\bx}$, then its optimal objective value is the tightest lower bound known for \ref{mod:conic} (in \cref{sec:guar:exact:sol}, we prove $\bc^T \hat{\bx} \geq L$ is a consequence of the $\K^*$ cuts), so \cref{ln:lplb} updates $L$ to $\bc^T \hat{\bx}$. 
\Cref{ln:lpfbnd} fathoms the node by bound if $L$ is no better than the incumbent value $U$, otherwise if $\hat{\bx}$ is fractional (i.e. it violates an integrality constraint \cref{con:mi}), \cref{ln:lpbrnch} branches on it.%
\footnote{We could instead remove lines 15-17 and solve the conic subproblem even if the LP solution is fractional, rather than branching. This variation may perform better if the conic subproblem solves are quite fast in practice.}
The branch procedure strictly partitions the node's integer bounds $\bl$ and $\bu$ by picking an $i \in \iI : \hat{x}_i \notin \bbZ$ and adding two child nodes to $\mathcal{N}$: $(\bl, (u_1, \ldots, \lfloor \hat{x}_i \rfloor, \ldots, u_N), L)$ and $((l_1, \ldots, \lceil \hat{x}_i \rceil, \ldots, l_N), \bu, L)$.

If the node is not fathomed or branched on immediately after the LP solve (before \cref{ln:cp}), then \ref{mod:poly} is either unbounded or has an optimal solution $\hat{\bx}$ that is integral (i.e. $\hat{x}_i \in \bbZ, \forall i \in \iI$) with optimal value $\bc^T \hat{\bx} < U$. 
Then \cref{ln:cp} solves the conic subproblem \ref{mod:conic} with the primal-dual continuous conic solver.
Recall from \cref{sec:alg:sub} our assumption that the primal-dual subproblem pair \ref{mod:conic}--\ref{mod:conic*} is well-posed, so the conic solver returns one of the three possible certificates, which we handle as follows.
\begin{description}
\item[A dual improving ray] on \cref{ln:dir} provides a $\K^*$ point, which \cref{ln:cutinf} adds to $\cZ$ (as described in \cref{sec:alg:poly}). 
This certificate proves that \ref{mod:conic} is infeasible, so \cref{ln:finf} fathoms the node by infeasibility. 
\item[A primal improving ray and feasible point] on \cref{ln:pir} certifies that \ref{mod:conic} is unbounded. 
Since the primal improving ray conditions \crefrange{eq:certunb:start}{eq:certunb:end} are the same for any conic subproblem, every subproblem is infeasible or unbounded, so \ref{mod:micp} is either infeasible or unbounded. 
The incumbent solution set must be empty and $U = \infty$. 
\Cref{ln:intunb} checks whether the feasible point $\hat{\bx}$ is integral. If so, it is a feasible solution for \ref{mod:micp}, so \ref{mod:micp} is unbounded and \cref{ln:funb} terminates the main loop.
\item[A complementary solution] on \cref{ln:comp} provides a $\K^*$ point that \cref{ln:cutsd} adds to $\cZ$ (as described in \cref{sec:alg:poly}) and an optimal solution $\hat{\bx}$ for \ref{mod:conic}. 
The optimal objective value gives the tightest lower bound known for the node, so \cref{ln:cplb} updates $L$ to $\bc^T \hat{\bx}$, and \cref{ln:fsdbd} fathoms by bound if this value is no better than $U$. 
\Cref{ln:intsd} checks if $\hat{\bx}$ is integral, in which case it becomes the new incumbent solution for \ref{mod:micp} on \cref{ln:upsdop}, and the node is fathomed by integrality on \cref{ln:fsdop}.
\end{description}

If the node is not fathomed immediately after the conic solve (before \cref{ln:brnch}), then $\hat{\bx}$ is a feasible solution for \ref{mod:conic} that is fractional. 
$L$ is either $\infty$ (in the primal improving ray case) or finite (in the complementary solution case), and is the best known lower bound for the node. \Cref{ln:brnch} branches on $\hat{\bx}$ using the same branch procedure we describe above for \cref{ln:lpbrnch}.

Since the initial bounds on the integer variables are finite, and the main loop of \cref{alg:oa} either fathoms each node or strictly partitions its integer bounds or terminates the algorithm, it follows that \cref{alg:oa} terminates finitely. 
From the fact that \ref{mod:poly} is a valid polyhedral relaxation of \ref{mod:conic}, and from the correctness of our inferences from the subproblem certificates, it is clear \cref{alg:oa} terminates correctly, under the assumption of well-posed conic subproblems. 
On \cref{ln:ret}, if $U = \infty$, then \ref{mod:micp} is proven infeasible, otherwise if $U$ is finite, then $\cX$ contains an optimal solution for \ref{mod:micp}, otherwise $U = -\infty$ and \ref{mod:micp} is proven unbounded. 

We note that without using the LP \ref{mod:poly} (i.e. removing \crefrange{ln:lp}{ln:lpbrnch} and not creating and augmenting $\cZ$ on \cref{ln:zinit,ln:cutinf,ln:cutsd}), we get a simple conic-certificate-based B\&B-NL algorithm for \ref{mod:conic}, for which finite termination guarantees and correctness follow from the same assumptions and arguments.
We have omitted any discussion of node selection or fractional variable selection for branching. 
MILP solvers can use LP certificates to make intelligent selections, and we expect that some of these LP-based criteria are generalizeable to the conic case, as conic duality theory is a simple extension of LP duality under the well-posed assumption.
\section{Polyhedral Relaxation Guarantees From Conic Certificates}
\label{sec:guar}

Recall from \cref{sec:alg:poly} that $\K^*$ cuts yield \textit{valid} polyhedral relaxations of the conic constraint $\bb - \bA \bx \in \K$, and a \textit{certificate $\K^*$ cut} can be obtained directly from the conic certificate for an infeasible or bounded and feasible subproblem \ref{mod:conic}. 
We demonstrate in \cref{sec:guar:exact} that a certificate $\K^*$ cut implies useful \textit{guarantees} about the infeasibility or optimal objective of the LP OAs, suggesting that \cref{alg:oa} can often fathom a node immediately after solving the LP OA rather than proceeding to the expensive conic subproblem solve.%
\footnote{By similar arguments, we expect that the certificate $\K^*$ cut may be useful at nearby nodes for duality based prepossessing such as reduced cost fixing \citep[sec.~7]{Gally2016} or conflict analysis \citep{Witzig2017}.}
In \cref{sec:guar:tol}, we consider how these guarantees may be lost in the more realistic setting of an LP solver with a positive feasibility tolerance, and propose a practical methodology for \textit{scaling} a certificate $\K^*$ cut to recover similar guarantees.

\subsection{Under An Exact LP Solver}
\label{sec:guar:exact}

We continue to assume well-posedness of every conic subproblem at every node. 
We consider what a certificate $\K^*$ from the conic subproblem \ref{mod:conic} at a node with integer variable bounds $\bl, \bu$ implies about the LP OA \hyperref[mod:poly]{$\poly (\cZ, \tibl, \tibu)$} at a different node with bounds $\tibl, \tibu$.

\subsubsection{Certificate Cuts From Dual Improving Rays}
\label{sec:guar:exact:ray}

Suppose $(\bar{\bz}, \bar{\bmu}, \bar{\bnu})$ is an improving ray of the dual subproblem \ref{mod:conic*}, certifying infeasibility of \ref{mod:conic}.
Using properties \crefrange{eq:certinf:obj}{eq:certinf:ray} of this certificate, any point $\bx \in \bbR^N$ satisfying the integer variable bounds $\til_i \leq x_i \leq \tiu_i, \forall i \in \iI$ at the new node also satisfies:
\begin{subequations}
\begin{align}
\bar{\bz}^T (\bb - \bA \bx) 
\label{eq:guar:ray1}
& = \bb^T \bar{\bz} - \bx^T \bA^T \bar{\bz}
\\
& = \bb^T \bar{\bz} + \bx^T \bar{\bmu}' + \bx^T \bar{\bnu}'
\\
& \leq \bb^T \bar{\bz} + \bx^T \bar{\bmu}' + \bx^T \bar{\bnu}' + \sum_{\mathclap{i \in \iI}} \, ( (\til_i - x_i) \bar{\mu}_i + (\tiu_i - x_i) \bar{\nu}_i )
\label{eq:guar:ray2}
\\
& = \bb^T \bar{\bz} + \tibl^T \bar{\bmu} + \tibu^T \bar{\bnu}
\label{eq:guar:ray3}
\\
& = (\bb^T \bar{\bz} + \bl^T \bar{\bmu} + \bu^T \bar{\bnu}) - (\bl - \tibl)^T \bar{\bmu} - (\bu - \tibu)^T \bar{\bnu}
\label{eq:guar:ray4}
.
\end{align}
\end{subequations}

From property \cref{eq:certinf:obj} of the certificate, $\bb^T \bar{\bz} + \bl^T \bar{\bmu} + \bu^T \bar{\bnu} < 0$. 
If $l_i \leq \til_i \leq \tiu_i \leq u_i, \forall i \in \iI$, then $(\bl - \tibl)^T \bar{\bmu} \geq 0$ and $(\bu - \tibu)^T \bar{\bnu} \geq 0$. In this case, the value \cref{eq:guar:ray4} is negative, so from \crefrange{eq:guar:ray1}{eq:guar:ray4}, the certificate $\K^*$ cut $\bar{\bz}^T (\bb - \bA \bx) \geq 0$ is violated. 
Therefore, the certificate $\K^*$ cut from the infeasible subproblem \ref{mod:conic} guarantees infeasibility of any LP OA \hyperref[mod:poly]{$\poly (\cZ, \tibl, \tibu)$} in the subtree of the node with integer variable bounds $\bl, \bu$.

More importantly for \cref{alg:oa}, the certificate $\K^*$ cut is likely to remain violated at `nearby' nodes outside of this subtree, as the conditions \crefrange{eq:guar:ray1}{eq:guar:ray4} have a natural interpretation from global sensitivity analysis. 
Perturbing the integer variable bounds from $\bl, \bu$ to $\tibl, \tibu$ changes the upper bound on $\bar{\bz}^T (\bb - \bA \bx)$ through a linear dependence on the values $\bmu \leq 0$ and $\bnu \geq 0$ of the dual variables in the improving ray of \ref{mod:conic*}.

\subsubsection{Certificate Cuts From Dual Optimal Solutions}
\label{sec:guar:exact:sol}

Suppose $(\hat{\bx}, (\hat{\bz}, \hat{\bmu}, \hat{\bnu}))$ is a complementary solution pair for the subproblem \ref{mod:conic}, certifying optimality of the solution pair.
Using the strong duality conditions (property \cref{eq:certsd} and feasibility for \ref{mod:conic} and \ref{mod:conic*}), any point $\bx \in \bbR^N$ satisfying the integer variable bounds $\til_i \leq x_i \leq \tiu_i, \forall i \in \iI$ at the new node and the certificate $\K^*$ cut $\hat{\bz}^T (\bb - \bA \bx) \geq 0$ has objective value:
\begin{subequations}
\begin{align}
\bc^T \bx 
\label{eq:guar:sol1}
& = -(\bA^T \hat{\bz} + \hat{\bmu}' + \hat{\bnu}')^T \bx 
\\
& = -\hat{\bz}^T \bA \bx - \bx^T (\hat{\bmu}' + \hat{\bnu}')
\\
& = -\bb^T \hat{\bz} + \hat{\bz}^T (\bb - \bA \bx)^T - \bx^T (\hat{\bmu}' + \hat{\bnu}') 
\\
& \geq -\bb^T \hat{\bz} - \bx^T (\hat{\bmu}' + \hat{\bnu}')
\\
& \geq -\bb^T \hat{\bz} - \bx^T (\hat{\bmu}' + \hat{\bnu}') - \sum_{\mathclap{i \in \iI}} \, ( (\til_i - x_i) \hat{\mu}_i + (\tiu_i - x_i) \hat{\nu}_i )
\label{eq:guar:sol2}
\\
& = -\bb^T \hat{\bz} - \tibl^T \hat{\bmu} - \tibu^T \hat{\bnu}
\label{eq:guar:sol3}
\\
& = (-\bb^T \hat{\bz} - \bl^T \hat{\bmu} - \bu^T \hat{\bnu}) + (\bl - \tibl)^T \hat{\bmu} + (\bu - \tibu)^T \hat{\bnu}
\\
& = \bc^T \hat{\bx} + (\bl - \tibl)^T \hat{\bmu} + (\bu - \tibu)^T \hat{\bnu}
\label{eq:guar:sol4}
.
\end{align}
\end{subequations}

If $l_i \leq \til_i \leq \tiu_i \leq u_i, \forall i \in \iI$, then $(\bl - \tibl)^T \hat{\bmu} \geq 0$ and $(\bu - \tibu)^T \hat{\bnu} \geq 0$. 
In this case, the value \cref{eq:guar:sol4} is no smaller than $\bc^T \hat{\bx}$, the lower bound from the subproblem \ref{mod:conic}. 
Therefore, the certificate $\K^*$ cut from the feasible subproblem \ref{mod:conic} guarantees that the optimal value of any LP OA \hyperref[mod:poly]{$\poly (\cZ, \tibl, \tibu)$} in the subtree of the node with integer variable bounds $\bl, \bu$ does not decrease, but may actually improve.%
\footnote{If \cref{alg:oa} branches on \cref{ln:brnch} after solving a bounded and feasible conic subproblem to get the tightest lower bound, then when examining a child node, this objective guarantee ensures the node's lower bound $L$ does not decrease when we update it to the optimal value of the LP OA on \cref{ln:lplb}.}

More importantly for \cref{alg:oa}, at `nearby' nodes outside of this subtree, the objective bounds implied by the certificate $\K^*$ cut in the LP OA model are likely to remain fairly tight. 
Perturbing the integer variable bounds from $\bl, \bu$ to $\tibl, \tibu$ changes the lower bound on $\bc^T \bx$ through a linear dependence on the values $\bmu \leq 0$ and $\bnu \geq 0$ of the dual variables in the complementary solution pair for \ref{mod:conic}.

\subsection{Under An LP Solver With A Feasibility Tolerance}
\label{sec:guar:tol}

So far, we have been assuming that the LP solver computes a solution that satisfies all the $\K^*$ cuts in the LP OAs exactly. In practice, LP solvers based on the Simplex method (except those that use rational arithmetic) enforce constraints up to an absolute constraint-wise violation tolerance $\delta > 0$ (typically set by the user). 
Therefore, a more realistic assumption is that any solution returned by the LP solver does not violate any $\K^*$ cut by more than $\delta$, i.e. a $\K^*$ point $\bz$ effectively yields a `relaxed $\K^*$ cut' $\bz^T (\bb - \bA \bx) \geq -\delta$. 
Under this relaxed condition, we may lose the `within-subtree' guarantees described in \cref{sec:guar:exact}.
However, noting that any positive scaling of a $\K^*$ point is still a $\K^*$ point, we demonstrate how to recover the infeasibility guarantee from \cref{sec:guar:exact:ray} exactly, and the objective bound guarantee from \cref{sec:guar:exact:sol} to within a given relative objective gap tolerance. Such an analysis appears to be novel in the MI-convex literature.

\subsubsection{Certificate Cuts From Dual Improving Rays}
\label{sec:guar:tol:ray}

Suppose $(\bar{\bz}, \bar{\bmu}, \bar{\bnu})$ is an improving ray of the dual subproblem \ref{mod:conic*}.
From the property \cref{eq:certinf:obj} of the certificate and the conditions \crefrange{eq:guar:ray1}{eq:guar:ray4}, any point $\bx \in \bbR^N$ satisfying the integer variable bounds $l_i \leq x_i \leq u_i, \forall i \in \iI$ and the relaxed certificate $\K^*$ cut condition $\bar{\bz}^T (\bb - \bA \bx) \geq -\delta$ must satisfy:
\begin{equation}
0 > \bb^T \bar{\bz} + \bl^T \bar{\bmu} + \bu^T \bar{\bnu} \geq \bar{\bz}^T (\bb - \bA \bx) \geq -\delta
.
\end{equation}
Therefore, if $\delta > 0$ is sufficiently large, the relaxed certificate $\K^*$ cut condition fails to enforce the infeasibility guarantee from \cref{sec:guar:exact:ray}. 

However, for a positive multiplier $\bar{\gamma} > 0$ satisfying:
\begin{equation}
\bar{\gamma} > \frac{\delta}{-\bb^T \bar{\bz} - \bl^T \bar{\bmu} - \bu^T \bar{\bnu}} > 0
\label{eq:tol:gamray}
,
\end{equation}
we have $\bar{\gamma} (\bb^T \bar{\bz} + \bl^T \bar{\bmu} + \bu^T \bar{\bnu}) < -\delta$. 
Therefore, the relaxed \textit{scaled} certificate $\K^*$ cut condition $\bar{\gamma} \bar{\bz}^T (\bb - \bA \bx) \geq -\delta$ recovers the infeasibility guarantee within the subtree of the node from which the certificate is obtained. 
Note that the scaling factor \cref{eq:tol:gamray} depends only on $\delta$, problem data, and the certificate for the infeasible subproblem \ref{mod:conic}. We can modify \cref{alg:oa} on \cref{ln:cutinf} to add the scaled $\K^*$ point $\bar{\gamma} \bar{\bz}$ to $\cZ$.

\subsubsection{Certificate Cuts From Dual Optimal Solutions}
\label{sec:guar:tol:sol}

Suppose $(\hat{\bx}, (\hat{\bz}, \hat{\bmu}, \hat{\bnu}))$ is a complementary solution pair for the subproblem \ref{mod:conic}.
From the conditions \crefrange{eq:guar:sol1}{eq:guar:sol4}, any point $\bx \in \bbR^N$ satisfying the integer variable bounds $l_i \leq x_i \leq u_i, \forall i \in \iI$ and the relaxed certificate $\K^*$ cut condition $\hat{\bz}^T (\bb - \bA \bx) \geq -\delta$ has objective value:
\begin{equation}
\bc^T \bx \geq -\bb^T \hat{\bz} + \hat{\bz}^T (\bb - \bA \bx) - \bl^T \hat{\bmu} - \bu^T \hat{\bnu} \geq L - \delta
\label{eq:tol:sd}
.
\end{equation}
Recall $L = \bc^T \hat{\bx} = -\bb^T \hat{\bz} - \bl^T \hat{\bmu} - \bu^T \hat{\bnu}$ is the optimal objective value of \ref{mod:conic} and \ref{mod:conic*}.
Therefore, the relaxed certificate $\K^*$ cut condition only enforces the objective guarantee from \cref{sec:guar:exact:sol} to an absolute tolerance of $\delta$. 
In general, it makes little sense for an objective guarantee to depend on the the LP solver's feasibility tolerance.

Instead, for a relative optimality gap tolerance $\epsilon > 0$, we can easily motivate a relative objective gap condition such as:
\begin{equation}
\frac{L - \bc^T \bx}{\lvert L \rvert + \theta} \leq \epsilon
\label{eq:tol:gap}
.
\end{equation}
Consider a positive multiplier $\hat{\gamma} > 0$ satisfying:
\begin{equation}
\hat{\gamma} \geq \frac{\delta}{\epsilon (\lvert L \rvert + \theta)} > 0
\label{eq:tol:gamsol}
.
\end{equation}
Modifying the conditions \cref{eq:tol:sd} for the relaxed \textit{scaled} certificate $\K^*$ cut condition $\hat{\gamma} \hat{\bz}^T (\bb - \bA \bx) \geq -\delta$, we get $\bc^T \bx \geq L - \sfrac{\delta}{\hat{\gamma}}$. Rearranging, this implies:
\begin{equation}
\frac{L - \bc^T \bx}{\lvert L \rvert + \theta} \leq \frac{\delta}{\hat{\gamma} (\lvert L \rvert + \theta)} \leq \epsilon
,
\end{equation}
so by scaling the certificate $\K^*$ cut by $\hat{\gamma}$, we achieve the relative objective gap guarantee \cref{eq:tol:gap} within the subtree of the node from which the certificate is obtained.
Note that the scaling factor \cref{eq:tol:gamsol} depends only on $\epsilon$, $\delta$, problem data, and the certificate for the bounded and feasible subproblem \ref{mod:conic}. We can modify \cref{alg:oa} on \cref{ln:cutsd} to add the scaled $\K^*$ point $\hat{\gamma} \hat{\bz}$ to $\cZ$.

\section{Tightening Polyhedral Relaxations}
\label{sec:tight}

In \cref{sec:tight:extr}, we outline a two-stage procedure for \textit{disaggregating} $\K^*$ cuts to get stronger polyhedral relaxations, and show how to maintain the certificate $\K^*$ cut guarantees from \cref{sec:guar}.
In \cref{sec:tight:init}, we argue for initializing the polyhedral relaxations using \textit{initial fixed $\K^*$ cuts}, and in \cref{sec:tight:sep}, we describe a procedure for cheaply obtaining \textit{separation $\K^*$ cuts} to cut off an infeasible LP OA solution.
All of our proposed techniques for tightening the LP OAs require minimal modifications to \cref{alg:oa} and are practical to implement.

\subsection{Extreme Ray Disaggregation}
\label{sec:tight:extr}

Consider a set of $\K^*$ points $\cZ = \{\bz^1, \ldots, \bz^J\} \subset \K^*$. 
By \textit{aggregating} the corresponding $\K^*$ cuts, we see they imply infinitely many $\K^*$ cuts:
\begin{equation}
\bz^T (\bb - \bA \bx) \geq 0 \qquad \forall \bz \in \cone(\cZ) 
,
\end{equation}
where $\cone(\cZ)$ is the \textit{conic hull} of $\cZ$, i.e. the set of conic (nonnegative) combinations of $\bz^1, \ldots, \bz^J$:
\begin{equation}
\cone(\cZ) = \{ \alpha^1 \bz^1 + \cdots \alpha^J \bz^J : \alpha^1, \ldots, \alpha^j \geq 0 \} \subset \K^*
.
\end{equation}
Thus for a redundant $\K^*$ point $\bz^{J+1} \in \cone(\cZ)$, the polyhedral relaxation of the conic constraint $\bb - \bA \bx \in \K$ implied by $\cZ \cup \{\bz^{J+1}\}$ is no stronger than that implied by $\cZ$ alone.
An \textit{extreme ray} of $\K^*$ is a point $\bz \in \K^*$ that cannot be written as a nontrivial conic combination of other points in $\K^*$ that are not positive rescalings of $\bz$. 
To maximize the efficiency of our polyhedral relaxations, we propose adding only extreme rays of $\K^*$ to the $\K^*$ point set $\cZ$ maintained by \cref{alg:oa}. 

Recall from \cref{sec:intro:conic} that our closed convex cone $\K$ is encoded as a Cartesian product $\K = \K_1 \times \cdots \times \K_K$ of standard primitive cones $\K_1, \ldots, \K_K$ (e.g. nonnegative, second-order, exponential, and positive semidefinite cones). 
A primitive closed convex cone cannot be written as a Cartesian product of two or more lower-dimensional closed convex cones \citep{Friberg2016}.
If $\K$ is separable, then its dual cone $\K^*$ is also separable:
\begin{equation}
\K^* = (\K_1 \times \cdots \times \K_K)^* = \K^*_1 \times \cdots \times \K^*_K 
\label{eq:sepcone}
.
\end{equation}
We exploit this separability and our understanding of the structure of the standard primitive cones to disaggregate a $\K^*$ point $\bz$ into extreme rays of $\K^*$.

First, we note that $\bz = (\tilde{\bz}^1, \ldots, \tilde{\bz}^K) \in \K^*$, where $\tilde{\bz}^k \in \K^*_k, \forall k \in \iK$. 
Second, for each $k \in \iK$, we disaggregate $\tilde{\bz}^k$ into extreme rays of the primitive standard dual cone $\K^*_k$. 
This step is trivial for linear cones. For second-order, positive semidefinite, and exponential cones, we describe practical computational procedures for dual disaggregation in \cref{sec:spec}.%
\footnote{For example, if $\K_k$ is a positive semidefinite cone, we disaggregate $\tilde{\bz}^k \in \K^*_k$ by performing an eigendecomposition on it; see \cref{sec:spec:psd:extr}.}
We have $\tilde{\bz}^k = \sum_{j \in \iJk} \tilde{\bz}^{k,j}$, where $\tilde{\bz}^{k,j} \neq 0$ is an extreme ray of $\K^*_k$, for all $j \in \iJk$. 
We choose these extreme rays so that none is a positive scaling of another, and $J_k$ does not exceed $\dim(\K^*_k)$. Note that $J_k = 0$ if $\tilde{\bz}^k = \bm{0}$. 

For some $k \in \iK$ and $j \in \iJk$, consider the point $\bz^{k,j} = (0, \ldots, 0, \tilde{\bz}^{k,j}, 0, \ldots, 0)$, which is nonzero only on the elements corresponding to the $k$th primitive dual cone. 
Since any cone contains the origin $\bm{0}$, and $\tilde{\bz}^{k,j} \in \K_k^*$, $\bz^{k,j} \in \K^*$ by equation \cref{eq:sepcone}. 
Furthermore, since $\tilde{\bz}^{k,j}$ is an extreme ray of $\K_k^*$, it cannot be written as a nontrivial sum of extreme rays of $\K^*_k$, and so $\bz^{k,j}$ cannot be written as a nontrivial sum of extreme rays of $\K^*$. Thus $\bz^{k,j}$ is an extreme ray of $\K^*$.

Our two-stage disaggregation procedure for $\bz \in \K^*$ yields $\sum_{k \in \iK} J_k \leq \dim(\K) = M$ extreme rays of $\K^*$:
\begin{equation}
\bz = \sum_{\mathclap{k \in \iK}} \;\;\;\, \sum_{\mathclap{j \in \iJk}} \bz^{k,j}
\label{eq:pracagg}
.
\end{equation} 
Besides adding potentially multiple $\K^*$ points to $\cZ$, no modifications are needed to the description of \cref{alg:oa}.
Since $\bz$ is clearly contained in the conic hull of these $\K^*$ points, there is no loss of strength in the polyhedral relaxations, so the certificate $\K^*$ guarantees from \cref{sec:guar:exact} are maintained.
The polyhedral relaxations are potentially much tighter, improving the power of the LP OA for fathoming a node by infeasibility or objective bound without proceeding to an expensive conic subproblem solve.%
\footnote{The LP solver may need to deal with more cuts at nodes visited early in the search tree, but is ultimately likely to need to examine fewer nodes overall and solve fewer expensive conic subproblems, so the tradeoff can be worthwhile.}

We can also recover the guarantees from \cref{sec:guar:tol} for an LP solver with a feasibility tolerance $\delta > 0$. 
We assume $\bz$ is a certificate $\K^*$ point that has already been scaled according to \cref{sec:guar:tol}. After disaggregating $\bz$, we scale each extreme ray up by $J = \sum_{k \in \iK} J_k$ before adding it to $\cZ$. 
The $J$ relaxed scaled disaggregated $\K^*$ cut conditions are:
\begin{equation}
(J \bz^{k,j})^T (\bb - \bA \bx) \geq -\delta \qquad \forall k \in \iK, j \in \iJk
.
\end{equation}
Summing and using equation \cref{eq:pracagg}, and dividing by $J$, we see that these conditions imply the relaxed scaled original $\K^*$ cut condition $\bz^T (\bb - \bA \bx) \geq -\delta$.

\subsection{Initial Fixed Polyhedral Relaxations}
\label{sec:tight:init}

We can modify \cref{alg:oa} on \cref{ln:zinit} to initialize a nonempty set $\cZ$ of \textit{initial fixed} $\K^*$ extreme rays that are not derived from subproblem certificates, but depend only on the geometry of $\K^*$. 
If $\K$ is a separable product of standard primitive cones, we can obtain initial fixed $\K^*$ extreme rays by treating each primitive cone constraint separately.
In particular, a linear cone constraint need not be relaxed at all, since it is equivalent to one $\K^*$ cut (for a nonnegative or nonpositive cone) or two $\K^*$ cuts (for the zero cone). 
In \cref{sec:spec}, we describe simple sets of initial fixed $\K^*$ extreme rays for second-order, positive semidefinite, or exponential primitive cones.%
\footnote{For example, for a positive semidefinite cone, we use the extreme rays of the polyhedral cone of diagonally dominant symmetric matrices as initial fixed $\K^*$ extreme rays; see \cref{sec:spec:psd:init}.} 
We show in \cref{sec:spec} how knowledge of the initial fixed $\K^*$ extreme rays allows us to tailor our extreme ray disaggregation procedures from \cref{sec:tight:extr} for certificate $\K^*$ points to further increase the strength of the polyhedral relaxations and reduce redundancy in $\cZ$.%
\footnote{However, to be able to recover the guarantees from \cref{sec:guar:tol} under an LP solver with a feasibility tolerance, we would need the ability to dynamically scale up the initial fixed $\K^*$ points.}

\subsection{Separation Of Infeasible Points}
\label{sec:tight:sep}

Inspired by separation-based OA algorithms, we can modify \cref{alg:oa} on \cref{ln:lpbrnch0} to add \textit{separation $\K^*$ points} to $\cZ$ that cut off a fractional optimal LP solution $\hat{\bx}$ that violates the conic constraint, right before branching on $\hat{\bx}$. 
We show that a separation $\K^*$ point exists when $\bb - \bA \hat{\bx} \notin \K$. 
Since $\K$ is closed and convex, there exists a hyperplane $(\bz, \theta)$ that separates $\hat{\by} = \bb - \bA \hat{\bx}$ from $\K$, i.e. $\bz^T \hat{\by} < \theta$ and $\bz^T \by \geq \theta, \forall \by \in \K$. 
Since the problem $\inf_{\by \in \K} \bz^T \by$ is homogeneous (as $\K$ is a cone) and the optimal value is bounded below by finite $\theta$, the optimal value must equal zero. 
So $\theta \leq 0$, implying $\bz^T \hat{\by} < 0$ and $\bz^T \by \geq 0, \forall \by \in \K$. 
Thus $\bz \in \K^*$ (by definition \cref{eq:kd} of $\K^*$), and it implies a $\K^*$ cut that separates $\hat{\bx}$ from the feasible set of the conic constraint.

A separation $\K^*$ point may fail to improve the objective lower bound from the LP OA, and does not in general possess the sort of guarantees from \cref{sec:guar} that a certificate $\K^*$ point implies. 
However, deriving a separation $\K^*$ point can be much cheaper than solving a continuous conic subproblem.
If $\K$ is a separable product of standard primitive cones, we can obtain separation $\K^*$ extreme rays easily by treating each primitive cone constraint separately. 
In \cref{sec:spec}, we describe practical computational methods for obtaining separation $\K^*$ extreme rays for primitive conic constraints involving second-order, positive semidefinite, or exponential cones.%
\footnote{For example, we obtain separation $\K^*$ extreme rays for a point that violates a positive semidefinite cone constraint by performing an eigendecomposition on it; see \cref{sec:spec:psd:sep}.}

\section{Pajarito Solver And Related Software}
\label{sec:soft}

We describe the software architecture and algorithmic implementation of Pajarito, our open source MI-convex solver.%
\footnote{In Pajarito's readme file (\href{https://github.com/JuliaOpt/Pajarito.jl}{github.com/JuliaOpt/Pajarito.jl}) we provide more guidance on the recommended ways of using the solver, as well as default options and tolerances.}
This section may be of particular interest to advanced users and developers of mathematical optimization software. 
We emphasize that our implementations diverge from the idealized description of \cref{alg:oa} in \cref{sec:alg:alg}, because of our decision to leverage powerful external mixed-integer linear (MILP/MIP) solvers through limited, solver-independent interfaces.
Developers of MI-conic software with low-level control of the MIP search tree are able to implement features of \cref{alg:oa} that we are not capable of in Pajarito.

\subsection{Julia And MathProgBase}
\label{sec:soft:mpb}

Pajarito is the first MI-convex solver written in the relatively young Julia language \citep{Bezanson2017}. 
MI-convex solvers such as $\alpha$-ECP, Artelys Knitro, Bonmin, DICOPT, FilMINT, MINLP\_BB, and SBB, which are reviewed by \citet{Bonami2012}, are to our knowledge written in C, C++, or Fortran. 
Julia is a high-level programming language that can match the performance of these lower-level languages for writing solvers with much less boilerplate code \citep{Lubin2015a}. 
Pajarito's compact codebase is thoroughly commented, and conveniently reusable and extensible by other researchers. 
We implement an extensive testing infrastructure with hundreds of unit tests. 
Since Pajarito's first release, several other MINLP solvers have been written in Julia and are available through MathProgBase, such as POD \citep{Nagarajan2017}, Juniper \citep{Kroger2018}, and Katana.%
\footnote{See \href{https://github.com/lanl-ansi/POD.jl}{github.com/lanl-ansi/POD.jl}, \href{https://github.com/lanl-ansi/Juniper.jl}{github.com/lanl-ansi/Juniper.jl}, and \href{https://github.com/lanl-ansi/Katana.jl}{github.com/lanl-ansi/Katana.jl}.}

Pajarito is integrated with the powerful MathProgBase abstraction layer. MathProgBase is a standardized API in Julia for interacting with optimization solvers, designed in part to allow the user to write solver-independent code.%
\footnote{MathProgBase is being replaced by a redesigned API, MathOptInterface. The process of building Pajarito has motivated many of the planned improvements in MathOptInterface.}
The breadth of problem classes covered by MathProgBase is described at \href{http://www.juliaopt.org/}{juliaopt.org} and distinguishes it from similar abstraction layers such as OSI \citep{OSI2004}, a COIN-OR library in C++. 
It includes specifications for continuous and mixed-integer solvers that use linear/quadratic, conic, or oracle-based NLP (nonlinear programming) forms. 

\begin{figure}[!htb]
\begin{tikzpicture}[
font = \small,
node distance = 0.8cm,
auto,
align = center,
mynode/.style = {text centered, draw=black, inner sep=0.15cm},
solver/.style = {mynode, rectangle, rounded corners, fill=red!15},
input/.style = {mynode, rectangle, fill=blue!15},
inter/.style = {text centered, fill=green!15, inner sep=1.3mm, anchor=center},
line/.style = {thick, <->, >={Stealth[inset=0pt,length=5pt]}},
]
\node (inp) [input] {\textbf{MI-convex model:}\\CBF, Convex.jl, CVXPY, JuMP};
\node (paj) [solver, below=of inp, yshift=-0.6cm] {\textbf{MI-conic solver:}\\Pajarito};
\node (con) [solver, below left=of paj, yshift=-1cm] {\textbf{Continuous solver:}\\CSDP, ECOS,\\MOSEK, SCS, SDPA};
\node (mip) [solver, below right=of paj, yshift=-1cm] {\textbf{MILP solver:}\\CBC, CPLEX, GLPK,\\Gurobi, MOSEK, SCIP};
\draw [line] (inp) -- node (ip) [inter, yshift=0.05cm] {conic interface} (paj);
\draw [line] (paj) -- node (pc) [inter, xshift=-0.7cm] {conic interface} (con);
\draw [line] (paj) -- node (pm) [inter, xshift=0.2cm] {linear/quadratic interface\\(through JuMP)} (mip);
\end{tikzpicture}
\caption{Pajarito's integration with MathProgBase.}
\label{fig:soft:arch}
\end{figure}
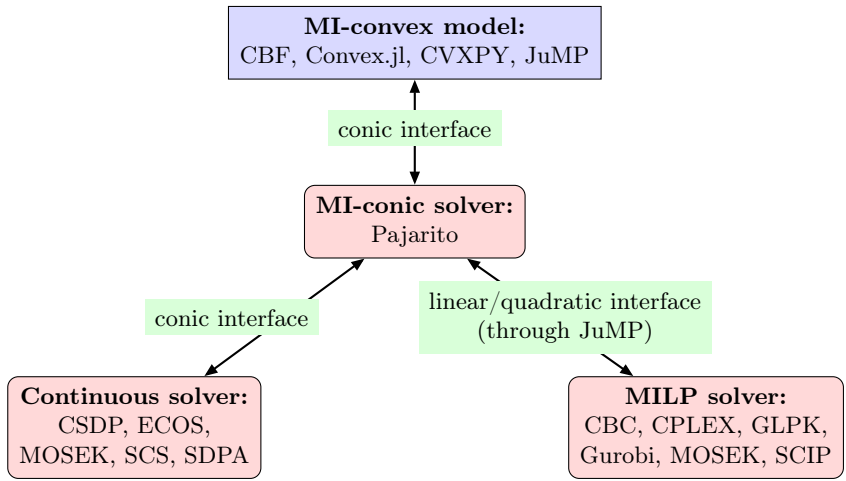

In \cref{sec:soft:access}, we describe accessing Pajarito through MathProgBase's conic interface (see top of \cref{fig:soft:arch}). 
The user specifies external MIP and continuous primal-dual conic solvers (including solver options) from the available solvers (i.e., those accessible though MathProgBase) and passes each solver object as an option into a function that creates a Pajarito solver object.
In \cref{sec:soft:alg}, we summarize Pajarito's main algorithmic implementations. 
Pajarito uses the modeling package JuMP to conveniently build and manage the external MIP solver's OA model. 
JuMP itself interacts with the MIP solver via MathProgBase's linear/quadratic interface (see bottom right of \cref{fig:soft:arch}).
To solve a continuous conic subproblem for a conic certificate, Pajarito calls the external primal-dual conic solver through the conic interface (see bottom left of \cref{fig:soft:arch}).%
\footnote{MathProgBase documents the conic and linear/quadratic interfaces at \href{http://mathprogbasejl.readthedocs.io/en/latest/}{mathprogbasejl.readthedocs.io/en/latest}. JuMP is documented at \href{http://www.juliaopt.org/JuMP.jl/0.18/}{juliaopt.org/JuMP.jl/0.18/}.}

\subsection{Accessing Pajarito}
\label{sec:soft:access}

Pajarito's use of conic form is a significant architectural difference from most existing MI-convex solvers, which interact with a MI-convex instance almost exclusively through oracles to query values and derivatives of the constraint and objective functions.
MathProgBase conic form can be described compactly from a constraint matrix in sparse or dense format, right-hand side and objective coefficient vectors, variable and constraint cones expressed as lists of standard primitive cones ($1$-dimensional vector sets) with corresponding ordered row indices, and a vector of variable types (each continuous, binary, or general integer).
In addition to the basic linear cones (nonnegative, nonpositive, zero, and free cones), Pajarito recognizes three standard primitive nonpolyhedral cones introduced in \cref{sec:intro:conic}: exponential cones (see \cref{sec:spec:exp}), second-order cones (see \cref{sec:spec:soc}), and positive semidefinite cones (see \cref{sec:spec:psd}).%
\footnote{As we note in \cref{sec:spec:soc}, Pajarito also recognizes rotated second-order cones, but for simplicity converts them to second-order cones during preprocessing.} 

\citet{Friberg2016} designed the Conic Benchmark Format (CBF) as a file format originally to support mixed-integer second-order cone (SOCP) and positive semidefinite cone (SDP) instances. 
In collaboration with Henrik Friberg, we extended the format to support exponential cones in Version 2, and developed a Julia interface ConicBenchmarkUtilities.jl to provide utilities for translating between CBF and MathProgBase conic format.%
\footnote{Pajarito's extensive unit tests rely on small example instances loaded from CBF files.}
One may use Pajarito to solve any instance in the Conic Benchmark Library (CBLIB), which contains thousands of benchmark problems from a wide variety of sources. 

\citet{Lubin2016} demonstrate that all 333 known MI-convex instances in MINLPLib2 \citep{GAMS2018} are representable with linear, second-order, exponential, and \textit{power} cones. 
Since a power cone constraint is representable with linear and exponential cone constraints, Pajarito can be used to solve any of the MI-convex instances in MINLPLib2.
We translated 115 instances from the MINLPLIB2 library to CBF and contributed them to CBLIB.%
\footnote{\citet{Lubin2016} first translated these instances from the MINLPLIB2 library into Convex.jl models. We used ConicBenchmarkUtilities.jl to translate these to CBF. The instances, available at \href{http://github.com/mlubin/MICPExperiments}{github.com/mlubin/MICPExperiments}, are 48 `rsyn' instances, 48 `syn' instances, 6 `tls' instances, 12 `clay' instances, and the challenging `gams01' instance.} 
Many of the MINLPLIB2 instances have tiny values artificially-introduced in order to work around potential numerical issues with smooth derivative-based NLP solvers \citep{Guenluek2012}, which we manually removed before converting to conic form.
For example, the instance `rsyn0805h' has a constraint:
\begin{equation}
\left( \frac{x_{289}}{10^{-6} + b_{306}} - \frac{6}{5} \log \left( 1 + \frac{x_{285}}{10^{-6} + b_{306}} \right) \right) (10^{-6} + b_{306}) \leq 0
\label{con:rsyn}
,
\end{equation}
where $b_{306}, x_{285}, x_{289}$ are scalar variables. 
Without the artificial $10^{-6}$ values, a conic encoding of the NLP constraint \cref{con:rsyn} in terms of the exponential cone $\KE$ is:
\begin{equation}
(b_{306} + x_{285}, b_{306}, \sfrac{5}{6} \, x_{289}) \in \KE
.
\end{equation}

Within Julia, the modeling packages JuMP \citep{Dunning2017} and Convex.jl \citep{Udell2014} each provide a convenient way for users to specify MI-convex problems, call Pajarito solver, and interpret solutions. 
JuMP is particularly useful for a large, sparse problem involving complex indexing schemes for variables, expressions, or constraints. 
It efficiently builds a MathProgBase conic form representation of a problem involving second-order or positive semidefinite cones, but currently does not recognize exponential cones.

Convex.jl, unlike JuMP, is a Disciplined Convex Programming (DCP) modeling package. It defines a list of \textit{atoms} for the user to model a MI-convex problem with and performs automatic verification of convexity of the continuous relaxation by applying simple composition rules described by \citet{Grant2006}.
It converts the problem into a MI-conic instance in MathProgBase conic form through epigraph and perspective transformations that introduce additional variables and constraints in conic form using only the standard primitive cones recognized by Pajarito.
CVXPY \citep{Diamond2016} is a Python-based DCP modeling package analogous to Convex.jl. In collaboration with Steven Diamond and Baris Ungun, we developed cmpb.jl (\href{http://github.com/mlubin/cmpb}{github.com/mlubin/cmpb}), a prototype C API to MathProgBase that enables Pajarito to be called on a problem modeled with CVXPY. 

We illustrate Convex.jl and JuMP modeling using a simple MI-convex example described by \citet[ch.~7.5]{Boyd2004}: `E-optimal experimental design'.%
\footnote{More Pajarito examples can be found at \href{http://github.com/JuliaOpt/Pajarito.jl/tree/master/examples}{github.com/JuliaOpt/Pajarito.jl/blob/master/examples}.}
While we can solve E-optimal experimental design exactly using Pajarito, \citet[ch.~7.5]{Boyd2004} choose to relax the integrality constraints in order to use a continuous convex solver before rounding the fractional solution heuristically.
To begin, we set up the Pajarito solver object \mintinline{julia}{mysolver} using a GLPK MILP solver object and a SCS conic solver object, each with internal options set.
\begin{minted}[breaklines,autogobble,fontsize=\small]{julia}
using Pajarito, GLPKMathProgInterface, SCS #load packages
mysolver = PajaritoSolver(log_level = 3, #use verbose output
    mip_solver = GLPKSolverMIP(msg_lev = GLPK.MSG_OFF), #set MIP solver 
    cont_solver = SCSSolver(eps = 1e-6, verbose = 0)) #set conic solver
\end{minted}
Pajarito performs a sanity check on the combination of options and solvers specified.%
\footnote{MathProgBase does not attempt to provide an abstraction for solver parameters like convergence tolerances. In cases where we need certain tolerances on the continuous conic and MIP solvers in order for Pajarito to converge to a requested tolerance, it is the user's responsibility to set these tolerances. For example, we ask users to manually adjust the MIP solver's linear feasibility tolerance and integer feasibility tolerance for improved convergence behavior. These cases are documented in Pajarito's readme file.}
Next, we model and solve the problem using Convex.jl as follows, where $p, m, n \in \bbR$ and $V \in \bbR^{n \times p}$ are problem data.
\begin{minted}[breaklines,autogobble,fontsize=\small]{julia}
using Convex
mp = Variable(p, Positive(), :Int) #create p nonneg. integer variables
eOpt = maximize(lambdamin(V * diagm(mp./m) * V'), #max. min. eigenvalue
    sum(mp) <= m) #add linear constraint
solve!(eOpt, mysolver) #solve model using Pajarito solver
@show eOpt.status, eOpt.optval, mp.value #show solve status and results
\end{minted}
Alternatively, we model and solve the problem using JuMP as follows.
\begin{minted}[breaklines,autogobble,fontsize=\small]{julia}
using JuMP
eOpt = Model(solver = mysolver) #initialize model using Pajarito solver
@variable(eOpt, mp[1:p] >= 0, Int) #create p nonneg. integer variables
@constraint(eOpt, sum(mp) <= m) #add linear constraint
@variable(eOpt, t) #create auxiliary variable
FI = V * diagm(mp./m) * V' #create linear expression matrix
@SDconstraint(eOpt, FI - t * eye(n) >= 0) #add PSD constraint on matrix
@objective(eOpt, Max, t) #maximize linear objective
@show solve(eOpt) #solve model and show status
@show getobjectivevalue(eOpt), getvalue(mp) #show solve results
\end{minted}
Pajarito manipulates the conic data and performs sanity checks. We refer to the resulting preprocessed representation of the instance as \ref{mod:micp}. 
After Pajarito executes one of the OA algorithms described in \cref{sec:soft:alg} on \ref{mod:micp}, the user can use Convex.jl or JuMP to conveniently query information such as Pajarito's solve status, objective bound, objective value, and solution.

\subsection{Basic Algorithmic Implementations}
\label{sec:soft:alg}

We discuss the main conic-certificate-based methods Pajarito uses to solve the preprocessed MI-conic model \ref{mod:micp}. 
We omit many options, enhancements, and numerical details that can be understood from the Pajarito readme file and from browsing the high-level Julia code and comments.%
\footnote{For explaining computational experiments, \cref{sec:exp:algs} briefly introduces several other algorithmic variants that we do not discuss here, such as separation-based methods that do not utilize conic certificates.}
In \cref{sec:soft:init}, we summarize the initialization procedure for the OA model, an MILP relaxation of \ref{mod:micp} that Pajarito constructs and later refines (with extreme ray $\K^*$ cuts) using JuMP.
In \cref{sec:soft:iter}, we describe the `iterative' method, an extension of the simple sequential OA algorithm by \citet{Lubin2017}. 
In \cref{sec:soft:msd}, we describe the `MIP-solver-driven' (MSD) method, so-called because it relies on the power of the branch-and-cut MIP solver to manage convergence in a single tree. 
Since MathProgBase's solver-independent abstraction for MIP solver callbacks is designed primarily around shared behavior between CPLEX and Gurobi, Pajarito is limited to interacting with the MIP solver through a \textit{lazy cut callback} function and a \textit{heuristic callback} function.
Although the MSD method is generally much faster than the iterative method, the latter may be used with MILP solvers for which callback functionality is unavailable or unreliable.

\subsubsection{Initializing The MIP OA Model}
\label{sec:soft:init}

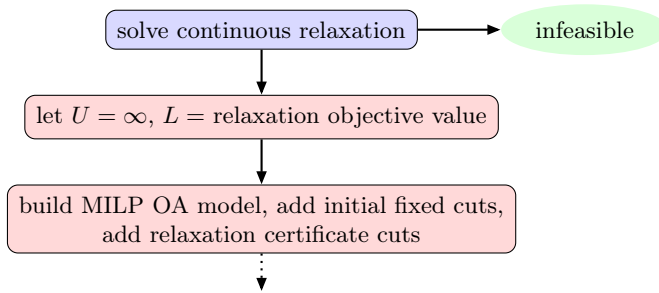
\begin{figure}[!htb]
\begin{tikzpicture}[
font = \small,
node distance = 0.6cm,
auto,
align = center,
mynode/.style = {text centered, draw=black, inner sep=0.15cm},
process/.style = {mynode, rectangle, rounded corners, fill=red!15},
solve/.style = {mynode, rectangle, rounded corners, fill=blue!15},
status/.style = {text centered, ellipse, fill=green!15, inner sep=1.3mm, anchor=center},
arrow/.style = {thick, ->, >={Stealth[inset=0pt,length=5pt]}},
]
\node (rel) [solve] {solve continuous relaxation};
\node (infs) [status, right=of rel, xshift=0.5cm] {infeasible};
\node (bnd) [process, below=of rel] {let $U = \infty$, $L =$ relaxation objective value};
\node (ini) [process, below=of bnd] {build MILP OA model, add initial fixed cuts,\\add relaxation certificate cuts};
\draw [arrow] (rel) -- (bnd);
\draw [arrow] (bnd) -- (ini);
\draw [arrow] (rel) -- (infs);
\draw [arrow, dotted] (ini.south) to node[auto] {} ++ (0,-5mm);
\end{tikzpicture}
\caption{Pajarito's OA model initialization.}
\label{fig:soft:init}
\end{figure}

We first solve the continuous relaxation of \ref{mod:micp} (in which only the integrality constraints are relaxed), using the primal-dual conic solver (see top of \cref{fig:soft:init}).%
\footnote{The conic relaxation is analogous to the first node subproblem in \cref{alg:oa}, but without finite integer variable bounds. We preprocess this conic model slightly to tighten any non-integral bounds on the integer variables.}
If the conic solver indicates this relaxation is infeasible, then \ref{mod:micp} must be infeasible, so we terminate with an `infeasible' status.
If the conic solver returns a complementary solution pair, the optimal value gives an objective lower bound $L > -\infty$ for \ref{mod:micp}. 
Otherwise, we set $L = -\infty$.
We initialize the objective upper bound $U$ for \ref{mod:micp} to $\infty$.

Using JuMP, we build the initial OA model, adding the variables and integrality constraints and setting the objective (see bottom of \cref{fig:soft:init}). 
We then add initial fixed cuts for each primitive cone, as we describe in \cref{sec:tight:init}. 
Primitive linear cone constraints are imposed entirely (as equivalent LP equality or inequality constraints), and for each primitive nonpolyhedral cone, we add a small number of initial fixed cuts (defined in \cref{sec:spec}).

A complementary solution pair from the conic relaxation solve yields a $\K^*$ point, so we perform an extreme ray disaggregation from \cref{sec:tight:extr} and add certificate cuts for each primitive nonpolyhedral cone (using the procedures in \cref{sec:spec}).
These continuous relaxation certificate cuts technically guarantee that the root node of the OA model has an optimal value no smaller than $L$.%
\footnote{This can be seen from a simple modification of the complementary solution case polyhedral relaxation guarantee we prove in \cref{sec:guar:exact:sol}, with trivial integer variable bounds.}
This is important because we cannot handle unboundedness of the OA model.

\subsubsection{Iterative Method}
\label{sec:soft:iter}

\begin{figure}[!htb]
\begin{tikzpicture}[
font = \small,
node distance = 0.6cm,
auto,
align = center,
mynode/.style = {text centered, draw=black, inner sep=0.15cm},
process/.style = {mynode, rectangle, rounded corners, fill=red!15},
solve/.style = {mynode, rectangle, rounded corners, fill=blue!15},
status/.style = {text centered, ellipse, fill=green!15, inner sep=1.3mm, anchor=center},
arrow/.style = {thick, ->, >={Stealth[inset=0pt, length=5pt]}},
]
\node (mip) [solve] {solve OA model};
\node (conv) [process, below=of mip] {let $L =$ objective bound;\\$L, U$ converged?};
\node (opts) [status, right=of conv, xshift=0.5cm] {optimal};
\node (fails) [status, above=of opts, yshift=-0.4cm] {OA fail};
\node (infs) [status, above=of fails, yshift=-0.4cm] {infeasible};
\node (int) [process, below=of conv] {integer sub-solution repeated?};
\node (sub) [solve, below left=of int, xshift=2.9cm, yshift=-0.2cm] {solve subproblem; feasible?};
\node (subup) [process, below=of sub] {update $U$ and incumbent}; 
\node (subcut) [process, below=of subup] {add certificate cuts}; 
\node (sep) [process, below right=of int, xshift=-1.4cm, yshift=-0.2cm] {solution conic feasible?};
\node (sepup) [process, below=of sep] {update $U$ and incumbent}; 
\node (sepcut) [process, below=of sepup] {add separation cuts}; 
\node (warm) [process, below=of int, yshift=-3.5cm] {warm-start OA model with incumbent}; 
\draw [arrow, <-, dotted] (mip.north) to node [auto] {} ++ (0, 5mm);
\draw [arrow] (mip) -- (infs);
\draw [arrow] (mip) -- (fails);
\draw [arrow] (mip) -- (conv);
\draw [arrow] (conv) -- node [anchor=north] {yes} (opts);
\draw [arrow] (conv) -- node [anchor=west] {no} (int);
\draw [arrow] (int) -- node [anchor=north west, near start] {no} (sub);
\draw [arrow] (int) -- node [anchor=north east, near start] {yes} (sep);
\draw [arrow] (sub.south east) to [out=300, in=35] node [anchor=east, pos=0.135] {no} (subcut.east);
\draw [arrow] (sub) -- node [anchor=west] {yes} (subup);
\draw [arrow] (subup) -- (subcut);
\draw [arrow] (subcut) -- (warm);
\draw [arrow] (sep.south west) to [out=220, in=155] node [anchor=west, pos=0.16] {no} (sepcut.west);
\draw [arrow] (sep) -- node [anchor=east] {yes} (sepup);
\draw [arrow] (sepup.north east) to [out=80, in=290] node [anchor=west, very near start] {} (opts.south east);
\draw [arrow] (sepcut) -- (warm);
\draw [arrow] (warm.west) to [out=150, in=190] (mip.west);
\end{tikzpicture}
\caption{Pajarito's iterative method, following initialization.}
\label{fig:soft:iter}
\end{figure}
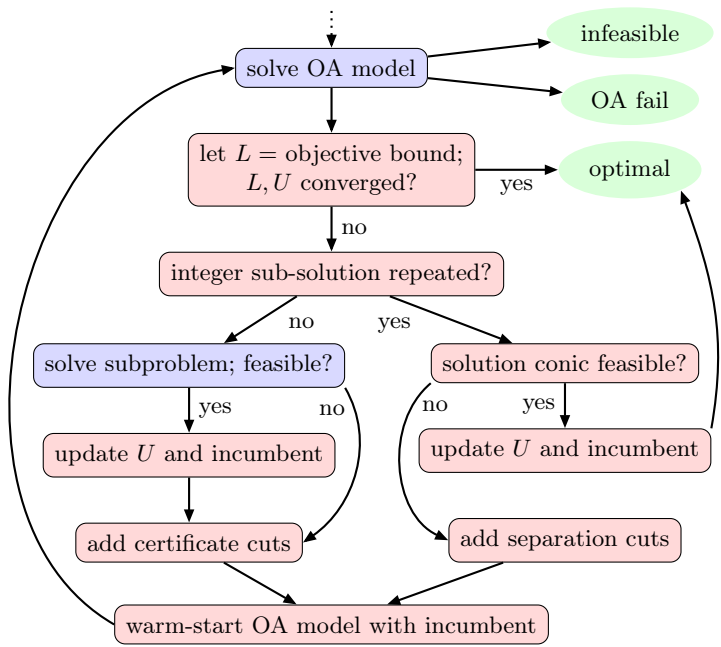

The iterative method, following initialization in \cref{fig:soft:init}, is outlined in \cref{fig:soft:iter}. At each iteration of the main loop, Pajarito solves the current OA model using the MIP solver.%
\footnote{We suggest the user set the MIP solver's relative optimality gap tolerance to its smallest possible value.}
If the OA model is infeasible, \ref{mod:micp} must be infeasible, so we terminate with an `infeasible' status. If it is unbounded, Pajarito fails with an `OA fail' error status, as we are unable to handle unbounded rays.
If the MIP solver returns an optimal solution to the OA model, this OA solution satisfies the integrality constraints and initial fixed cuts, but in general not all of nonpolyhedral primitive cone constraints. 
The MIP solver's objective bound provides a lower bound for \ref{mod:micp}, so we update $L$. Pajarito terminates with an `optimal' status if the relative optimality gap condition \cref{eq:tol:gap} on $L, U$ is satisfied.%
\footnote{Pajarito uses $\theta = 10^{-5}$ (to avoid division by $0$). The gap tolerance $\epsilon > 0$ is specified by the user, but defaults to $10^{-5}$.}

If after solving the OA model we have an optimal OA solution and the objective bounds haven't converged, we check whether the OA sub-solution on the integer variables has been encountered before. 
If so, we check the conic feasibility of the OA solution. We calculate the absolute violation on each primitive nonpolyhedral cone constraint as the violation of the appropriate separation cut (defined in \cref{sec:spec}). 
If the worst absolute violation does not exceed Pajarito's feasibility tolerance (set by the user), then the OA solution is considered feasible.
In this case, since the solution is optimal for the OA model, we can consider it optimal for \ref{mod:micp}, so we update the incumbent and upper bound and terminate the solve immediately. 
If the OA solution is not considered feasible, we add all of the separation cuts that are (significantly) violated to the OA model. 

If the integer sub-solution has not already been encountered at a previous iteration, then we solve a continuous conic subproblem in which the integer variables are fixed to their values in the integer sub-solution. 
This subproblem is analogous to \ref{mod:conic} from \cref{sec:alg:sub}, with $\bl = \bu$.%
\footnote{In preprocessing, we remove any subproblem equality constraints that effectively have no variables when an integer sub-solution is fixed. For efficient loading of the subproblem data at each iteration, we only change the constant vector $\bb$ of the preprocessed conic subproblem, as this is the only data that changes.}
Since it is more constrained than the OA model, the conic subproblem is bounded or infeasible.
If the conic subproblem solver fails to return a certificate, we backtrack and perform the separation procedure (as if the integer sub-solution repeated). 
Otherwise, we scale the certificate's dual solution or dual ray according to \cref{sec:guar:tol} (using the tolerance values set as Pajarito options), then disaggregate the scaled $\K^*$ point and add extreme ray certificate cuts to the OA model (as we described for the continuous relaxation certificate in \cref{sec:soft:init}).
In the case of a complementary solution pair, the primal solution yields a feasible point for \ref{mod:micp}, since it satisfies both the integrality and conic constraints.  
If it has an objective value better than $U$, we update $U$ and the incumbent solution and check the relative optimality gap condition again.%
\footnote{Conic solvers typically do not use an absolute primitive constraint-wise feasibility tolerance, as Pajarito does for checking feasibility of OA solutions for the conic constraint. Our incumbent may not satisfy this notion of feasibility, since we do not perform a feasibility check on the conic solver's primal subproblem solutions.}

After adding separation or certificate cuts, we warm-start the MIP solver with our incumbent and re-execute the main loop. 
The procedure in \cref{fig:soft:iter} is iterated until $L$ and $U$ converge or the MIP solver detects infeasibility.%
\footnote{If the user sets a time limit, Pajarito may terminate with the status `user limit'. Pajarito sets the time limit on each MIP or conic solve to the remaining time. Note that the vast majority of Pajarito execution time is spent in MIP or conic solves.}
Note that since we only add cuts to the OA model on every loop, if the first OA model is bounded, then all subsequent (refined) OA models are bounded or infeasible, and the sequence of lower bounds $L$ is nondecreasing.%
\footnote{In case of failures of strong duality at some conic subproblems, Pajarito may fail to converge, as there exists no finite set of cuts that can tighten the lower bound sufficiently to meet the upper bound. See \citet{Lubin2016} for a discussion of strong duality in OA.}

\subsubsection{MIP-Solver-Driven Method}
\label{sec:soft:msd}

\begin{figure}[!htb]
\begin{tikzpicture}[
font = \small,
node distance = 0.6cm,
auto,
align = center,
mynode/.style = {text centered, draw=black, inner sep=0.15cm},
process/.style = {mynode, rectangle, rounded corners, fill=red!15},
solve/.style = {mynode, rectangle, rounded corners, fill=blue!15},
status/.style = {text centered, ellipse, fill=green!15, inner sep=1.3mm, anchor=center},
arrow/.style = {thick, ->, >={Stealth[inset=0pt, length=5pt]}},
]
\node (mip) [solve] {call branch-and-cut\\solver on OA model};
\node (fails) [status, right=of mip, xshift=0.5cm] {OA fail};
\node (infs) [status, above=of fails, yshift=-0.4cm] {infeasible};
\node (opts) [status, below=of fails, yshift=0.4cm] {optimal};
\node (int) [process, below=of mip, yshift=-0.6cm] {integer sub-solution repeated?};
\node (sub) [solve, below left=of int, xshift=2.5cm, yshift=-0.2cm] {solve subproblem; feasible?};
\node (subcut) [process, below=of sub, xshift=1.2cm, yshift=-0.1cm] {add certificate cuts}; 
\node (subup) [process, left=of subcut, xshift=-0.1cm] {store solution}; 
\node (sep) [process, below right=of int, xshift=-1.8cm, yshift=-0.2cm] {solution conic feasible?};
\node (sepcut) [process, below=of sep, yshift=-0.1cm] {add separation cuts}; 
\node (heur) [process, anchor=west] at (subup.west |- mip) {add stored solution};
\begin{pgfonlayer}{background}
\node (heurcb) [rectangle, fill=gray!10, inner sep=0.2cm, rounded corners=3mm, fit=(heur)] {};
\node (lazycb) [rectangle, fill=gray!10, inner sep=0.2cm, rounded corners=3mm, fit=(int) (sub) (subup) (subcut) (sep) (sepcut)] {};
\end{pgfonlayer}
\node[above, inner sep=1pt] at (heurcb.north) {heuristic callback\vphantom{y}};
\node[above, inner sep=1pt] at (lazycb.north -| heurcb) {lazy callback};
\draw [arrow, <-, dotted] (mip.north) to node[auto] {} ++ (0,5mm);
\draw [arrow] (mip.east) -- (infs);
\draw [arrow] (mip.east) -- (fails);
\draw [arrow] (mip.east) -- (opts);
\draw [thick, dashed] (mip) -- node [anchor=east] {} (heur);
\draw [arrow, dashed] (mip) -- node [anchor=west, pos=0.4] {integral\\solution} (int);
\draw [arrow] (int) -- node [anchor=north west, near start] {no} (sub);
\draw [arrow] (int) -- node [anchor=north east, near start] {yes} (sep);
\draw [arrow] (sub) -- node [anchor=north east, near start] {no} (subcut);
\draw [arrow] (sub) -- node [anchor=north west, near start] {yes} (subup);
\draw [arrow] (subup) -- (subcut);
\draw [arrow] (sep) -- node [anchor=west] {no} (sepcut);
\draw [arrow, dotted] (subup.north west) to [out=105, in=255] node [anchor=east] {} (heur.south west);
\end{tikzpicture}
\caption{Pajarito's MIP-solver-driven (MSD) method, following initialization.}
\label{fig:soft:msd}
\end{figure}
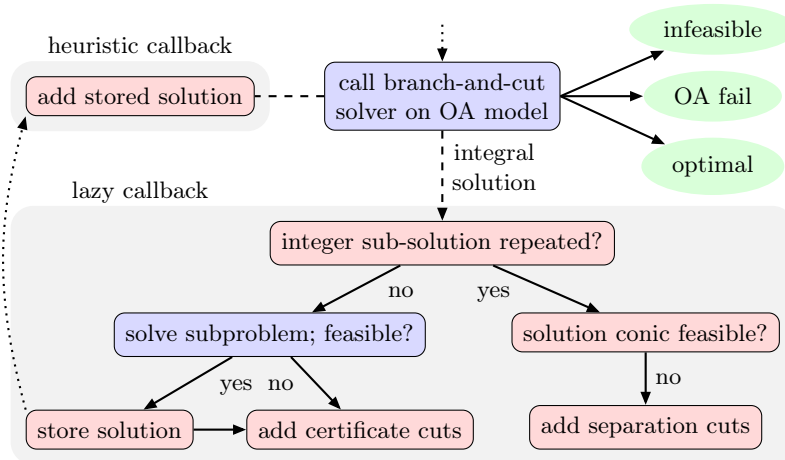

The MSD method, following initialization in \cref{fig:soft:init}, is outlined in \cref{fig:soft:msd}. 
As in the iterative method, Pajarito returns an `OA fail' status if the MIP solver detects unboundedness (as we are unable to handle unbounded rays), or an `infeasible' status in the case of infeasibility.%
\footnote{If the user sets a time limit, Pajarito sets a time limit on the MIP solver, and terminates with a `user limit' status if this limit is reached.}
The MIP-solver-independent callback interface allows us to pass in lazy cuts in during a lazy callback and feasible solutions during a heuristic callback, however we cannot exert any control over branching decisions, node selection, fathoming, or node lower bound updating. 

The MIP solver calls the lazy callback function whenever it finds an integer-feasible OA solution at a node. 
During a lazy callback, we first check whether the integral OA solution from the MIP solver is repeated. If so, we derive separation cuts to add as lazy constraints; if none can be added, the MIP solver considers the solution feasible and may update its incumbent. 
If the integer sub-solution is repeated, we solve a new (bounded or infeasible) conic subproblem. Since we lack the ability to query the node's integer variable bounds, we only solve subproblems with fixed integer sub-solutions, as in the iterative method. 
If the conic solver returns a certificate, we scale and disaggregate the $\K^*$ point (as we described for the iterative method in \cref{sec:soft:iter}), and add extreme ray cuts as lazy constraints.%
\footnote{The MIP solver is not guaranteed to respect the cuts that we add, and we may need to re-add the same cuts during multiple lazy callbacks (unlike in the iterative method, where cuts previously added are respected). We actually store a dictionary from the integer sub-solution to the cuts. For each repeated integer sub-solution, we re-add these saved certificate cuts, in addition to the new separation cuts.}
In the case of a complementary solution pair, the primal solution yields a feasible point for \ref{mod:micp}, which we store.
During a heuristic callback, if there is a stored feasible solution to \ref{mod:micp} that has never been added as a heuristic solution, we add it. 

Since there are no guarantees on when or how frequently the MIP solver calls the heuristic callback function, we may not be able to indirectly update the MIP solver's incumbent and upper bound when we are able to.
Partly for this reason, Pajarito maintains its own upper bound and incumbent (not illustrated in \cref{fig:soft:msd}, which we update during lazy callbacks. 
During each lazy callback, we ask the MIP solver for its lower bound and check our relative optimality gap condition (as we described for the iterative method). If the condition is met, we force the MIP solver to terminate early. 
In this case, or if the MIP solver terminates with an optimal solution and we verify that the relative optimality gap condition is met, we return our incumbent solution with an `optimal' status.%
\footnote{Note the user is responsible for setting the desired relative optimality gap tolerance on both the MIP solver and on Pajarito directly.}

\subsection{Some Advanced Algorithmic Enhancements}
\label{sec:soft:exten}

We conclude with several key optional OA enhancements we implemented in Pajarito.
First, Pajarito by default uses an \textit{extended formulation} for each second-order cone constraint. \citet{Vielma2016} demonstrate on a testset of mixed-integer second-order cone (MISOCP) problems that OA algorithms tend to converge much faster when using this extended representation for each second-order cone constraint.%
\footnote{DCP modeling software implementations such as Convex.jl do not perform this transformation because they are simply designed to access conic solvers. Pajarito keeps the original second-order cone formulation in the conic subproblems because conic solvers are likely to perform better with this representation than with the higher-dimensional extended formulation.}
In \cref{sec:socef}, we describe how to lift a $\K^*$ cut for the second-order cone into $\K^*$ cuts for the extended formulation. This technique also allows us to describe a much more economical set of initial fixed $\K^*$ cuts for the second-order cone.

Second, Pajarito can optionally use a MISOCP OA model instead of a MILP OA model. 
There exist several powerful MISOCP solvers that can be used, or Pajarito itself may be used. 
Since a second-order cone constraint can imply an infinite number of $\K^*$ cuts, Pajarito can achieve tighter relaxations of the conic constraint in the OA model. 
Of course, this can only make practical sense for certain types of problems that aren't pure MISOCP.
One potential use case is where the MI-conic problem has second-order cones as well as exponential and/or positive semidefinite cones, but the only stable and efficient continuous conic solvers we have access to are for SOCP problems. 
In this case, we can use Pajarito with the SOCP solver as an MISOCP solver, and pass this into a second Pajarito solver that uses a conic solver for mixed-cone problems. 
This arrangement helps minimize the number of calls to the less-effective conic solver. 
Another use case is for problems involving positive semidefinite cone constraints. 
We demonstrate in \cref{sec:psdsoc} how to strengthen $\K^*$ cuts for PSD constraints to rotated-second-order cone constraints.%
\footnote{For the MSD method, since most MISOCP solvers don't currently allow adding lazy quadratic constraints, only the initial fixed cuts can be strengthened in this way.}

\section{Computational Experiments}
\label{sec:exp}

Our computational experiments demonstrate the speed and robustness of our open source MI-conic solver Pajarito. 
As we emphasize in \cref{sec:soft:alg}, our algorithmic implementations differ from the description of \cref{alg:oa}, because of our practical decision to use branch-and-cut MILP solvers through a limited, solver-independent interface.
In \cref{sec:exp:pres}, we summarize our metrics for comparing the practical performance of different MI-conic solvers and describe our presentation of tables and performance profile plots. 
In \cref{sec:exp:solv}, we benchmark Pajarito (version 0.5.1) and several open source and commercial mixed-integer second-order cone (MISOCP) solver packages accessible through MathProgBase on a MISOCP library, and conclude that Pajarito is the fastest and most-reliable open source solver for MISOCP. 
In \cref{sec:exp:algs}, we compare the performance of several of Pajarito's algorithmic variants on MI-conic instances involving mixtures of positive semidefinite, second-order, and exponential cones, demonstrating practical advantages of the methodological extensions we describe in \cref{sec:guar,sec:tight} and \cref{sec:spec}. 
The scripts and data we use to run our experiments are available in the supplement \href{http://github.com/chriscoey/PajaritoSupplement}{github.com/chriscoey/PajaritoSupplement}.

\subsection{Presentation Of Results}
\label{sec:exp:pres}

We define a `solver' as a MathProgBase solver object given a particular complete set of algorithmic options.%
\footnote{Each solver we test is deterministic, i.e. it performs consistently across different runs on a particular dedicated system.}
We define an `instance' as a particular MI-conic problem (stored in CBF format; see \cref{sec:soft:access}) that is known to be feasible and bounded (but an optimal solution or the optimal objective value is not necessarily known). 
For a particular instance, a solver may return a `solution', which is a vector of real floating point numbers representing an assignment of the variables of the instance (not necessarily feasible for the constraints). 

First, we compare the performances of a group of MI-conic solvers on a particular testset of instances by counting the number of instances for which each solver returns and apparently proves `approximate optimality' of a solution. 
To be more precise, we use the following four categories to characterize a solver's apparent success or failure on an instance.
\begin{description}
\item[ex] (exclude) means either the solver incorrectly claims the instance is infeasible or unbounded, or the solver returns a solution it claims is approximately-optimal but we detect one of the following inconsistencies.
\begin{itemize}
\item The solution significantly violates at least one primitive cone constraint or integrality constraint.%
\footnote{Absolute violation of a primitive cone constraint is calculated as worst violation of the inequalities defining the standard cone (see \cref{sec:spec}), and our tolerances are $10^{-6}$ for linear cones, $10^{-5}$ for second-order and exponential cones, and $10^{-4}$ for positive semidefinite cones. Variable-wise integrality violation is calculated as distance to the nearest integer, and our tolerance is $10^{-6}$.}
\item The relative objective gap condition (equation \cref{eq:tol:gap}) for optimality is significantly violated.%
\footnote{Our optimality condition \cref{eq:tol:gap} matches that used by most MIP solvers. We set the constant $\theta = 10^{-5}$ (to avoid division by zero) and use the tolerance $\epsilon = 10^{-5}$. We ensure we do not exclude in the case that the gap we calculate is sensitive to a solver's different value of $\theta$.}
\item The objective value or objective bound significantly differs from that of a preponderance of other solvers.%
\footnote{This is assessed semi-manually from output of our scripts.}
\end{itemize}
\item[co] (converge) means the solver returns a solution that it claims is (approximately) optimal (and it is not excluded for the reasons above).
\item[li] (reach limit) means the solver does not terminate before the time limit, or (rarely) the solver reaches a memory limit and is forced to terminate.
\item[er] (error) means the solver crashes or terminates with an error message.
\end{description}

Second, we compare aggregate quantitative measures of solver performance.
We define the shifted geometric mean $\tilde{g}$ of $L$ positive values $p_1, \ldots, p_L$ as:
\begin{equation}
\tilde{g}(\bm{p}, q) = \prod_{\mathclap{l \in \iL}} \, (p_l + q)^\frac{1}{L} - q
,
\end{equation}
where $q > 0$ is the shift \citep{Achterberg2009}. 
Unlike the standard geometric mean $\tilde{g}(\bm{p}, 0)$, the shifted geometric mean decreases the relative influence of smaller values in $\bm{p}$, thus giving less weight to very `easy' instances (small values are preferable for all of our metrics). 
We shift by $q = 10$ seconds for execution times, $q = 1$ iterations for iteration counts, and $q = 10$ nodes for MIP-solver-reported node counts. 
For comparing a particular group of solvers $S_1, \ldots, S_n$ on a particular performance metric (such as execution time), we calculate for each solver $S_i$ the following three shifted geomeans, each over a different subset of the testset.
\begin{description}
\item[aco] (all solvers converge) is calculated over the instances for which $S_1, \ldots, S_n$ all have a `co' status.
\item[tco] (this solver converges) is calculated over the instances for which $S_i$ has a `co' status.
\item[all] (all instances) is calculated over all instances. Missing execution times are set to the time limit, and missing iteration/node counts are ignored.
\end{description}

Finally, we employ `performance profiles', described by \citet{Dolan2002,Gould2016}, to visually compare the relative execution times and iteration or node counts of pairs of solvers.
Again, we decrease the relative influence of very easy instances by shifting the metrics by the same shift values $q$ we use for shifted geomeans.
A performance profile is a plot that should be interpreted as follows: for a fixed factor $F$ on the horizontal axis (a linear scale from $1$ to the value at the bottom right of the plot), the level of solver $S_i$ on the vertical axis (a linear scale from $0$ to $1$) represents the proportion $P_i$ of instances (out of the instances for which at least one of the pair of solvers has a `co' status) for which $S_i$ has a `co' status and a (shifted) performance metric that is within a factor of $F$ of per-instance best achieved by either solver. 
So, at $F = 1$ (i.e. on the left vertical axis), $P_i$ is the fraction of solved instances on which solver $S_i$ has the best performance. As $F$ increases, we can infer that solver $S_i$ has reliably better performance than solver $S_j$ if $P_i$ remains above $P_j$.

\subsection{MISOCP Solver Performance Comparisons}
\label{sec:exp:solv}

Our open source Pajarito solvers, `Iter-GLPK' and `Iter-CBC', use the iterative method (see \cref{sec:soft:iter}) with ECOS \citep{Domahidi2013} for continuous conic subproblems and CBC or GLPK for MILPs.%
\footnote{We do not test the Pajarito's MIP-solver-driven method with CBC or GLPK MIP solvers because their support for MathProgBase callbacks is limited.}
Our two restricted-license Pajarito solvers, `Iter-CPLEX' (using the iterative method) and `MSD-CPLEX' (using the MIP-solver-driven method; see \cref{sec:soft:msd}), call MOSEK's continuous conic solver and CPLEX's MILP solver. 

The open source Bonmin solver package is described in detail by \citet{Bonami2008} and uses CBC to manage branching and Ipopt to solve continuous NLP (derivative-based nonlinear programming) subproblems.%
\footnote{We are unaware of any mainstream open source solvers designed for MISOCP. The functional representation of the second-order cone has points of nondifferentiability that may cause Bonmin to crash or suffer numerical issues.}
Our `Bonmin-BB' solver uses the nonlinear B\&B method (no polyhedral approximation), `Bonmin-OA' uses the B\&B OA method, and `Bonmin-OA-D' is equivalent to the `Bonmin-OA' solver but applied to transformed instances that use the second-order cone extended formulation we describe in \cref{sec:soft:exten}.
Our two restricted-license MISOCP solvers are `SCIP' and `CPLEX'. Unlike Bonmin, these MISOCP solvers use the second-order cone extended formulation internally.%
\footnote{CPLEX is available under an academic or commercial licence, and SCIP is an academic solver that is not released under an OSI-approved open source license. We use CPLEX version 12.7.0 and SCIP version 4.0.0.}

These nine MISOCP solvers are each given a relative optimality gap tolerance of $10^{-5}$. 
The `SCIP' and `CPLEX' solvers are given an absolute linear-constraint-wise feasibility tolerance of $10^{-8}$, and `CPLEX' is given an integrality tolerance of $10^{-9}$.
The MILP solvers used by Pajarito are given an absolute linear-constraint-wise feasibility tolerance of $10^{-8}$, an integrality tolerance of $10^{-9}$, and a relative optimality gap tolerance of $0$ for `Iter-GLPK', `Iter-CBC', and `Iter-CPLEX' and $10^{-5}$ for `MSD-CPLEX'. 
Due to limited resources, we set a one hour time limit for each run of a solver on an instance, and run all solvers (including the MILP and conic solvers called by Pajarito) in single-threaded mode. 

We use a testset of $120$ MISOCP instances drawn from the larger CBLIB library, recently compiled by \citet{Friberg2016}. 
The testset contains randomly selected subsets of most of the major families of models in CBLIB. 
We exclude instances that are not bounded and feasible, or are solved in under $5$ seconds by all solvers, or are unable to be solved by all solvers in under an hour.
Our computations are performed on the Amazon EC2 cloud computing platform with `m4.xlarge' computing nodes having $16$GB of RAM.%
\footnote{See \href{http://aws.amazon.com/ec2/instance-types/}{aws.amazon.com/ec2/instance-types}.}
As the computing nodes are virtual machines, timing results on EC2 are subject to random variability, but repeated runs suggest the variation is sufficiently small to avoid impacting our conclusions.
The nodes run Ubuntu 16.04 with Julia version 0.6.0. 
Version information for the Julia packages can be obtained from the supplement.

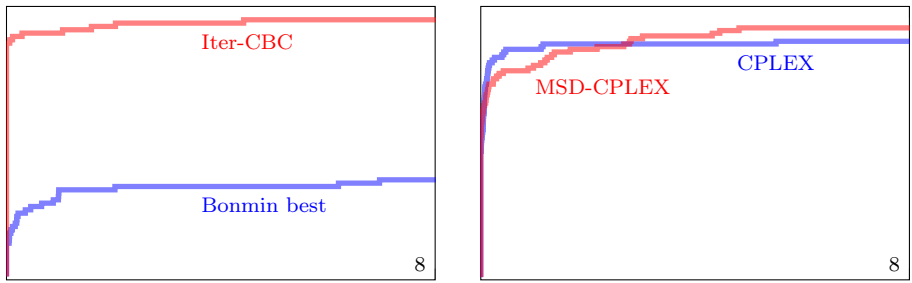
\begin{figure}[!htb]
\begin{subfigure}{0.48\textwidth}
\centering\begin{tikzpicture}[font = \footnotesize]
\begin{axis}[xmax=8, xtick={8}]
\addplot+ coordinates {
(1.0, 0.012195121951219513)
(1.0, 0.024390243902439025)
(1.0, 0.036585365853658534)
(1.0, 0.04878048780487805)
(1.0, 0.06097560975609756)
(1.0, 0.07317073170731707)
(1.0, 0.08536585365853659)
(1.0, 0.0975609756097561)
(1.0, 0.10975609756097561)
(1.0, 0.12195121951219512)
(1.0, 0.13414634146341464)
(1.0430685561336026, 0.14634146341463414)
(1.0472059475253042, 0.15853658536585366)
(1.0493362092141154, 0.17073170731707318)
(1.0718832873248636, 0.18292682926829268)
(1.096365762058061, 0.1951219512195122)
(1.1428984015156518, 0.2073170731707317)
(1.166355735832212, 0.21951219512195122)
(1.1677514934156539, 0.23170731707317074)
(1.1887708057843076, 0.24390243902439024)
(1.3073077319741815, 0.25609756097560976)
(1.3950871870245125, 0.2682926829268293)
(1.573439177302647, 0.2804878048780488)
(1.7641654667804547, 0.2926829268292683)
(1.8476806447795036, 0.3048780487804878)
(1.8505171590501808, 0.3170731707317073)
(1.8527637580791558, 0.32926829268292684)
(2.7707273874923715, 0.34146341463414637)
(6.4211353087368686, 0.35365853658536583)
(7.091632795658307, 0.36585365853658536)
(10.54526748270037, 0.3780487804878049)
(13.411921913188495, 0.3902439024390244)
(16.838748932768677, 0.4024390243902439)
(19.480552284240893, 0.4146341463414634)
(27.75561830651233, 0.4268292682926829)
(30.049659082720716, 0.43902439024390244)
(30.522663777479647, 0.45121951219512196)
(33.3038756414091, 0.4634146341463415)
(56.86728772094848, 0.47560975609756095)
(135.37232110363158, 0.4878048780487805)
(265.16667879716874, 0.5)
(530.3333575943375, 0.5121951219512195)
(530.3333575943375, 0.524390243902439)
(530.3333575943375, 0.5365853658536586)
(530.3333575943375, 0.5487804878048781)
(530.3333575943375, 0.5609756097560976)
(530.3333575943375, 0.573170731707317)
(530.3333575943375, 0.5853658536585366)
(530.3333575943375, 0.5975609756097561)
(530.3333575943375, 0.6097560975609756)
(530.3333575943375, 0.6219512195121951)
(530.3333575943375, 0.6341463414634146)
(530.3333575943375, 0.6463414634146342)
(530.3333575943375, 0.6585365853658537)
(530.3333575943375, 0.6707317073170732)
(530.3333575943375, 0.6829268292682927)
(530.3333575943375, 0.6951219512195121)
(530.3333575943375, 0.7073170731707317)
(530.3333575943375, 0.7195121951219512)
(530.3333575943375, 0.7317073170731707)
(530.3333575943375, 0.7439024390243902)
(530.3333575943375, 0.7560975609756098)
(530.3333575943375, 0.7682926829268293)
(530.3333575943375, 0.7804878048780488)
(530.3333575943375, 0.7926829268292683)
(530.3333575943375, 0.8048780487804879)
(530.3333575943375, 0.8170731707317073)
(530.3333575943375, 0.8292682926829268)
(530.3333575943375, 0.8414634146341463)
(530.3333575943375, 0.8536585365853658)
(530.3333575943375, 0.8658536585365854)
(530.3333575943375, 0.8780487804878049)
(530.3333575943375, 0.8902439024390244)
(530.3333575943375, 0.9024390243902439)
(530.3333575943375, 0.9146341463414634)
(530.3333575943375, 0.926829268292683)
(530.3333575943375, 0.9390243902439024)
(530.3333575943375, 0.9512195121951219)
(530.3333575943375, 0.9634146341463414)
(530.3333575943375, 0.975609756097561)
(530.3333575943375, 0.9878048780487805)
(530.3333575943375, 1.0)
(530.3333575943375, 1.0121951219512195)
} node[below right] at (4, 0.34146341) {Bonmin best};
\addplot+ coordinates {
(1.0, 0.012195121951219513)
(1.0, 0.024390243902439025)
(1.0, 0.036585365853658534)
(1.0, 0.04878048780487805)
(1.0, 0.06097560975609756)
(1.0, 0.07317073170731707)
(1.0, 0.08536585365853659)
(1.0, 0.0975609756097561)
(1.0, 0.10975609756097561)
(1.0, 0.12195121951219512)
(1.0, 0.13414634146341464)
(1.0, 0.14634146341463414)
(1.0, 0.15853658536585366)
(1.0, 0.17073170731707318)
(1.0, 0.18292682926829268)
(1.0, 0.1951219512195122)
(1.0, 0.2073170731707317)
(1.0, 0.21951219512195122)
(1.0, 0.23170731707317074)
(1.0, 0.24390243902439024)
(1.0, 0.25609756097560976)
(1.0, 0.2682926829268293)
(1.0, 0.2804878048780488)
(1.0, 0.2926829268292683)
(1.0, 0.3048780487804878)
(1.0, 0.3170731707317073)
(1.0, 0.32926829268292684)
(1.0, 0.34146341463414637)
(1.0, 0.35365853658536583)
(1.0, 0.36585365853658536)
(1.0, 0.3780487804878049)
(1.0, 0.3902439024390244)
(1.0, 0.4024390243902439)
(1.0, 0.4146341463414634)
(1.0, 0.4268292682926829)
(1.0, 0.43902439024390244)
(1.0, 0.45121951219512196)
(1.0, 0.4634146341463415)
(1.0, 0.47560975609756095)
(1.0, 0.4878048780487805)
(1.0, 0.5)
(1.0, 0.5121951219512195)
(1.0, 0.524390243902439)
(1.0, 0.5365853658536586)
(1.0, 0.5487804878048781)
(1.0, 0.5609756097560976)
(1.0, 0.573170731707317)
(1.0, 0.5853658536585366)
(1.0, 0.5975609756097561)
(1.0, 0.6097560975609756)
(1.0, 0.6219512195121951)
(1.0, 0.6341463414634146)
(1.0, 0.6463414634146342)
(1.0, 0.6585365853658537)
(1.0, 0.6707317073170732)
(1.0, 0.6829268292682927)
(1.0, 0.6951219512195121)
(1.0, 0.7073170731707317)
(1.0, 0.7195121951219512)
(1.0, 0.7317073170731707)
(1.0, 0.7439024390243902)
(1.0, 0.7560975609756098)
(1.0, 0.7682926829268293)
(1.0, 0.7804878048780488)
(1.0, 0.7926829268292683)
(1.0, 0.8048780487804879)
(1.0, 0.8170731707317073)
(1.0, 0.8292682926829268)
(1.0, 0.8414634146341463)
(1.0, 0.8536585365853658)
(1.0, 0.8658536585365854)
(1.0528809083736543, 0.8780487804878049)
(1.0984121687362984, 0.8902439024390244)
(1.2535714347142664, 0.9024390243902439)
(1.9154414551485588, 0.9146341463414634)
(2.386897975613603, 0.926829268292683)
(2.7793272902966972, 0.9390243902439024)
(4.872796852421928, 0.9512195121951219)
(530.3333575943375, 0.9634146341463414)
(530.3333575943375, 0.975609756097561)
(530.3333575943375, 0.9878048780487805)
(530.3333575943375, 1.0)
(530.3333575943375, 1.0121951219512195)
} node[below right] at (4, 0.9390243) {Iter-CBC};
\end{axis}
\end{tikzpicture}
\caption{Open source Bonmin (instance-wise best of $3$) and Pajarito iterative solvers.}
\label{fig:misocp:open}
\end{subfigure}%
\hspace{0.04\textwidth}%
\begin{subfigure}{0.48\textwidth}
\centering\begin{tikzpicture}[font = \footnotesize]
\begin{axis}[xmax=8, xtick={8}]
\addplot+ coordinates {
(1.0, 0.00980392156862745)
(1.0, 0.0196078431372549)
(1.0, 0.029411764705882353)
(1.0, 0.0392156862745098)
(1.0, 0.049019607843137254)
(1.0, 0.058823529411764705)
(1.0, 0.06862745098039216)
(1.0, 0.0784313725490196)
(1.0, 0.08823529411764706)
(1.0, 0.09803921568627451)
(1.0, 0.10784313725490197)
(1.0, 0.11764705882352941)
(1.0, 0.12745098039215685)
(1.0, 0.13725490196078433)
(1.0, 0.14705882352941177)
(1.0, 0.1568627450980392)
(1.0, 0.16666666666666666)
(1.0, 0.17647058823529413)
(1.0, 0.18627450980392157)
(1.0, 0.19607843137254902)
(1.0, 0.20588235294117646)
(1.0, 0.21568627450980393)
(1.0, 0.22549019607843138)
(1.0, 0.23529411764705882)
(1.0, 0.24509803921568626)
(1.0, 0.2549019607843137)
(1.0, 0.2647058823529412)
(1.0, 0.27450980392156865)
(1.0, 0.28431372549019607)
(1.0, 0.29411764705882354)
(1.0, 0.30392156862745096)
(1.0, 0.3137254901960784)
(1.0, 0.3235294117647059)
(1.0, 0.3333333333333333)
(1.0, 0.3431372549019608)
(1.0, 0.35294117647058826)
(1.0, 0.3627450980392157)
(1.0, 0.37254901960784315)
(1.0, 0.38235294117647056)
(1.0, 0.39215686274509803)
(1.0, 0.4019607843137255)
(1.0, 0.4117647058823529)
(1.0, 0.4215686274509804)
(1.0, 0.43137254901960786)
(1.0, 0.4411764705882353)
(1.0, 0.45098039215686275)
(1.0, 0.46078431372549017)
(1.0, 0.47058823529411764)
(1.0054195593758806, 0.4803921568627451)
(1.0063721903843719, 0.49019607843137253)
(1.0067018722107588, 0.5)
(1.0125170538283255, 0.5098039215686274)
(1.0138749389096744, 0.5196078431372549)
(1.0205709962460412, 0.5294117647058824)
(1.0248434386410081, 0.5392156862745098)
(1.0324277199648093, 0.5490196078431373)
(1.0344522448868743, 0.5588235294117647)
(1.035023152897776, 0.5686274509803921)
(1.0358150082987632, 0.5784313725490197)
(1.0429098317379137, 0.5882352941176471)
(1.0436488034307856, 0.5980392156862745)
(1.0459268924856082, 0.6078431372549019)
(1.052586719240808, 0.6176470588235294)
(1.0530731087224348, 0.6274509803921569)
(1.0535883778753254, 0.6372549019607843)
(1.0553122586765393, 0.6470588235294118)
(1.0554829827933068, 0.6568627450980392)
(1.0588366218015306, 0.6666666666666666)
(1.0608130616397111, 0.6764705882352942)
(1.0636645975143604, 0.6862745098039216)
(1.0652113295259356, 0.696078431372549)
(1.0682436155791384, 0.7058823529411765)
(1.0770827741125948, 0.7156862745098039)
(1.0919488545620897, 0.7254901960784313)
(1.0964124538243696, 0.7352941176470589)
(1.108577107244179, 0.7450980392156863)
(1.1099224596266763, 0.7549019607843137)
(1.115207612666624, 0.7647058823529411)
(1.12033135380324, 0.7745098039215687)
(1.1320916220506734, 0.7843137254901961)
(1.14252045206174, 0.7941176470588235)
(1.1863392295405868, 0.803921568627451)
(1.2238006598092874, 0.8137254901960784)
(1.3455607060301882, 0.8235294117647058)
(1.3748200665553088, 0.8333333333333334)
(1.4094164373893805, 0.8431372549019608)
(1.9860113428856183, 0.8529411764705882)
(2.016329439894827, 0.8627450980392157)
(5.822437674312484, 0.8725490196078431)
(36.47728043256067, 0.8823529411764706)
(72.95456086512134, 0.8921568627450981)
(72.95456086512134, 0.9019607843137255)
(72.95456086512134, 0.9117647058823529)
(72.95456086512134, 0.9215686274509803)
(72.95456086512134, 0.9313725490196079)
(72.95456086512134, 0.9411764705882353)
(72.95456086512134, 0.9509803921568627)
(72.95456086512134, 0.9607843137254902)
(72.95456086512134, 0.9705882352941176)
(72.95456086512134, 0.9803921568627451)
(72.95456086512134, 0.9901960784313726)
(72.95456086512134, 1.0)
(72.95456086512134, 1.0098039215686274)
} node[below right] at (5, 0.86274509) {CPLEX};
\addplot+ coordinates {
(1.0, 0.00980392156862745)
(1.0, 0.0196078431372549)
(1.0, 0.029411764705882353)
(1.0, 0.0392156862745098)
(1.0, 0.049019607843137254)
(1.0, 0.058823529411764705)
(1.0, 0.06862745098039216)
(1.0, 0.0784313725490196)
(1.0, 0.08823529411764706)
(1.0, 0.09803921568627451)
(1.0, 0.10784313725490197)
(1.0, 0.11764705882352941)
(1.0, 0.12745098039215685)
(1.0, 0.13725490196078433)
(1.0, 0.14705882352941177)
(1.0, 0.1568627450980392)
(1.0, 0.16666666666666666)
(1.0, 0.17647058823529413)
(1.0, 0.18627450980392157)
(1.0, 0.19607843137254902)
(1.0, 0.20588235294117646)
(1.0, 0.21568627450980393)
(1.0, 0.22549019607843138)
(1.0, 0.23529411764705882)
(1.0, 0.24509803921568626)
(1.0, 0.2549019607843137)
(1.0, 0.2647058823529412)
(1.0, 0.27450980392156865)
(1.0, 0.28431372549019607)
(1.0, 0.29411764705882354)
(1.0, 0.30392156862745096)
(1.0, 0.3137254901960784)
(1.0, 0.3235294117647059)
(1.0, 0.3333333333333333)
(1.0, 0.3431372549019608)
(1.0, 0.35294117647058826)
(1.0, 0.3627450980392157)
(1.0, 0.37254901960784315)
(1.0, 0.38235294117647056)
(1.0, 0.39215686274509803)
(1.0, 0.4019607843137255)
(1.0, 0.4117647058823529)
(1.0, 0.4215686274509804)
(1.0, 0.43137254901960786)
(1.0, 0.4411764705882353)
(1.0, 0.45098039215686275)
(1.0, 0.46078431372549017)
(1.0, 0.47058823529411764)
(1.0, 0.4803921568627451)
(1.0, 0.49019607843137253)
(1.0, 0.5)
(1.0, 0.5098039215686274)
(1.0, 0.5196078431372549)
(1.0, 0.5294117647058824)
(1.0011193930852056, 0.5392156862745098)
(1.0016761605235645, 0.5490196078431373)
(1.0041254686213292, 0.5588235294117647)
(1.0050593452607004, 0.5686274509803921)
(1.0134655298301956, 0.5784313725490197)
(1.017088510186184, 0.5882352941176471)
(1.0225313194195804, 0.5980392156862745)
(1.0257608397335773, 0.6078431372549019)
(1.035239179820909, 0.6176470588235294)
(1.0416283628572538, 0.6274509803921569)
(1.0464516704711262, 0.6372549019607843)
(1.0627453808357095, 0.6470588235294118)
(1.0809491837475353, 0.6568627450980392)
(1.0816685449186307, 0.6666666666666666)
(1.086291776533692, 0.6764705882352942)
(1.1058497238521305, 0.6862745098039216)
(1.1142119796120102, 0.696078431372549)
(1.121899158324991, 0.7058823529411765)
(1.1325793078110145, 0.7156862745098039)
(1.203845354901535, 0.7254901960784313)
(1.2063870901833145, 0.7352941176470589)
(1.2574700715421976, 0.7450980392156863)
(1.3005328723044858, 0.7549019607843137)
(1.367754905385992, 0.7647058823529411)
(1.7714241153931016, 0.7745098039215687)
(1.861241844937024, 0.7843137254901961)
(1.9778001283520517, 0.7941176470588235)
(2.0462350579719333, 0.803921568627451)
(2.0859179392254874, 0.8137254901960784)
(2.156951641099406, 0.8235294117647058)
(2.2009038395679372, 0.8333333333333334)
(2.461301047814494, 0.8431372549019608)
(2.899553897968147, 0.8529411764705882)
(3.334816658854219, 0.8627450980392157)
(3.429106845287982, 0.8725490196078431)
(3.4531282032708, 0.8823529411764706)
(3.641049795239174, 0.8921568627450981)
(4.779439183057116, 0.9019607843137255)
(4.902828717963059, 0.9117647058823529)
(5.216827540058391, 0.9215686274509803)
(9.977693860733334, 0.9313725490196079)
(13.48230784210268, 0.9411764705882353)
(20.459483014502023, 0.9509803921568627)
(72.95456086512134, 0.9607843137254902)
(72.95456086512134, 0.9705882352941176)
(72.95456086512134, 0.9803921568627451)
(72.95456086512134, 0.9901960784313726)
(72.95456086512134, 1.0)
(72.95456086512134, 1.0098039215686274)
} node[below right] at (1.7, 0.775) {MSD-CPLEX};
\end{axis}
\end{tikzpicture}
\caption{CPLEX MISOCP and Pajarito MSD solvers.}
\label{fig:misocp:comm}
\end{subfigure}
\caption{MISOCP solver execution time performance profiles.}
\end{figure}

\Cref{tab:misocp} summarizes the status counts and shifted geomeans of performance metrics on instance subsets (explained in \cref{sec:exp:pres}) for the nine MISOCP solvers on the $120$ MISOCP instances.
The Bonmin solvers fail on most instances, and overall solve significantly fewer instances than the open source Pajarito solvers. 
Pajarito tends perform faster using CBC rather than GLPK.%
\footnote{However, for most of the $9$ excluded instances from `Iter-CBC', we verify that CBC is responsible for the significant integrality violations that result in exclusion.}
\Cref{fig:misocp:open} is a performance profile (explained in \cref{sec:exp:pres}) comparing the execution times of the open source Pajarito (with CBC) solver and the instance-wise best of the three Bonmin solvers. 
From these results, we claim that Pajarito with ECOS and CBC is the fastest and most reliable open source MISOCP solver.

Using CPLEX, Pajarito's MSD method is significantly faster and more reliable than its iterative method.%
\footnote{For `MSD-CPLEX', the two errors occur where Pajarito claims a solution is suboptimal and has an objective gap no worse than $1.04 \times 10^{-5}$, and the one exclusion occurs where Pajarito's solution violates a linear constraint by $9.78 \times 10^{-6}$.}
The performance profile \cref{fig:misocp:comm} compares the execution times of Pajarito's MSD method using CPLEX's MILP solver against CPLEX's specialized MISOCP solver.
The execution time comparisons between `CPLEX' and `MSD-CPLEX' are ambiguous, however we argue that, at least by our metrics, Pajarito is a more reliable MISOCP solver. 

\begin{table}[!htb]
\footnotesize
\begin{tabular}{l l *{4}{S[table-format=2]} *{2}{S[table-format=2.1]} S[table-format=3]}
\toprule
& & \multicolumn{4}{l}{statuses} & \multicolumn{3}{l}{time (s)} 
\\
\cmidrule(lr){3-6} \cmidrule(lr){7-9}
& solver & {co} & {li} & {er} & {ex} & {aco} & {tco} & {all} \\
\midrule
\multirow{5}{*}{\rotatebox[origin=c]{90}{open source}}
& Bonmin-BB & 
34 & 44 & 11 & 31 & 38.0 & 83.8 & 463 \\
& Bonmin-OA & 
25 & 53 & 29 & 13 & 64.2 & 64.5 & 726 \\
& Bonmin-OA-D & 
30 & 48 & 29 & 13 & 15.1 & 61.6 & 610 \\
& Iter-GLPK & 
56 & 60 & 3 & 1 & 2.0 & 29.7 & 377 \\
& Iter-CBC & 
78 & 30 & 3 & 9 & 1.6 & 50.3 & 163 \\
\addlinespace
\multirow{4}{*}{\rotatebox[origin=c]{90}{restricted}}
& SCIP & 
74 & 35 & 8 & 3 & 3.2 & 41.5 & 160 \\
& CPLEX & 
90 & 16 & 5 & 9 & 0.9 & 16.1 & 50 \\
& Iter-CPLEX & 
86 & 26 & 0 & 8 & 0.4 & 37.0 & 106 \\
& MSD-CPLEX & 
97 & 20 & 2 & 1 & 0.4 & 18.2 & 56 \\
\bottomrule
\end{tabular}
\caption{MISOCP solver performance summary.}
\label{tab:misocp}
\end{table}

\subsection{Comparative Testing Of Algorithmic Variants}
\label{sec:exp:algs}

To compare the performance of several of Pajarito's algorithmic variants, we use a testset of $95$ MI-conic instances involving mixtures of positive semidefinite (PSD), second-order, and exponential cones. 
These instances are all bounded and feasible and come from the following four sources.%
\footnote{Formulations for instances we generated can be found at \href{https://github.com/JuliaOpt/Pajarito.jl/tree/master/examples}{github.com/JuliaOpt/Pajarito.jl/tree/master/examples}.}
\begin{description}
\item[Discrete experimental design] ($14$ instances). 
Recall from \cref{sec:soft:access} that \citet[Ch.~7.5]{Boyd2004} describes MI-convex experimental design problems. 
We generate `A-optimal' and `E-optimal' instances that include PSD cones, and `D-optimal' instances that include PSD and exponential cones.
\item[Portfolios with mixed risk constraints] ($16$ instances). 
We formulate a portfolio problem that maximizes expected returns subject to some combinatorial constraints on stocks and three types of convex risk constraints on subsets of stocks with known covariances. 
Each instance includes multiple exponential cones from entropy risk constraints, second-order cones from norm risk constraints, and PSD cones from robust norm risk constraints. 
\item[Retrofit-synthesis of process networks] ($32$ instances). 
We select a representative subset of the two CBLIB families `syn' and `rsyn'. 
Each instance includes exponential cones.
\item[A subset of the MISOCP testset] ($33$ instances). 
We select a representative subset of the CBLIB families `estein', `ccknapsack', `sssd', `uflquad', and `portfoliocard'.
\end{description}

We use Pajarito with Gurobi (version 7.5.2) as the MILP solver and MOSEK (version 9.0.0.29-alpha) as the continuous conic solver.%
\footnote{MOSEK 9 is the first version to recognize exponential cones.}
Pajarito is given a relative optimality gap tolerance of $10^{-5}$. 
Gurobi is given an absolute linear-constraint-wise feasibility tolerance of $10^{-8}$, an integrality tolerance of $10^{-9}$, and a relative optimality gap tolerance of $0$ when the iterative method is used and $10^{-5}$ when the MSD method is used. 
We set a one hour time limit for each run of a solver on an instance, and limit Gurobi and MOSEK to $8$ threads. 
We run the computations on dedicated hardware with $16$ Intel Xeon E5-2650 CPUs ($2$GHz) and $64$GB of RAM. 
Repeated runs suggest the variation is sufficiently small to avoid impacting our conclusions. 
The machine runs Ubuntu 17.10 and Julia 0.6.2. 
Version information for the Julia packages can be obtained from the supplement.

\subsubsection{Initial Fixed Cuts, Certificate Cuts, And Separation Cuts}
\label{sec:exp:cuts}

Recall from \cref{sec:soft:alg} that Pajarito by default uses three different types of $\K^*$ cuts: initial fixed cuts, certificate cuts, and separation cuts. 
For both the iterative and MSD methods, we compare the following four important algorithmic variants of OA that use different combinations of these three cut types.
\begin{description}
\item[c] means initial fixed cuts on linear primitive cones only, and certificate cuts on nonpolyhedral primitive cones.
\item[cs] means initial fixed cuts on linear primitive cones only, and certificate cuts on nonpolyhedral primitive cones, and separation cuts when apparently needed for convergence. 
The separation cuts allow us to cut off significantly infeasible OA solutions, so Pajarito can also obtain (approximately) feasible solutions from OA solutions found by the MILP solver.
\item[ics] means initial fixed cuts on all primitive cones, certificate cuts on nonpolyhedral primitive cones, and separation cuts when apparently needed for convergence. 
This is Pajarito's default approach, as described in \cref{sec:soft:alg}.
\item[is] means initial fixed cuts on all primitive cones, and separation cuts only. 
No conic solver is used, hence all (approximately) feasible solutions found are OA solutions.
\end{description}

\begin{figure}[!htbp]
\begin{subfigure}{0.48\textwidth}
\centering\begin{tikzpicture}[font = \footnotesize]
\begin{axis}[xmax=6, xtick={6}]
\addplot+ coordinates {
(1.0, 0.010869565217391304)
(1.0, 0.021739130434782608)
(1.0, 0.03260869565217391)
(1.0, 0.043478260869565216)
(1.0, 0.05434782608695652)
(1.0, 0.06521739130434782)
(1.0, 0.07608695652173914)
(1.0, 0.08695652173913043)
(1.0, 0.09782608695652174)
(1.0, 0.10869565217391304)
(1.0, 0.11956521739130435)
(1.0, 0.13043478260869565)
(1.0, 0.14130434782608695)
(1.0, 0.15217391304347827)
(1.0, 0.16304347826086957)
(1.0, 0.17391304347826086)
(1.0, 0.18478260869565216)
(1.0, 0.1956521739130435)
(1.0, 0.20652173913043478)
(1.0, 0.21739130434782608)
(1.0, 0.22826086956521738)
(1.0, 0.2391304347826087)
(1.0, 0.25)
(1.0, 0.2608695652173913)
(1.0, 0.2717391304347826)
(1.0028701885330489, 0.2826086956521739)
(1.0111721131052482, 0.29347826086956524)
(1.0138624290796316, 0.30434782608695654)
(1.0185235751219213, 0.31521739130434784)
(1.0214722914065741, 0.32608695652173914)
(1.0228783390718357, 0.33695652173913043)
(1.0253094063466426, 0.34782608695652173)
(1.0262745732177232, 0.358695652173913)
(1.03269546923343, 0.3695652173913043)
(1.0336469941953812, 0.3804347826086957)
(1.0467618302104935, 0.391304347826087)
(1.0630242900443991, 0.40217391304347827)
(1.0642246486926843, 0.41304347826086957)
(1.0681977284300215, 0.42391304347826086)
(1.0802471616235247, 0.43478260869565216)
(1.081730119117587, 0.44565217391304346)
(1.0949807946593466, 0.45652173913043476)
(1.096600138687109, 0.4673913043478261)
(1.1029527115652074, 0.4782608695652174)
(1.1043860045424123, 0.4891304347826087)
(1.1331035793124438, 0.5)
(1.1392868315157438, 0.5108695652173914)
(1.143002346102488, 0.5217391304347826)
(1.1482639280818911, 0.532608695652174)
(1.1488392722633851, 0.5434782608695652)
(1.1508464587808047, 0.5543478260869565)
(1.15279672258531, 0.5652173913043478)
(1.178841723586505, 0.5760869565217391)
(1.1834179998199905, 0.5869565217391305)
(1.2015016374427614, 0.5978260869565217)
(1.2250635545652646, 0.6086956521739131)
(1.2678437829952858, 0.6195652173913043)
(1.2758199968230621, 0.6304347826086957)
(1.3204090822738686, 0.6413043478260869)
(1.3306301498627102, 0.6521739130434783)
(1.3714052155592977, 0.6630434782608695)
(1.3733480590911893, 0.6739130434782609)
(1.4105506061567075, 0.6847826086956522)
(1.4202953691762235, 0.6956521739130435)
(1.464853298178419, 0.7065217391304348)
(1.472457143099408, 0.717391304347826)
(1.663057141019184, 0.7282608695652174)
(1.7137258948804277, 0.7391304347826086)
(1.7234867767872368, 0.75)
(1.7784422976730412, 0.7608695652173914)
(2.053129308867946, 0.7717391304347826)
(2.0607605556674726, 0.782608695652174)
(2.100923129444985, 0.7934782608695652)
(2.1234749035870806, 0.8043478260869565)
(2.2186679817878208, 0.8152173913043478)
(2.3265879499180366, 0.8260869565217391)
(2.5559526817319913, 0.8369565217391305)
(2.853940420367767, 0.8478260869565217)
(2.8991517110789893, 0.8586956521739131)
(2.986942286389168, 0.8695652173913043)
(3.6373227989225807, 0.8804347826086957)
(4.116840046374086, 0.8913043478260869)
(4.639117094915487, 0.9021739130434783)
(5.110798107927541, 0.9130434782608695)
(10.221596215855081, 0.9239130434782609)
(10.221596215855081, 0.9347826086956522)
(10.221596215855081, 0.9456521739130435)
(10.221596215855081, 0.9565217391304348)
(10.221596215855081, 0.967391304347826)
(10.221596215855081, 0.9782608695652174)
(10.221596215855081, 0.9891304347826086)
(10.221596215855081, 1.0)
(10.221596215855081, 1.0108695652173914)
} node[below right] at (1.5, 0.7173913) {is};
\addplot+ coordinates {
(1.0, 0.010869565217391304)
(1.0, 0.021739130434782608)
(1.0, 0.03260869565217391)
(1.0, 0.043478260869565216)
(1.0, 0.05434782608695652)
(1.0, 0.06521739130434782)
(1.0, 0.07608695652173914)
(1.0, 0.08695652173913043)
(1.0, 0.09782608695652174)
(1.0, 0.10869565217391304)
(1.0, 0.11956521739130435)
(1.0, 0.13043478260869565)
(1.0, 0.14130434782608695)
(1.0, 0.15217391304347827)
(1.0, 0.16304347826086957)
(1.0, 0.17391304347826086)
(1.0, 0.18478260869565216)
(1.0, 0.1956521739130435)
(1.0, 0.20652173913043478)
(1.0, 0.21739130434782608)
(1.0, 0.22826086956521738)
(1.0, 0.2391304347826087)
(1.0, 0.25)
(1.0, 0.2608695652173913)
(1.0, 0.2717391304347826)
(1.0, 0.2826086956521739)
(1.0, 0.29347826086956524)
(1.0, 0.30434782608695654)
(1.0, 0.31521739130434784)
(1.0, 0.32608695652173914)
(1.0, 0.33695652173913043)
(1.0, 0.34782608695652173)
(1.0, 0.358695652173913)
(1.0, 0.3695652173913043)
(1.0, 0.3804347826086957)
(1.0, 0.391304347826087)
(1.0, 0.40217391304347827)
(1.0, 0.41304347826086957)
(1.0, 0.42391304347826086)
(1.0, 0.43478260869565216)
(1.0, 0.44565217391304346)
(1.0, 0.45652173913043476)
(1.0, 0.4673913043478261)
(1.0, 0.4782608695652174)
(1.0, 0.4891304347826087)
(1.0, 0.5)
(1.0, 0.5108695652173914)
(1.0, 0.5217391304347826)
(1.0, 0.532608695652174)
(1.0, 0.5434782608695652)
(1.0, 0.5543478260869565)
(1.0, 0.5652173913043478)
(1.0, 0.5760869565217391)
(1.0, 0.5869565217391305)
(1.0, 0.5978260869565217)
(1.0, 0.6086956521739131)
(1.0, 0.6195652173913043)
(1.0, 0.6304347826086957)
(1.0, 0.6413043478260869)
(1.0, 0.6521739130434783)
(1.0, 0.6630434782608695)
(1.0, 0.6739130434782609)
(1.0, 0.6847826086956522)
(1.0, 0.6956521739130435)
(1.0, 0.7065217391304348)
(1.0, 0.717391304347826)
(1.0, 0.7282608695652174)
(1.0009899568910234, 0.7391304347826086)
(1.001602642567803, 0.75)
(1.0073475238469867, 0.7608695652173914)
(1.0089702540227063, 0.7717391304347826)
(1.0118645258578638, 0.782608695652174)
(1.01308266613699, 0.7934782608695652)
(1.0133066516864286, 0.8043478260869565)
(1.0158217668714256, 0.8152173913043478)
(1.016625251749096, 0.8260869565217391)
(1.0183217976773244, 0.8369565217391305)
(1.0200345869250387, 0.8478260869565217)
(1.0216051917499667, 0.8586956521739131)
(1.0335396451692136, 0.8695652173913043)
(1.0445054491249388, 0.8804347826086957)
(1.0515101941286291, 0.8913043478260869)
(1.0620376716049382, 0.9021739130434783)
(1.0880415405691723, 0.9130434782608695)
(1.1454391279724099, 0.9239130434782609)
(1.3535024683417325, 0.9347826086956522)
(1.8136425581824687, 0.9456521739130435)
(2.532171405555403, 0.9565217391304348)
(3.139644925337071, 0.967391304347826)
(10.221596215855081, 0.9782608695652174)
(10.221596215855081, 0.9891304347826086)
(10.221596215855081, 1.0)
(10.221596215855081, 1.0108695652173914)
} node[below right] at (1.5, 0.939) {ics};
\end{axis}
\end{tikzpicture}
\caption{Iter - execution time.}
\label{fig:types:itertime}
\end{subfigure}%
\hspace{0.04\textwidth}%
\begin{subfigure}{0.48\textwidth}
\centering\begin{tikzpicture}[font = \footnotesize]
\begin{axis}[xmax=4, xtick={4}]
\addplot+ coordinates {
(1.0, 0.010752688172043012)
(1.0, 0.021505376344086023)
(1.0, 0.03225806451612903)
(1.0, 0.043010752688172046)
(1.0, 0.053763440860215055)
(1.0, 0.06451612903225806)
(1.0, 0.07526881720430108)
(1.0, 0.08602150537634409)
(1.0, 0.0967741935483871)
(1.0, 0.10752688172043011)
(1.0, 0.11827956989247312)
(1.0, 0.12903225806451613)
(1.0, 0.13978494623655913)
(1.0, 0.15053763440860216)
(1.0, 0.16129032258064516)
(1.0, 0.17204301075268819)
(1.0, 0.1827956989247312)
(1.0, 0.1935483870967742)
(1.0, 0.20430107526881722)
(1.0, 0.21505376344086022)
(1.0, 0.22580645161290322)
(1.0, 0.23655913978494625)
(1.0, 0.24731182795698925)
(1.0, 0.25806451612903225)
(1.0, 0.26881720430107525)
(1.0, 0.27956989247311825)
(1.0, 0.2903225806451613)
(1.0, 0.3010752688172043)
(1.0, 0.3118279569892473)
(1.0, 0.3225806451612903)
(1.0, 0.3333333333333333)
(1.0, 0.34408602150537637)
(1.0, 0.3548387096774194)
(1.0, 0.3655913978494624)
(1.0, 0.3763440860215054)
(1.0, 0.3870967741935484)
(1.0, 0.3978494623655914)
(1.0, 0.40860215053763443)
(1.0, 0.41935483870967744)
(1.0, 0.43010752688172044)
(1.0, 0.44086021505376344)
(1.0, 0.45161290322580644)
(1.0, 0.46236559139784944)
(1.0, 0.4731182795698925)
(1.0, 0.4838709677419355)
(1.0, 0.4946236559139785)
(1.0, 0.5053763440860215)
(1.0, 0.5161290322580645)
(1.0, 0.5268817204301075)
(1.0, 0.5376344086021505)
(1.0, 0.5483870967741935)
(1.0, 0.5591397849462365)
(1.0, 0.5698924731182796)
(1.0, 0.5806451612903226)
(1.0, 0.5913978494623656)
(1.0, 0.6021505376344086)
(1.0, 0.6129032258064516)
(1.0, 0.6236559139784946)
(1.0, 0.6344086021505376)
(1.0005235851294498, 0.6451612903225806)
(1.0023236483013245, 0.6559139784946236)
(1.002480151417035, 0.6666666666666666)
(1.0048338271666606, 0.6774193548387096)
(1.009114712175315, 0.6881720430107527)
(1.0093004162831642, 0.6989247311827957)
(1.0113319546447344, 0.7096774193548387)
(1.01556489798526, 0.7204301075268817)
(1.0188200721041882, 0.7311827956989247)
(1.032133334276268, 0.7419354838709677)
(1.0687467866443194, 0.7526881720430108)
(1.0828695678910443, 0.7634408602150538)
(1.1084546634254429, 0.7741935483870968)
(1.1400770058223118, 0.7849462365591398)
(1.1804337191976282, 0.7956989247311828)
(1.29269286878582, 0.8064516129032258)
(1.342245864086877, 0.8172043010752689)
(1.9312314765281247, 0.8279569892473119)
(1.9918007909364603, 0.8387096774193549)
(2.0179390477617583, 0.8494623655913979)
(2.4334297285334516, 0.8602150537634409)
(2.479757997996459, 0.8709677419354839)
(3.2069490908519946, 0.8817204301075269)
(3.4407258195523833, 0.8924731182795699)
(3.515965368779121, 0.9032258064516129)
(7.031930737558242, 0.9139784946236559)
(7.031930737558242, 0.9247311827956989)
(7.031930737558242, 0.9354838709677419)
(7.031930737558242, 0.946236559139785)
(7.031930737558242, 0.956989247311828)
(7.031930737558242, 0.967741935483871)
(7.031930737558242, 0.978494623655914)
(7.031930737558242, 0.989247311827957)
(7.031930737558242, 1.0)
(7.031930737558242, 1.010752688172043)
} node[below right] at (1.5, 0.81720430) {is};
\addplot+ coordinates {
(1.0, 0.010752688172043012)
(1.0, 0.021505376344086023)
(1.0, 0.03225806451612903)
(1.0, 0.043010752688172046)
(1.0, 0.053763440860215055)
(1.0, 0.06451612903225806)
(1.0, 0.07526881720430108)
(1.0, 0.08602150537634409)
(1.0, 0.0967741935483871)
(1.0, 0.10752688172043011)
(1.0, 0.11827956989247312)
(1.0, 0.12903225806451613)
(1.0, 0.13978494623655913)
(1.0, 0.15053763440860216)
(1.0, 0.16129032258064516)
(1.0, 0.17204301075268819)
(1.0, 0.1827956989247312)
(1.0, 0.1935483870967742)
(1.0, 0.20430107526881722)
(1.0, 0.21505376344086022)
(1.0, 0.22580645161290322)
(1.0, 0.23655913978494625)
(1.0, 0.24731182795698925)
(1.0, 0.25806451612903225)
(1.0, 0.26881720430107525)
(1.0, 0.27956989247311825)
(1.0, 0.2903225806451613)
(1.0, 0.3010752688172043)
(1.0, 0.3118279569892473)
(1.0, 0.3225806451612903)
(1.0, 0.3333333333333333)
(1.0, 0.34408602150537637)
(1.0, 0.3548387096774194)
(1.0, 0.3655913978494624)
(1.0026383737819247, 0.3763440860215054)
(1.0027311256583658, 0.3870967741935484)
(1.0031003423304061, 0.3978494623655914)
(1.0042193807320337, 0.40860215053763443)
(1.0049813717932612, 0.41935483870967744)
(1.0051514049554422, 0.43010752688172044)
(1.0060374818356408, 0.44086021505376344)
(1.0066987947040242, 0.45161290322580644)
(1.0067701905667932, 0.46236559139784944)
(1.0092294773649557, 0.4731182795698925)
(1.0098136068153896, 0.4838709677419355)
(1.0105835120698798, 0.4946236559139785)
(1.0118044216779292, 0.5053763440860215)
(1.0121808853016934, 0.5161290322580645)
(1.0122517101926973, 0.5268817204301075)
(1.0129075118323905, 0.5376344086021505)
(1.023637245501315, 0.5483870967741935)
(1.0266333876310334, 0.5591397849462365)
(1.0273589719087979, 0.5698924731182796)
(1.0287866689056466, 0.5806451612903226)
(1.0300636348945378, 0.5913978494623656)
(1.0313616009941788, 0.6021505376344086)
(1.0321338349606388, 0.6129032258064516)
(1.0334677806645096, 0.6236559139784946)
(1.0352715935184842, 0.6344086021505376)
(1.0353064780204475, 0.6451612903225806)
(1.0358897224490056, 0.6559139784946236)
(1.0375472849127572, 0.6666666666666666)
(1.0383727075627964, 0.6774193548387096)
(1.0398791421250968, 0.6881720430107527)
(1.0400410462105867, 0.6989247311827957)
(1.0453569648801653, 0.7096774193548387)
(1.0482133747265603, 0.7204301075268817)
(1.0512873020486593, 0.7311827956989247)
(1.0525418927341212, 0.7419354838709677)
(1.0557362774835768, 0.7526881720430108)
(1.0643361888786105, 0.7634408602150538)
(1.0656796232683066, 0.7741935483870968)
(1.0687576064098094, 0.7849462365591398)
(1.0698234806483315, 0.7956989247311828)
(1.0782030693736275, 0.8064516129032258)
(1.0834276461475116, 0.8172043010752689)
(1.1081432099653685, 0.8279569892473119)
(1.1147326936893587, 0.8387096774193549)
(1.118778349398442, 0.8494623655913979)
(1.1263565941948772, 0.8602150537634409)
(1.1283261045642092, 0.8709677419354839)
(1.1401985937271983, 0.8817204301075269)
(1.1442438835503363, 0.8924731182795699)
(1.1457770167480141, 0.9032258064516129)
(1.1820321421736273, 0.9139784946236559)
(1.1999212958455523, 0.9247311827956989)
(1.255315108139884, 0.9354838709677419)
(1.269881499859377, 0.946236559139785)
(1.30627349360533, 0.956989247311828)
(1.5500639712419038, 0.967741935483871)
(1.6538037949904356, 0.978494623655914)
(1.7135736321374966, 0.989247311827957)
(7.031930737558242, 1.0)
(7.031930737558242, 1.010752688172043)
} node[below right] at (1.5, 0.965) {ics};
\end{axis}
\end{tikzpicture}
\caption{MSD - execution time.}
\label{fig:types:msdtime}
\end{subfigure}
\\[0.7em]
\begin{subfigure}{0.48\textwidth}
\centering\begin{tikzpicture}[font = \footnotesize]
\begin{axis}[xmax=6, xtick={6}]
\addplot+ coordinates {
(1.0, 0.010869565217391304)
(1.0, 0.021739130434782608)
(1.0, 0.03260869565217391)
(1.0, 0.043478260869565216)
(1.0, 0.05434782608695652)
(1.0, 0.06521739130434782)
(1.0, 0.07608695652173914)
(1.0, 0.08695652173913043)
(1.0, 0.09782608695652174)
(1.0, 0.10869565217391304)
(1.0, 0.11956521739130435)
(1.015625, 0.13043478260869565)
(1.1111111111111112, 0.14130434782608695)
(1.125, 0.15217391304347827)
(1.1428571428571428, 0.16304347826086957)
(1.1428571428571428, 0.17391304347826086)
(1.2, 0.18478260869565216)
(1.2, 0.1956521739130435)
(1.3506493506493507, 0.20652173913043478)
(1.4, 0.21739130434782608)
(1.5, 0.22826086956521738)
(1.5, 0.2391304347826087)
(1.5, 0.25)
(1.7333333333333334, 0.2608695652173913)
(1.75, 0.2717391304347826)
(1.875, 0.2826086956521739)
(1.9259259259259258, 0.29347826086956524)
(2.0, 0.30434782608695654)
(2.0454545454545454, 0.31521739130434784)
(2.25, 0.32608695652173914)
(2.263157894736842, 0.33695652173913043)
(2.2666666666666666, 0.34782608695652173)
(2.3076923076923075, 0.358695652173913)
(2.3333333333333335, 0.3695652173913043)
(2.4, 0.3804347826086957)
(2.4, 0.391304347826087)
(2.5, 0.40217391304347827)
(2.5, 0.41304347826086957)
(2.6666666666666665, 0.42391304347826086)
(2.6666666666666665, 0.43478260869565216)
(3.0, 0.44565217391304346)
(3.0, 0.45652173913043476)
(3.0, 0.4673913043478261)
(3.0, 0.4782608695652174)
(3.0, 0.4891304347826087)
(3.0, 0.5)
(3.1, 0.5108695652173914)
(3.2, 0.5217391304347826)
(3.2, 0.532608695652174)
(3.2, 0.5434782608695652)
(3.2, 0.5543478260869565)
(3.25, 0.5652173913043478)
(3.75, 0.5760869565217391)
(3.75, 0.5869565217391305)
(3.75, 0.5978260869565217)
(3.875, 0.6086956521739131)
(4.0, 0.6195652173913043)
(4.0, 0.6304347826086957)
(4.0, 0.6413043478260869)
(4.0, 0.6521739130434783)
(4.166666666666667, 0.6630434782608695)
(4.333333333333333, 0.6739130434782609)
(4.333333333333333, 0.6847826086956522)
(4.333333333333333, 0.6956521739130435)
(4.333333333333333, 0.7065217391304348)
(4.333333333333333, 0.717391304347826)
(4.333333333333333, 0.7282608695652174)
(4.333333333333333, 0.7391304347826086)
(4.333333333333333, 0.75)
(4.333333333333333, 0.7608695652173914)
(4.333333333333333, 0.7717391304347826)
(4.333333333333333, 0.782608695652174)
(4.333333333333333, 0.7934782608695652)
(4.333333333333333, 0.8043478260869565)
(4.375, 0.8152173913043478)
(4.5, 0.8260869565217391)
(4.5, 0.8369565217391305)
(4.5, 0.8478260869565217)
(5.0, 0.8586956521739131)
(5.0, 0.8695652173913043)
(5.0, 0.8804347826086957)
(5.0, 0.8913043478260869)
(5.333333333333333, 0.9021739130434783)
(6.0, 0.9130434782608695)
(12.0, 0.9239130434782609)
(12.0, 0.9347826086956522)
(12.0, 0.9456521739130435)
(12.0, 0.9565217391304348)
(12.0, 0.967391304347826)
(12.0, 0.9782608695652174)
(12.0, 0.9891304347826086)
(12.0, 1.0)
(12.0, 1.0108695652173914)
} node[below right] at (3.5, 0.56521739) {is};
\addplot+ coordinates {
(1.0, 0.010869565217391304)
(1.0, 0.021739130434782608)
(1.0, 0.03260869565217391)
(1.0, 0.043478260869565216)
(1.0, 0.05434782608695652)
(1.0, 0.06521739130434782)
(1.0, 0.07608695652173914)
(1.0, 0.08695652173913043)
(1.0, 0.09782608695652174)
(1.0, 0.10869565217391304)
(1.0, 0.11956521739130435)
(1.0, 0.13043478260869565)
(1.0, 0.14130434782608695)
(1.0, 0.15217391304347827)
(1.0, 0.16304347826086957)
(1.0, 0.17391304347826086)
(1.0, 0.18478260869565216)
(1.0, 0.1956521739130435)
(1.0, 0.20652173913043478)
(1.0, 0.21739130434782608)
(1.0, 0.22826086956521738)
(1.0, 0.2391304347826087)
(1.0, 0.25)
(1.0, 0.2608695652173913)
(1.0, 0.2717391304347826)
(1.0, 0.2826086956521739)
(1.0, 0.29347826086956524)
(1.0, 0.30434782608695654)
(1.0, 0.31521739130434784)
(1.0, 0.32608695652173914)
(1.0, 0.33695652173913043)
(1.0, 0.34782608695652173)
(1.0, 0.358695652173913)
(1.0, 0.3695652173913043)
(1.0, 0.3804347826086957)
(1.0, 0.391304347826087)
(1.0, 0.40217391304347827)
(1.0, 0.41304347826086957)
(1.0, 0.42391304347826086)
(1.0, 0.43478260869565216)
(1.0, 0.44565217391304346)
(1.0, 0.45652173913043476)
(1.0, 0.4673913043478261)
(1.0, 0.4782608695652174)
(1.0, 0.4891304347826087)
(1.0, 0.5)
(1.0, 0.5108695652173914)
(1.0, 0.5217391304347826)
(1.0, 0.532608695652174)
(1.0, 0.5434782608695652)
(1.0, 0.5543478260869565)
(1.0, 0.5652173913043478)
(1.0, 0.5760869565217391)
(1.0, 0.5869565217391305)
(1.0, 0.5978260869565217)
(1.0, 0.6086956521739131)
(1.0, 0.6195652173913043)
(1.0, 0.6304347826086957)
(1.0, 0.6413043478260869)
(1.0, 0.6521739130434783)
(1.0, 0.6630434782608695)
(1.0, 0.6739130434782609)
(1.0, 0.6847826086956522)
(1.0, 0.6956521739130435)
(1.0, 0.7065217391304348)
(1.0, 0.717391304347826)
(1.0, 0.7282608695652174)
(1.0, 0.7391304347826086)
(1.0, 0.75)
(1.0, 0.7608695652173914)
(1.0, 0.7717391304347826)
(1.0, 0.782608695652174)
(1.0, 0.7934782608695652)
(1.0, 0.8043478260869565)
(1.0, 0.8152173913043478)
(1.0, 0.8260869565217391)
(1.0, 0.8369565217391305)
(1.0, 0.8478260869565217)
(1.0, 0.8586956521739131)
(1.0, 0.8695652173913043)
(1.0, 0.8804347826086957)
(1.0, 0.8913043478260869)
(1.0, 0.9021739130434783)
(1.0434782608695652, 0.9130434782608695)
(1.25, 0.9239130434782609)
(1.303030303030303, 0.9347826086956522)
(1.3857142857142857, 0.9456521739130435)
(1.4642857142857142, 0.9565217391304348)
(1.8484848484848484, 0.967391304347826)
(12.0, 0.9782608695652174)
(12.0, 0.9891304347826086)
(12.0, 1.0)
(12.0, 1.0108695652173914)
} node[below right] at (3.5, 0.9673913) {ics};
\end{axis}
\end{tikzpicture}
\caption{Iter - iteration count.}
\label{fig:types:iteriters}
\end{subfigure}%
\hspace{0.04\textwidth}%
\begin{subfigure}{0.48\textwidth}
\centering\begin{tikzpicture}[font = \footnotesize]
\begin{axis}[xmax=12, xtick={12}]
\addplot+ coordinates {
(1.0, 0.010752688172043012)
(1.0, 0.021505376344086023)
(1.0, 0.03225806451612903)
(1.0, 0.043010752688172046)
(1.0, 0.053763440860215055)
(1.0, 0.06451612903225806)
(1.0, 0.07526881720430108)
(1.0, 0.08602150537634409)
(1.0, 0.0967741935483871)
(1.0, 0.10752688172043011)
(1.0, 0.11827956989247312)
(1.0, 0.12903225806451613)
(1.0, 0.13978494623655913)
(1.0, 0.15053763440860216)
(1.0, 0.16129032258064516)
(1.0, 0.17204301075268819)
(1.0, 0.1827956989247312)
(1.0, 0.1935483870967742)
(1.0, 0.20430107526881722)
(1.0, 0.21505376344086022)
(1.0, 0.22580645161290322)
(1.0, 0.23655913978494625)
(1.0, 0.24731182795698925)
(1.0, 0.25806451612903225)
(1.0, 0.26881720430107525)
(1.0, 0.27956989247311825)
(1.0, 0.2903225806451613)
(1.000402738622634, 0.3010752688172043)
(1.0022918258212377, 0.3118279569892473)
(1.0045931758530184, 0.3225806451612903)
(1.0506838294448915, 0.3333333333333333)
(1.0847668545494473, 0.34408602150537637)
(1.1180281245062411, 0.3548387096774194)
(1.184269107098689, 0.3655913978494624)
(1.2142857142857142, 0.3763440860215054)
(1.2528397139251157, 0.3870967741935484)
(1.2794561933534743, 0.3978494623655914)
(1.3440993788819875, 0.40860215053763443)
(1.378787878787879, 0.41935483870967744)
(1.455681818181818, 0.43010752688172044)
(1.4854771784232366, 0.44086021505376344)
(1.5302197802197801, 0.45161290322580644)
(1.5316455696202531, 0.46236559139784944)
(1.5454545454545454, 0.4731182795698925)
(1.6363636363636365, 0.4838709677419355)
(1.6656995486782722, 0.4946236559139785)
(1.6928131138657454, 0.5053763440860215)
(1.7191011235955056, 0.5161290322580645)
(1.7239002635313196, 0.5268817204301075)
(1.7307692307692308, 0.5376344086021505)
(1.7384615384615385, 0.5483870967741935)
(1.8066604995374653, 0.5591397849462365)
(1.8181818181818181, 0.5698924731182796)
(1.9534883720930232, 0.5806451612903226)
(2.1341463414634148, 0.5913978494623656)
(2.3125, 0.6021505376344086)
(2.3406079501169135, 0.6129032258064516)
(2.3690944881889764, 0.6236559139784946)
(2.3850267379679146, 0.6344086021505376)
(2.5174603174603174, 0.6451612903225806)
(2.6137463697967087, 0.6559139784946236)
(2.6266924564796907, 0.6666666666666666)
(2.727272727272727, 0.6774193548387096)
(2.7889908256880735, 0.6881720430107527)
(2.8368336025848144, 0.6989247311827957)
(2.9, 0.7096774193548387)
(3.1973969631236443, 0.7204301075268817)
(3.265193370165746, 0.7311827956989247)
(3.4141689373297, 0.7419354838709677)
(3.4628879892037787, 0.7526881720430108)
(4.0, 0.7634408602150538)
(4.281524926686217, 0.7741935483870968)
(4.7272727272727275, 0.7849462365591398)
(5.242424242424242, 0.7956989247311828)
(5.9, 0.8064516129032258)
(6.576923076923077, 0.8172043010752689)
(7.363636363636363, 0.8279569892473119)
(8.090909090909092, 0.8387096774193549)
(8.433333333333334, 0.8494623655913979)
(8.545454545454545, 0.8602150537634409)
(10.636363636363637, 0.8709677419354839)
(25.09090909090909, 0.8817204301075269)
(25.727272727272727, 0.8924731182795699)
(44.90909090909091, 0.9032258064516129)
(89.81818181818181, 0.9139784946236559)
(89.81818181818181, 0.9247311827956989)
(89.81818181818181, 0.9354838709677419)
(89.81818181818181, 0.946236559139785)
(89.81818181818181, 0.956989247311828)
(89.81818181818181, 0.967741935483871)
(89.81818181818181, 0.978494623655914)
(89.81818181818181, 0.989247311827957)
(89.81818181818181, 1.0)
(89.81818181818181, 1.010752688172043)
} node[below right] at (6.5, 0.8064516) {is};
\addplot+ coordinates {
(1.0, 0.010752688172043012)
(1.0, 0.021505376344086023)
(1.0, 0.03225806451612903)
(1.0, 0.043010752688172046)
(1.0, 0.053763440860215055)
(1.0, 0.06451612903225806)
(1.0, 0.07526881720430108)
(1.0, 0.08602150537634409)
(1.0, 0.0967741935483871)
(1.0, 0.10752688172043011)
(1.0, 0.11827956989247312)
(1.0, 0.12903225806451613)
(1.0, 0.13978494623655913)
(1.0, 0.15053763440860216)
(1.0, 0.16129032258064516)
(1.0, 0.17204301075268819)
(1.0, 0.1827956989247312)
(1.0, 0.1935483870967742)
(1.0, 0.20430107526881722)
(1.0, 0.21505376344086022)
(1.0, 0.22580645161290322)
(1.0, 0.23655913978494625)
(1.0, 0.24731182795698925)
(1.0, 0.25806451612903225)
(1.0, 0.26881720430107525)
(1.0, 0.27956989247311825)
(1.0, 0.2903225806451613)
(1.0, 0.3010752688172043)
(1.0, 0.3118279569892473)
(1.0, 0.3225806451612903)
(1.0, 0.3333333333333333)
(1.0, 0.34408602150537637)
(1.0, 0.3548387096774194)
(1.0, 0.3655913978494624)
(1.0, 0.3763440860215054)
(1.0, 0.3870967741935484)
(1.0, 0.3978494623655914)
(1.0, 0.40860215053763443)
(1.0, 0.41935483870967744)
(1.0, 0.43010752688172044)
(1.0, 0.44086021505376344)
(1.0, 0.45161290322580644)
(1.0, 0.46236559139784944)
(1.0, 0.4731182795698925)
(1.0, 0.4838709677419355)
(1.0, 0.4946236559139785)
(1.0, 0.5053763440860215)
(1.0, 0.5161290322580645)
(1.0, 0.5268817204301075)
(1.0, 0.5376344086021505)
(1.0, 0.5483870967741935)
(1.0, 0.5591397849462365)
(1.0, 0.5698924731182796)
(1.0, 0.5806451612903226)
(1.0, 0.5913978494623656)
(1.0, 0.6021505376344086)
(1.0, 0.6129032258064516)
(1.0, 0.6236559139784946)
(1.0, 0.6344086021505376)
(1.0, 0.6451612903225806)
(1.0, 0.6559139784946236)
(1.0, 0.6666666666666666)
(1.0, 0.6774193548387096)
(1.0, 0.6881720430107527)
(1.0, 0.6989247311827957)
(1.0, 0.7096774193548387)
(1.0, 0.7204301075268817)
(1.0, 0.7311827956989247)
(1.0248202563301032, 0.7419354838709677)
(1.0685268334402396, 0.7526881720430108)
(1.0953110634385534, 0.7634408602150538)
(1.1265859962406015, 0.7741935483870968)
(1.1974615064502705, 0.7849462365591398)
(1.2074417624117213, 0.7956989247311828)
(1.2100047361630546, 0.8064516129032258)
(1.2504524571905888, 0.8172043010752689)
(1.2506781874811614, 0.8279569892473119)
(1.3667628134013756, 0.8387096774193549)
(1.4039983340274886, 0.8494623655913979)
(1.5055633802816901, 0.8602150537634409)
(1.5405405405405406, 0.8709677419354839)
(1.6655844155844155, 0.8817204301075269)
(1.7015415216310292, 0.8924731182795699)
(1.7087880528431936, 0.9032258064516129)
(1.7897196261682242, 0.9139784946236559)
(1.7985520007827023, 0.9247311827956989)
(1.907744583898022, 0.9354838709677419)
(1.98, 0.946236559139785)
(2.566693613581245, 0.956989247311828)
(3.1479565442317643, 0.967741935483871)
(4.035714285714286, 0.978494623655914)
(5.4575043461464166, 0.989247311827957)
(89.81818181818181, 1.0)
(89.81818181818181, 1.010752688172043)
} node[below right] at (6.5, 0.993) {ics};
\end{axis}
\end{tikzpicture}
\caption{MSD - node count.}
\label{fig:types:msdnodes}
\end{subfigure}
\caption{$\K^*$ cut types performance profiles.}
\end{figure}
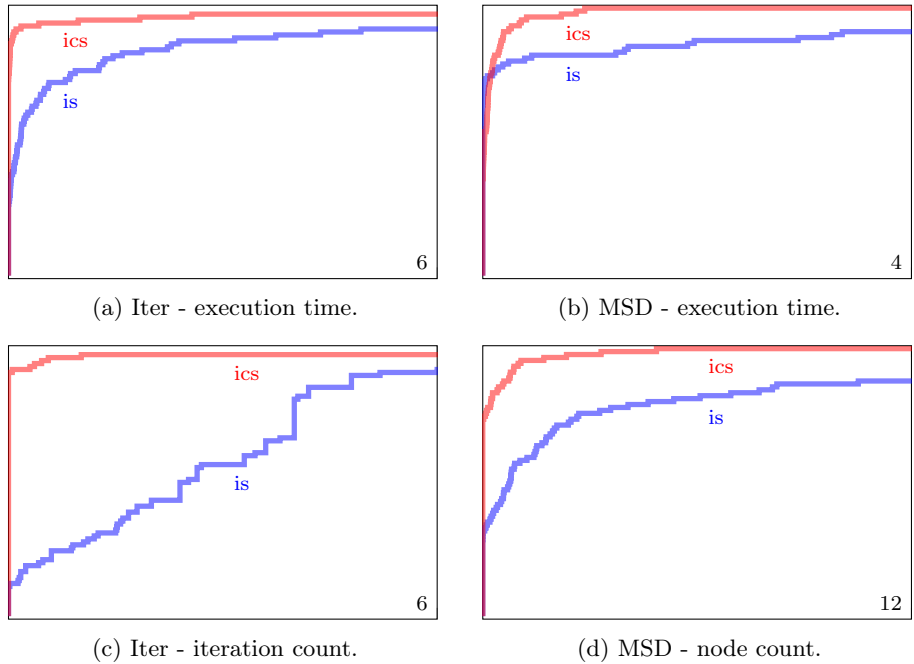

\Cref{tab:types:summ} summarizes the status counts and shifted geomeans of performance metrics on instance subsets. 
Although the MSD method is significantly faster than the iterative method, we see similar relative performances for the four types of cuts under iterative versus MSD.
When using certificate cuts only (`c'), Pajarito often failed to converge to the desired optimality gap (though it typically came very close), likely due to the inexactness of the certificates from the numerical continuous conic solver.
By also using separation cuts on repeated integer sub-solutions and accepting (approximately) conic feasible OA solutions as incumbents, Pajarito is able to converge on many more instances. 
Starting with initial fixed cuts (`ics') further increases Pajarito's robustness, particularly for the MSD method.
Comparing the `ics' variant with the separation-based variant with initial fixed cuts (`is'), we see significantly faster overall performance and fewer iterations or nodes when using the continuous conic solver and adding certificate cuts. 
The performance profiles \crefrange{fig:types:itertime}{fig:types:msdnodes} compare the execution times or iteration/node counts for the `ics' and `is' solvers, unambiguously demonstrating superiority of Pajarito's default `ics' method.

\begin{table}[!htb]
\footnotesize
\begin{tabular}{l c *{4}{S[table-format=2]} *{3}{S[table-format=2.2]} *{3}{S[table-format=2.2]} *{3}{S[table-format=2.2]}}
\toprule
& & \multicolumn{4}{l}{statuses} & \multicolumn{3}{l}{time (s)} & \multicolumn{3}{l}{subproblems} & \multicolumn{3}{l}{iters or nodes} 
\\
\cmidrule(lr){3-6} \cmidrule(lr){7-9} \cmidrule(lr){10-12} \cmidrule(lr){13-15}
& \rotatebox{90}{\rlap{cuts}} & {co} & {li} & {er} & {ex} & {aco} & {tco} & {all} & {aco} & {tco} & {all} & {aco} & {tco} & {all} \\
\midrule
\multirow{4}{*}{\rotatebox[origin=c]{90}{Iter}} 
& {c} & 
72 & 1 & 21 & 1 & 5.59 & 6.47 & 7.17 & 5.31 & 5.36 & 4.23 & 5.48 & 5.53 & 4.46 \\
& {cs} & 
88 & 1 & 3 & 3 & 5.56 & 12.72 & 14.38 & 5.31 & 6.97 & 6.36 & 5.48 & 7.22 & 6.84 \\
& {ics} & 
89 & 2 & 0 & 4 & 4.73 & 11.57 & 14.77 & 4.17 & 5.93 & 6.03 & 4.32 & 6.15 & 6.32 \\
& {is} &  
84 & 1 & 0 & 10 & 8.35 & 14.53 & 22.08 & {-} & {-} & {-} & 13.41 & 16.52 & 18.07 \\
\addlinespace
\multirow{4}{*}{\rotatebox[origin=c]{90}{MSD}} 
& {c} & 
76 & 0 & 18 & 1 & 2.37 & 3.40 & 3.50 & 12.63 & 15.80 & 12.70 & {$223$} & {$438$} & {$348$} \\
& {cs} & 
88 & 0 & 5 & 2 & 3.33 & 6.47 & 7.76 & 18.96 & 26.77 & 24.87 & {$295$} & {$843$} & {$815$} \\
& {ics} & 
92 & 0 & 1 & 2 & 2.20 & 6.31 & 6.52 & 15.62 & 24.58 & 24.95 & {$273$} & {$796$} & {$857$} \\
& {is} & 
84 & 1 & 0 & 10 & 3.12 & 5.29 & 7.49 & {-} & {-} & {-} & {$522$} & {$932$} & {$\llap{1}345$} \\
\bottomrule
\end{tabular}
\caption{$\K^*$ cut types performance summary.}
\label{tab:types:summ}
\end{table}

\subsubsection{Extreme Ray Disaggregation}
\label{sec:exp:disagg}

To test the efficacy of the $\K^*$ extreme ray disaggregation technique we describe in \cref{sec:tight:extr} and \cref{sec:spec}, we run Pajarito using only certificate cuts (the `c' variant described in \cref{sec:exp:cuts}), with and without disaggregation. 
Note that disabling disaggregation disables use of the second-order cone extended formulation, which has no benefit without disaggregation.

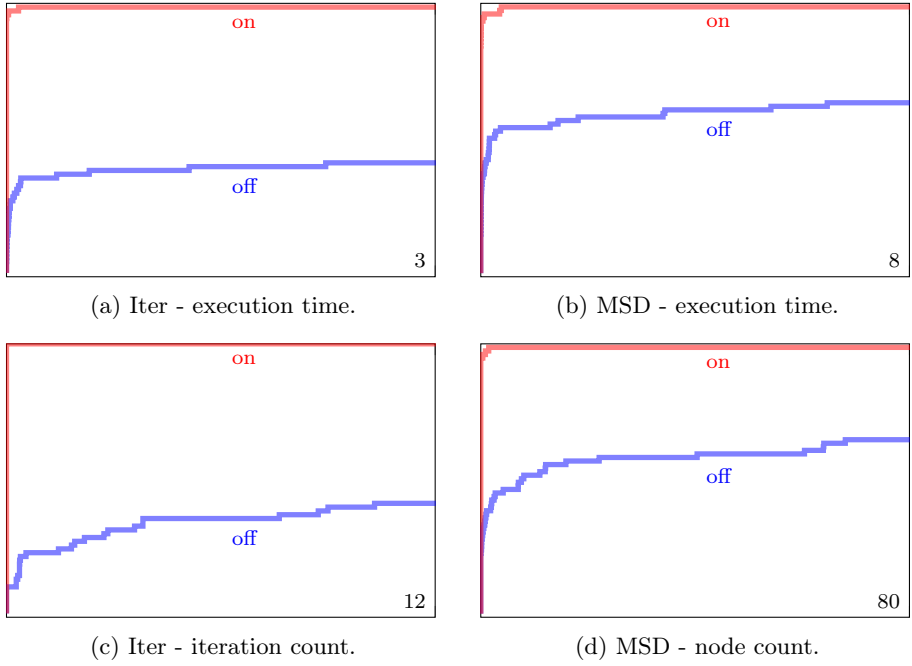
\begin{figure}[!htbp]
\begin{subfigure}{0.48\textwidth}
\centering\begin{tikzpicture}[font = \footnotesize]
\begin{axis}[xmax=3, xtick={3}]
\addplot+ coordinates {
(1.0, 0.013888888888888888)
(1.0, 0.027777777777777776)
(1.0, 0.041666666666666664)
(1.0, 0.05555555555555555)
(1.0006735047489936, 0.06944444444444445)
(1.0007315451368426, 0.08333333333333333)
(1.0009382644328861, 0.09722222222222222)
(1.001150431462039, 0.1111111111111111)
(1.0015948768081298, 0.125)
(1.0019249512903123, 0.1388888888888889)
(1.0025485438530253, 0.1527777777777778)
(1.0071752000118717, 0.16666666666666666)
(1.007650472066426, 0.18055555555555555)
(1.0091257327649903, 0.19444444444444445)
(1.0092464759914992, 0.20833333333333334)
(1.010700276529313, 0.2222222222222222)
(1.0121874846904153, 0.2361111111111111)
(1.0128201466184492, 0.25)
(1.0178600121489347, 0.2638888888888889)
(1.0178841988033156, 0.2777777777777778)
(1.0295779950470902, 0.2916666666666667)
(1.0367022485647948, 0.3055555555555556)
(1.0480367893416622, 0.3194444444444444)
(1.0561478934156476, 0.3333333333333333)
(1.0652350900567356, 0.3472222222222222)
(1.066676500012974, 0.3611111111111111)
(1.2338050015543585, 0.375)
(1.3861029071343212, 0.3888888888888889)
(1.85185737121092, 0.4027777777777778)
(2.489175256180064, 0.4166666666666667)
(51.8397099002809, 0.4305555555555556)
(52.0466808220347, 0.4444444444444444)
(54.62564221212095, 0.4583333333333333)
(109.2512844242419, 0.4722222222222222)
(109.2512844242419, 0.4861111111111111)
(109.2512844242419, 0.5)
(109.2512844242419, 0.5138888888888888)
(109.2512844242419, 0.5277777777777778)
(109.2512844242419, 0.5416666666666666)
(109.2512844242419, 0.5555555555555556)
(109.2512844242419, 0.5694444444444444)
(109.2512844242419, 0.5833333333333334)
(109.2512844242419, 0.5972222222222222)
(109.2512844242419, 0.6111111111111112)
(109.2512844242419, 0.625)
(109.2512844242419, 0.6388888888888888)
(109.2512844242419, 0.6527777777777778)
(109.2512844242419, 0.6666666666666666)
(109.2512844242419, 0.6805555555555556)
(109.2512844242419, 0.6944444444444444)
(109.2512844242419, 0.7083333333333334)
(109.2512844242419, 0.7222222222222222)
(109.2512844242419, 0.7361111111111112)
(109.2512844242419, 0.75)
(109.2512844242419, 0.7638888888888888)
(109.2512844242419, 0.7777777777777778)
(109.2512844242419, 0.7916666666666666)
(109.2512844242419, 0.8055555555555556)
(109.2512844242419, 0.8194444444444444)
(109.2512844242419, 0.8333333333333334)
(109.2512844242419, 0.8472222222222222)
(109.2512844242419, 0.8611111111111112)
(109.2512844242419, 0.875)
(109.2512844242419, 0.8888888888888888)
(109.2512844242419, 0.9027777777777778)
(109.2512844242419, 0.9166666666666666)
(109.2512844242419, 0.9305555555555556)
(109.2512844242419, 0.9444444444444444)
(109.2512844242419, 0.9583333333333334)
(109.2512844242419, 0.9722222222222222)
(109.2512844242419, 0.9861111111111112)
(109.2512844242419, 1.0)
(109.2512844242419, 1.0138888888888888)
} node[below right] at (2, 0.40277777) {off};
\addplot+ coordinates {
(1.0, 0.013888888888888888)
(1.0, 0.027777777777777776)
(1.0, 0.041666666666666664)
(1.0, 0.05555555555555555)
(1.0, 0.06944444444444445)
(1.0, 0.08333333333333333)
(1.0, 0.09722222222222222)
(1.0, 0.1111111111111111)
(1.0, 0.125)
(1.0, 0.1388888888888889)
(1.0, 0.1527777777777778)
(1.0, 0.16666666666666666)
(1.0, 0.18055555555555555)
(1.0, 0.19444444444444445)
(1.0, 0.20833333333333334)
(1.0, 0.2222222222222222)
(1.0, 0.2361111111111111)
(1.0, 0.25)
(1.0, 0.2638888888888889)
(1.0, 0.2777777777777778)
(1.0, 0.2916666666666667)
(1.0, 0.3055555555555556)
(1.0, 0.3194444444444444)
(1.0, 0.3333333333333333)
(1.0, 0.3472222222222222)
(1.0, 0.3611111111111111)
(1.0, 0.375)
(1.0, 0.3888888888888889)
(1.0, 0.4027777777777778)
(1.0, 0.4166666666666667)
(1.0, 0.4305555555555556)
(1.0, 0.4444444444444444)
(1.0, 0.4583333333333333)
(1.0, 0.4722222222222222)
(1.0, 0.4861111111111111)
(1.0, 0.5)
(1.0, 0.5138888888888888)
(1.0, 0.5277777777777778)
(1.0, 0.5416666666666666)
(1.0, 0.5555555555555556)
(1.0, 0.5694444444444444)
(1.0, 0.5833333333333334)
(1.0, 0.5972222222222222)
(1.0, 0.6111111111111112)
(1.0, 0.625)
(1.0, 0.6388888888888888)
(1.0, 0.6527777777777778)
(1.0, 0.6666666666666666)
(1.0, 0.6805555555555556)
(1.0, 0.6944444444444444)
(1.0, 0.7083333333333334)
(1.0, 0.7222222222222222)
(1.0, 0.7361111111111112)
(1.0, 0.75)
(1.0, 0.7638888888888888)
(1.0, 0.7777777777777778)
(1.0, 0.7916666666666666)
(1.0, 0.8055555555555556)
(1.0, 0.8194444444444444)
(1.0, 0.8333333333333334)
(1.0, 0.8472222222222222)
(1.0, 0.8611111111111112)
(1.0, 0.875)
(1.0, 0.8888888888888888)
(1.0, 0.9027777777777778)
(1.0, 0.9166666666666666)
(1.0, 0.9305555555555556)
(1.0, 0.9444444444444444)
(1.0009381791616012, 0.9583333333333334)
(1.0096403567555166, 0.9722222222222222)
(1.0562127413972395, 0.9861111111111112)
(3.724432375157983, 1.0)
(109.2512844242419, 1.0138888888888888)
} node[below right] at (2, 0.9861111) {on};
\end{axis}
\end{tikzpicture}
\caption{Iter - execution time.}
\label{fig:disagg:itertime}
\end{subfigure}%
\hspace{0.04\textwidth}%
\begin{subfigure}{0.48\textwidth}
\centering\begin{tikzpicture}[font = \footnotesize]
\begin{axis}[xmax=8, xtick={8}]
\addplot+ coordinates {
(1.0, 0.012987012987012988)
(1.0, 0.025974025974025976)
(1.0, 0.03896103896103896)
(1.0, 0.05194805194805195)
(1.0, 0.06493506493506493)
(1.0, 0.07792207792207792)
(1.0, 0.09090909090909091)
(1.0, 0.1038961038961039)
(1.0, 0.11688311688311688)
(1.0, 0.12987012987012986)
(1.0, 0.14285714285714285)
(1.0, 0.15584415584415584)
(1.0008550732891703, 0.16883116883116883)
(1.0010549774084303, 0.18181818181818182)
(1.0014000245692518, 0.19480519480519481)
(1.0017318952396594, 0.2077922077922078)
(1.002141742295396, 0.22077922077922077)
(1.002421530302084, 0.23376623376623376)
(1.0031187241873625, 0.24675324675324675)
(1.0038626455670165, 0.2597402597402597)
(1.0040112035641793, 0.2727272727272727)
(1.0061555889473066, 0.2857142857142857)
(1.0096132031252611, 0.2987012987012987)
(1.0108347025101643, 0.3116883116883117)
(1.0127530846161021, 0.3246753246753247)
(1.0151156164160138, 0.33766233766233766)
(1.018683154216385, 0.35064935064935066)
(1.0209642045054885, 0.36363636363636365)
(1.0373846168667198, 0.37662337662337664)
(1.0592242628059598, 0.38961038961038963)
(1.0614470569207528, 0.4025974025974026)
(1.0684769472829063, 0.4155844155844156)
(1.0890309046455515, 0.42857142857142855)
(1.1226389710541216, 0.44155844155844154)
(1.1244139828098596, 0.45454545454545453)
(1.1363996616910474, 0.4675324675324675)
(1.139271552075885, 0.4805194805194805)
(1.1410060371931967, 0.4935064935064935)
(1.1413288279685692, 0.5064935064935064)
(1.2261758979348347, 0.5194805194805194)
(1.2423199325883911, 0.5324675324675324)
(1.3215654843460076, 0.5454545454545454)
(2.1330699031644893, 0.5584415584415584)
(2.2617228061023424, 0.5714285714285714)
(2.5906010065644365, 0.5844155844155844)
(3.9763344339807607, 0.5974025974025974)
(4.005595024024893, 0.6103896103896104)
(5.738577618987365, 0.6233766233766234)
(6.658047132383878, 0.6363636363636364)
(20.162130335449895, 0.6493506493506493)
(39.50687330694465, 0.6623376623376623)
(79.0137466138893, 0.6753246753246753)
(79.0137466138893, 0.6883116883116883)
(79.0137466138893, 0.7012987012987013)
(79.0137466138893, 0.7142857142857143)
(79.0137466138893, 0.7272727272727273)
(79.0137466138893, 0.7402597402597403)
(79.0137466138893, 0.7532467532467533)
(79.0137466138893, 0.7662337662337663)
(79.0137466138893, 0.7792207792207793)
(79.0137466138893, 0.7922077922077922)
(79.0137466138893, 0.8051948051948052)
(79.0137466138893, 0.8181818181818182)
(79.0137466138893, 0.8311688311688312)
(79.0137466138893, 0.8441558441558441)
(79.0137466138893, 0.8571428571428571)
(79.0137466138893, 0.8701298701298701)
(79.0137466138893, 0.8831168831168831)
(79.0137466138893, 0.8961038961038961)
(79.0137466138893, 0.9090909090909091)
(79.0137466138893, 0.922077922077922)
(79.0137466138893, 0.935064935064935)
(79.0137466138893, 0.948051948051948)
(79.0137466138893, 0.961038961038961)
(79.0137466138893, 0.974025974025974)
(79.0137466138893, 0.987012987012987)
(79.0137466138893, 1.0)
(79.0137466138893, 1.0129870129870129)
} node[below right] at (4.5, 0.61038961) {off};
\addplot+ coordinates {
(1.0, 0.012987012987012988)
(1.0, 0.025974025974025976)
(1.0, 0.03896103896103896)
(1.0, 0.05194805194805195)
(1.0, 0.06493506493506493)
(1.0, 0.07792207792207792)
(1.0, 0.09090909090909091)
(1.0, 0.1038961038961039)
(1.0, 0.11688311688311688)
(1.0, 0.12987012987012986)
(1.0, 0.14285714285714285)
(1.0, 0.15584415584415584)
(1.0, 0.16883116883116883)
(1.0, 0.18181818181818182)
(1.0, 0.19480519480519481)
(1.0, 0.2077922077922078)
(1.0, 0.22077922077922077)
(1.0, 0.23376623376623376)
(1.0, 0.24675324675324675)
(1.0, 0.2597402597402597)
(1.0, 0.2727272727272727)
(1.0, 0.2857142857142857)
(1.0, 0.2987012987012987)
(1.0, 0.3116883116883117)
(1.0, 0.3246753246753247)
(1.0, 0.33766233766233766)
(1.0, 0.35064935064935066)
(1.0, 0.36363636363636365)
(1.0, 0.37662337662337664)
(1.0, 0.38961038961038963)
(1.0, 0.4025974025974026)
(1.0, 0.4155844155844156)
(1.0, 0.42857142857142855)
(1.0, 0.44155844155844154)
(1.0, 0.45454545454545453)
(1.0, 0.4675324675324675)
(1.0, 0.4805194805194805)
(1.0, 0.4935064935064935)
(1.0, 0.5064935064935064)
(1.0, 0.5194805194805194)
(1.0, 0.5324675324675324)
(1.0, 0.5454545454545454)
(1.0, 0.5584415584415584)
(1.0, 0.5714285714285714)
(1.0, 0.5844155844155844)
(1.0, 0.5974025974025974)
(1.0, 0.6103896103896104)
(1.0, 0.6233766233766234)
(1.0, 0.6363636363636364)
(1.0, 0.6493506493506493)
(1.0, 0.6623376623376623)
(1.0, 0.6753246753246753)
(1.0, 0.6883116883116883)
(1.0, 0.7012987012987013)
(1.0, 0.7142857142857143)
(1.0, 0.7272727272727273)
(1.0, 0.7402597402597403)
(1.0, 0.7532467532467533)
(1.0, 0.7662337662337663)
(1.0, 0.7792207792207793)
(1.0, 0.7922077922077922)
(1.0, 0.8051948051948052)
(1.0, 0.8181818181818182)
(1.0, 0.8311688311688312)
(1.0, 0.8441558441558441)
(1.0003788359033186, 0.8571428571428571)
(1.0009517674677446, 0.8701298701298701)
(1.0022386449556642, 0.8831168831168831)
(1.0036786284746515, 0.8961038961038961)
(1.0054708653420503, 0.9090909090909091)
(1.0061929847960693, 0.922077922077922)
(1.006260230230252, 0.935064935064935)
(1.0350684865632789, 0.948051948051948)
(1.0565213025146172, 0.961038961038961)
(1.3083087425140294, 0.974025974025974)
(1.3354293337074559, 0.987012987012987)
(79.0137466138893, 1.0)
(79.0137466138893, 1.0129870129870129)
} node[below right] at (4.5, 0.9870129) {on};
\end{axis}
\end{tikzpicture}
\caption{MSD - execution time.}
\label{fig:disagg:msdtime}
\end{subfigure}
\\[0.7em]
\begin{subfigure}{0.48\textwidth}
\centering\begin{tikzpicture}[font = \footnotesize]
\begin{axis}[xmax=12, xtick={12}]
\addplot+ coordinates {
(1.0, 0.013888888888888888)
(1.0, 0.027777777777777776)
(1.0, 0.041666666666666664)
(1.0, 0.05555555555555555)
(1.0, 0.06944444444444445)
(1.0, 0.08333333333333333)
(1.0, 0.09722222222222222)
(1.0, 0.1111111111111111)
(1.25, 0.125)
(1.25, 0.1388888888888889)
(1.2875, 0.1527777777777778)
(1.3333333333333333, 0.16666666666666666)
(1.3333333333333333, 0.18055555555555555)
(1.3333333333333333, 0.19444444444444445)
(1.336734693877551, 0.20833333333333334)
(1.3571428571428572, 0.2222222222222222)
(1.5, 0.2361111111111111)
(2.3333333333333335, 0.25)
(2.6666666666666665, 0.2638888888888889)
(2.75, 0.2777777777777778)
(3.0, 0.2916666666666667)
(3.5, 0.3055555555555556)
(3.6, 0.3194444444444444)
(4.285714285714286, 0.3333333333333333)
(4.5, 0.3472222222222222)
(4.5, 0.3611111111111111)
(8.0, 0.375)
(9.0, 0.3888888888888889)
(9.25, 0.4027777777777778)
(10.444444444444445, 0.4166666666666667)
(23.583333333333332, 0.4305555555555556)
(28.0, 0.4444444444444444)
(43.0, 0.4583333333333333)
(86.0, 0.4722222222222222)
(86.0, 0.4861111111111111)
(86.0, 0.5)
(86.0, 0.5138888888888888)
(86.0, 0.5277777777777778)
(86.0, 0.5416666666666666)
(86.0, 0.5555555555555556)
(86.0, 0.5694444444444444)
(86.0, 0.5833333333333334)
(86.0, 0.5972222222222222)
(86.0, 0.6111111111111112)
(86.0, 0.625)
(86.0, 0.6388888888888888)
(86.0, 0.6527777777777778)
(86.0, 0.6666666666666666)
(86.0, 0.6805555555555556)
(86.0, 0.6944444444444444)
(86.0, 0.7083333333333334)
(86.0, 0.7222222222222222)
(86.0, 0.7361111111111112)
(86.0, 0.75)
(86.0, 0.7638888888888888)
(86.0, 0.7777777777777778)
(86.0, 0.7916666666666666)
(86.0, 0.8055555555555556)
(86.0, 0.8194444444444444)
(86.0, 0.8333333333333334)
(86.0, 0.8472222222222222)
(86.0, 0.8611111111111112)
(86.0, 0.875)
(86.0, 0.8888888888888888)
(86.0, 0.9027777777777778)
(86.0, 0.9166666666666666)
(86.0, 0.9305555555555556)
(86.0, 0.9444444444444444)
(86.0, 0.9583333333333334)
(86.0, 0.9722222222222222)
(86.0, 0.9861111111111112)
(86.0, 1.0)
(86.0, 1.0138888888888888)
} node[below right] at (6.5, 0.36111111) {off};
\addplot+ coordinates {
(1.0, 0.013888888888888888)
(1.0, 0.027777777777777776)
(1.0, 0.041666666666666664)
(1.0, 0.05555555555555555)
(1.0, 0.06944444444444445)
(1.0, 0.08333333333333333)
(1.0, 0.09722222222222222)
(1.0, 0.1111111111111111)
(1.0, 0.125)
(1.0, 0.1388888888888889)
(1.0, 0.1527777777777778)
(1.0, 0.16666666666666666)
(1.0, 0.18055555555555555)
(1.0, 0.19444444444444445)
(1.0, 0.20833333333333334)
(1.0, 0.2222222222222222)
(1.0, 0.2361111111111111)
(1.0, 0.25)
(1.0, 0.2638888888888889)
(1.0, 0.2777777777777778)
(1.0, 0.2916666666666667)
(1.0, 0.3055555555555556)
(1.0, 0.3194444444444444)
(1.0, 0.3333333333333333)
(1.0, 0.3472222222222222)
(1.0, 0.3611111111111111)
(1.0, 0.375)
(1.0, 0.3888888888888889)
(1.0, 0.4027777777777778)
(1.0, 0.4166666666666667)
(1.0, 0.4305555555555556)
(1.0, 0.4444444444444444)
(1.0, 0.4583333333333333)
(1.0, 0.4722222222222222)
(1.0, 0.4861111111111111)
(1.0, 0.5)
(1.0, 0.5138888888888888)
(1.0, 0.5277777777777778)
(1.0, 0.5416666666666666)
(1.0, 0.5555555555555556)
(1.0, 0.5694444444444444)
(1.0, 0.5833333333333334)
(1.0, 0.5972222222222222)
(1.0, 0.6111111111111112)
(1.0, 0.625)
(1.0, 0.6388888888888888)
(1.0, 0.6527777777777778)
(1.0, 0.6666666666666666)
(1.0, 0.6805555555555556)
(1.0, 0.6944444444444444)
(1.0, 0.7083333333333334)
(1.0, 0.7222222222222222)
(1.0, 0.7361111111111112)
(1.0, 0.75)
(1.0, 0.7638888888888888)
(1.0, 0.7777777777777778)
(1.0, 0.7916666666666666)
(1.0, 0.8055555555555556)
(1.0, 0.8194444444444444)
(1.0, 0.8333333333333334)
(1.0, 0.8472222222222222)
(1.0, 0.8611111111111112)
(1.0, 0.875)
(1.0, 0.8888888888888888)
(1.0, 0.9027777777777778)
(1.0, 0.9166666666666666)
(1.0, 0.9305555555555556)
(1.0, 0.9444444444444444)
(1.0, 0.9583333333333334)
(1.0, 0.9722222222222222)
(1.0, 0.9861111111111112)
(1.0, 1.0)
(86.0, 1.0138888888888888)
} node[below right] at (6.5, 1.0) {on};
\end{axis}
\end{tikzpicture}
\caption{Iter - iteration count.}
\label{fig:disagg:iteriters}
\end{subfigure}%
\hspace{0.04\textwidth}%
\begin{subfigure}{0.48\textwidth}
\centering\begin{tikzpicture}[font = \footnotesize]
\begin{axis}[xmax=80, xtick={80}]
\addplot+ coordinates {
(1.0, 0.012987012987012988)
(1.0, 0.025974025974025976)
(1.0, 0.03896103896103896)
(1.0, 0.05194805194805195)
(1.0, 0.06493506493506493)
(1.0, 0.07792207792207792)
(1.0, 0.09090909090909091)
(1.0, 0.1038961038961039)
(1.0, 0.11688311688311688)
(1.0, 0.12987012987012986)
(1.0, 0.14285714285714285)
(1.0, 0.15584415584415584)
(1.0, 0.16883116883116883)
(1.0, 0.18181818181818182)
(1.0, 0.19480519480519481)
(1.0, 0.2077922077922078)
(1.0, 0.22077922077922077)
(1.0, 0.23376623376623376)
(1.0230127979335446, 0.24675324675324675)
(1.1428571428571428, 0.2597402597402597)
(1.1538461538461537, 0.2727272727272727)
(1.1666666666666667, 0.2857142857142857)
(1.2372491638795986, 0.2987012987012987)
(1.3513513513513513, 0.3116883116883117)
(1.4102564102564104, 0.3246753246753247)
(1.571150097465887, 0.33766233766233766)
(1.7407407407407407, 0.35064935064935066)
(1.792655772351587, 0.36363636363636365)
(1.9526011560693641, 0.37662337662337664)
(2.0785413744740535, 0.38961038961038963)
(2.8181818181818183, 0.4025974025974026)
(2.956022944550669, 0.4155844155844156)
(3.0196918588750115, 0.42857142857142855)
(3.4545454545454546, 0.44155844155844154)
(3.6578947368421053, 0.45454545454545453)
(5.133177570093458, 0.4675324675324675)
(7.909090909090909, 0.4805194805194805)
(7.929713455149502, 0.4935064935064935)
(8.48913043478261, 0.5064935064935064)
(8.928571428571429, 0.5194805194805194)
(11.887218045112782, 0.5324675324675324)
(12.87987987987988, 0.5454545454545454)
(12.94862227324914, 0.5584415584415584)
(16.687599364069953, 0.5714285714285714)
(22.816650148662042, 0.5844155844155844)
(40.86959328599096, 0.5974025974025974)
(60.659303882195445, 0.6103896103896104)
(64.13348008805283, 0.6233766233766234)
(64.21348314606742, 0.6363636363636364)
(68.18172305545866, 0.6493506493506493)
(303.4278059785674, 0.6623376623376623)
(606.8556119571348, 0.6753246753246753)
(606.8556119571348, 0.6883116883116883)
(606.8556119571348, 0.7012987012987013)
(606.8556119571348, 0.7142857142857143)
(606.8556119571348, 0.7272727272727273)
(606.8556119571348, 0.7402597402597403)
(606.8556119571348, 0.7532467532467533)
(606.8556119571348, 0.7662337662337663)
(606.8556119571348, 0.7792207792207793)
(606.8556119571348, 0.7922077922077922)
(606.8556119571348, 0.8051948051948052)
(606.8556119571348, 0.8181818181818182)
(606.8556119571348, 0.8311688311688312)
(606.8556119571348, 0.8441558441558441)
(606.8556119571348, 0.8571428571428571)
(606.8556119571348, 0.8701298701298701)
(606.8556119571348, 0.8831168831168831)
(606.8556119571348, 0.8961038961038961)
(606.8556119571348, 0.9090909090909091)
(606.8556119571348, 0.922077922077922)
(606.8556119571348, 0.935064935064935)
(606.8556119571348, 0.948051948051948)
(606.8556119571348, 0.961038961038961)
(606.8556119571348, 0.974025974025974)
(606.8556119571348, 0.987012987012987)
(606.8556119571348, 1.0)
(606.8556119571348, 1.0129870129870129)
} node[below right] at (40.5, 0.5844155) {off};
\addplot+ coordinates {
(1.0, 0.012987012987012988)
(1.0, 0.025974025974025976)
(1.0, 0.03896103896103896)
(1.0, 0.05194805194805195)
(1.0, 0.06493506493506493)
(1.0, 0.07792207792207792)
(1.0, 0.09090909090909091)
(1.0, 0.1038961038961039)
(1.0, 0.11688311688311688)
(1.0, 0.12987012987012986)
(1.0, 0.14285714285714285)
(1.0, 0.15584415584415584)
(1.0, 0.16883116883116883)
(1.0, 0.18181818181818182)
(1.0, 0.19480519480519481)
(1.0, 0.2077922077922078)
(1.0, 0.22077922077922077)
(1.0, 0.23376623376623376)
(1.0, 0.24675324675324675)
(1.0, 0.2597402597402597)
(1.0, 0.2727272727272727)
(1.0, 0.2857142857142857)
(1.0, 0.2987012987012987)
(1.0, 0.3116883116883117)
(1.0, 0.3246753246753247)
(1.0, 0.33766233766233766)
(1.0, 0.35064935064935066)
(1.0, 0.36363636363636365)
(1.0, 0.37662337662337664)
(1.0, 0.38961038961038963)
(1.0, 0.4025974025974026)
(1.0, 0.4155844155844156)
(1.0, 0.42857142857142855)
(1.0, 0.44155844155844154)
(1.0, 0.45454545454545453)
(1.0, 0.4675324675324675)
(1.0, 0.4805194805194805)
(1.0, 0.4935064935064935)
(1.0, 0.5064935064935064)
(1.0, 0.5194805194805194)
(1.0, 0.5324675324675324)
(1.0, 0.5454545454545454)
(1.0, 0.5584415584415584)
(1.0, 0.5714285714285714)
(1.0, 0.5844155844155844)
(1.0, 0.5974025974025974)
(1.0, 0.6103896103896104)
(1.0, 0.6233766233766234)
(1.0, 0.6363636363636364)
(1.0, 0.6493506493506493)
(1.0, 0.6623376623376623)
(1.0, 0.6753246753246753)
(1.0, 0.6883116883116883)
(1.0, 0.7012987012987013)
(1.0, 0.7142857142857143)
(1.0, 0.7272727272727273)
(1.0, 0.7402597402597403)
(1.0, 0.7532467532467533)
(1.0, 0.7662337662337663)
(1.0, 0.7792207792207793)
(1.0, 0.7922077922077922)
(1.0, 0.8051948051948052)
(1.0, 0.8181818181818182)
(1.0, 0.8311688311688312)
(1.0, 0.8441558441558441)
(1.0, 0.8571428571428571)
(1.0, 0.8701298701298701)
(1.0, 0.8831168831168831)
(1.0, 0.8961038961038961)
(1.0, 0.9090909090909091)
(1.0, 0.922077922077922)
(1.0, 0.935064935064935)
(1.0, 0.948051948051948)
(1.375, 0.961038961038961)
(1.8571428571428572, 0.974025974025974)
(2.466666666666667, 0.987012987012987)
(606.8556119571348, 1.0)
(606.8556119571348, 1.0129870129870129)
} node[below right] at (40.5, 0.9870129) {on};
\end{axis}
\end{tikzpicture}
\caption{MSD - node count.}
\label{fig:disagg:msdnodes}
\end{subfigure}
\caption{$\K^*$ cut disaggregation performance profiles.}
\end{figure}

\Cref{tab:disagg:summ} summarizes the status counts and shifted geomeans of performance metrics on instance subsets, and the performance profiles \crefrange{fig:disagg:itertime}{fig:disagg:msdnodes} compare the execution times or iteration/node counts. 
For both the iterative and MSD methods, disaggregation improves performance on nearly every solved instance. 
For the iterative method, it enables the pure-certificate-based variant to converge on more than double the number of instances, and it more than halves the execution time and iteration count.
Without disaggregation, the MSD method manages to converge on many more instances than the iterative method. 
Disaggregation greatly improves the performance of the MSD method, though the comparison is not quite as striking as for the iterative method.

\begin{table}[!htb]
\footnotesize
\begin{tabular}{l c *{4}{S[table-format=2]} *{3}{S[table-format=2.2]} *{3}{S[table-format=2.2]} *{3}{S[table-format=2.2]}}
\toprule
& & \multicolumn{4}{l}{statuses} & \multicolumn{3}{l}{time (s)} & \multicolumn{3}{l}{subproblems} & \multicolumn{3}{l}{iters or nodes} 
\\
\cmidrule(lr){3-6} \cmidrule(lr){7-9} \cmidrule(lr){10-12} \cmidrule(lr){13-15}
& \rotatebox{90}{\rlap{disag}} & {co} & {li} & {er} & {ex} & {aco} & {tco} & {all} & {aco} & {tco} & {all} & {aco} & {tco} & {all} \\
\midrule
\multirow{2}{*}{\rotatebox[origin=c]{90}{Iter}} 
& off & 
33 & 10 & 52 & 0 & 10.95 & 11.11 & 17.20 & 11.63 & 12.41 & 11.82 & 11.95 & 12.74 & 12.51 \\
& on & 
72 & 1 & 21 & 1 & 4.47 & 6.47 & 7.17 & 4.04 & 5.36 & 4.23 & 4.25 & 5.53 & 4.46 \\
\addlinespace
\multirow{2}{*}{\rotatebox[origin=c]{90}{MSD}} 
& off & 
51 & 3 & 41 & 0 & 1.74 & 6.18 & 6.71 & 15.51 & 50.51 & 27.98 & {$\hphantom{0}70$} & {$613$} & {$261$} \\
& on & 
76 & 0 & 18 & 1 & 1.06 & 3.40 & 3.50 & 7.57 & 15.80 & 12.70 & {$\hphantom{0}36$} & {$438$} & {$348$} \\
\bottomrule
\end{tabular}
\caption{$\K^*$ cut disaggregation performance summary.}
\label{tab:disagg:summ}
\end{table}

\subsubsection{Certificate-Based Scaling}
\label{sec:exp:scale}

To test the efficacy of the $\K^*$ certificate cut scaling technique for an LP solver with a feasibility tolerance we describe in \cref{sec:guar:tol}, we run Pajarito using only certificate cuts (the `c' variant described in \cref{sec:exp:cuts}), with and without scaling.
We set a larger feasibility tolerance on these four Pajarito solvers ($\delta = 10^{-6}$ instead of $10^{-8}$, which we used for all other tests), to reduce the chance that any observed effects are caused by numerical issues near machine epsilon.

\begin{figure}[!htbp]
\begin{subfigure}{0.48\textwidth}
\centering\begin{tikzpicture}[font = \footnotesize]
\begin{axis}[xmax=2, xtick={2}]
\addplot+ coordinates {
(1.0, 0.014084507042253521)
(1.0, 0.028169014084507043)
(1.0, 0.04225352112676056)
(1.0, 0.056338028169014086)
(1.0, 0.07042253521126761)
(1.0, 0.08450704225352113)
(1.0, 0.09859154929577464)
(1.0, 0.11267605633802817)
(1.0, 0.1267605633802817)
(1.0, 0.14084507042253522)
(1.0, 0.15492957746478872)
(1.0, 0.16901408450704225)
(1.0, 0.18309859154929578)
(1.0, 0.19718309859154928)
(1.0, 0.2112676056338028)
(1.0, 0.22535211267605634)
(1.0, 0.23943661971830985)
(1.0, 0.2535211267605634)
(1.0, 0.2676056338028169)
(1.0, 0.28169014084507044)
(1.0, 0.29577464788732394)
(1.0, 0.30985915492957744)
(1.0, 0.323943661971831)
(1.0, 0.3380281690140845)
(1.0, 0.352112676056338)
(1.0, 0.36619718309859156)
(1.0, 0.38028169014084506)
(1.0, 0.39436619718309857)
(1.0, 0.4084507042253521)
(1.0, 0.4225352112676056)
(1.0000102293744093, 0.43661971830985913)
(1.000374106751856, 0.4507042253521127)
(1.0003829928297188, 0.4647887323943662)
(1.0003960262251892, 0.4788732394366197)
(1.000639062318694, 0.49295774647887325)
(1.0006393391580313, 0.5070422535211268)
(1.0014061391451605, 0.5211267605633803)
(1.0028191381274747, 0.5352112676056338)
(1.0034409978945944, 0.5492957746478874)
(1.0035109631524575, 0.5633802816901409)
(1.0037433430763625, 0.5774647887323944)
(1.0038794413691934, 0.5915492957746479)
(1.0051304563875902, 0.6056338028169014)
(1.007580705472373, 0.6197183098591549)
(1.008160063309274, 0.6338028169014085)
(1.0082620559180133, 0.647887323943662)
(1.012694965473071, 0.6619718309859155)
(1.0185208612965642, 0.676056338028169)
(1.0283248152676157, 0.6901408450704225)
(1.0296822385022564, 0.704225352112676)
(1.0312536058458128, 0.7183098591549296)
(1.034705098635625, 0.7323943661971831)
(1.0383598005028662, 0.7464788732394366)
(1.0434359805895792, 0.7605633802816901)
(1.046403664268208, 0.7746478873239436)
(1.0472531414759447, 0.7887323943661971)
(1.0472830452021222, 0.8028169014084507)
(1.0535292034881967, 0.8169014084507042)
(1.066360855456763, 0.8309859154929577)
(1.0792374667416544, 0.8450704225352113)
(1.0989770935259877, 0.8591549295774648)
(1.099853251096174, 0.8732394366197183)
(1.2986808400657006, 0.8873239436619719)
(2.5973616801314012, 0.9014084507042254)
(2.5973616801314012, 0.9154929577464789)
(2.5973616801314012, 0.9295774647887324)
(2.5973616801314012, 0.9436619718309859)
(2.5973616801314012, 0.9577464788732394)
(2.5973616801314012, 0.971830985915493)
(2.5973616801314012, 0.9859154929577465)
(2.5973616801314012, 1.0)
(2.5973616801314012, 1.0140845070422535)
} node[below right] at (1.05, 0.83) {off};
\addplot+ coordinates {
(1.0, 0.014084507042253521)
(1.0, 0.028169014084507043)
(1.0, 0.04225352112676056)
(1.0, 0.056338028169014086)
(1.0, 0.07042253521126761)
(1.0, 0.08450704225352113)
(1.0, 0.09859154929577464)
(1.0, 0.11267605633802817)
(1.0, 0.1267605633802817)
(1.0, 0.14084507042253522)
(1.0, 0.15492957746478872)
(1.0, 0.16901408450704225)
(1.0, 0.18309859154929578)
(1.0, 0.19718309859154928)
(1.0, 0.2112676056338028)
(1.0, 0.22535211267605634)
(1.0, 0.23943661971830985)
(1.0, 0.2535211267605634)
(1.0, 0.2676056338028169)
(1.0, 0.28169014084507044)
(1.0, 0.29577464788732394)
(1.0, 0.30985915492957744)
(1.0, 0.323943661971831)
(1.0, 0.3380281690140845)
(1.0, 0.352112676056338)
(1.0, 0.36619718309859156)
(1.0, 0.38028169014084506)
(1.0, 0.39436619718309857)
(1.0, 0.4084507042253521)
(1.0, 0.4225352112676056)
(1.0, 0.43661971830985913)
(1.0, 0.4507042253521127)
(1.0, 0.4647887323943662)
(1.0, 0.4788732394366197)
(1.0, 0.49295774647887325)
(1.0, 0.5070422535211268)
(1.0, 0.5211267605633803)
(1.0, 0.5352112676056338)
(1.0, 0.5492957746478874)
(1.0, 0.5633802816901409)
(1.0, 0.5774647887323944)
(1.0000148477160564, 0.5915492957746479)
(1.000171776270362, 0.6056338028169014)
(1.000322146179331, 0.6197183098591549)
(1.0004345991783967, 0.6338028169014085)
(1.0006934458810874, 0.647887323943662)
(1.0010184322754987, 0.6619718309859155)
(1.001254055655628, 0.676056338028169)
(1.0017251288334548, 0.6901408450704225)
(1.002280044949674, 0.704225352112676)
(1.0039909002624328, 0.7183098591549296)
(1.0042510678082335, 0.7323943661971831)
(1.0055721003712397, 0.7464788732394366)
(1.0064828646463944, 0.7605633802816901)
(1.0082046909787366, 0.7746478873239436)
(1.0097409896888765, 0.7887323943661971)
(1.0097879004001824, 0.8028169014084507)
(1.00990778138025, 0.8169014084507042)
(1.0102384388553594, 0.8309859154929577)
(1.020670290773735, 0.8450704225352113)
(1.025417516630579, 0.8591549295774648)
(1.0260662164134073, 0.8732394366197183)
(1.0269968425703626, 0.8873239436619719)
(1.0385522609145057, 0.9014084507042254)
(1.0476260132725803, 0.9154929577464789)
(1.0599080714152582, 0.9295774647887324)
(1.1477643303389515, 0.9436619718309859)
(1.151615366005451, 0.9577464788732394)
(1.2930622408394061, 0.971830985915493)
(2.5973616801314012, 0.9859154929577465)
(2.5973616801314012, 1.0)
(2.5973616801314012, 1.0140845070422535)
} node[above right] at (1.05, 0.91) {on};
\end{axis}
\end{tikzpicture}
\caption{Iter - execution time.}
\label{fig:scale:itertime}
\end{subfigure}%
\hspace{0.04\textwidth}%
\begin{subfigure}{0.48\textwidth}
\centering\begin{tikzpicture}[font = \footnotesize]
\begin{axis}[xmax=2, xtick={2}]
\addplot+ coordinates {
(1.0, 0.014084507042253521)
(1.0, 0.028169014084507043)
(1.0, 0.04225352112676056)
(1.0, 0.056338028169014086)
(1.0, 0.07042253521126761)
(1.0, 0.08450704225352113)
(1.0, 0.09859154929577464)
(1.0, 0.11267605633802817)
(1.0, 0.1267605633802817)
(1.0, 0.14084507042253522)
(1.0, 0.15492957746478872)
(1.0, 0.16901408450704225)
(1.0, 0.18309859154929578)
(1.0, 0.19718309859154928)
(1.0, 0.2112676056338028)
(1.0, 0.22535211267605634)
(1.0, 0.23943661971830985)
(1.0, 0.2535211267605634)
(1.0, 0.2676056338028169)
(1.0, 0.28169014084507044)
(1.0, 0.29577464788732394)
(1.0, 0.30985915492957744)
(1.0, 0.323943661971831)
(1.0, 0.3380281690140845)
(1.0, 0.352112676056338)
(1.0, 0.36619718309859156)
(1.0, 0.38028169014084506)
(1.0, 0.39436619718309857)
(1.0, 0.4084507042253521)
(1.0, 0.4225352112676056)
(1.0, 0.43661971830985913)
(1.0, 0.4507042253521127)
(1.0, 0.4647887323943662)
(1.0, 0.4788732394366197)
(1.0, 0.49295774647887325)
(1.0, 0.5070422535211268)
(1.0, 0.5211267605633803)
(1.0, 0.5352112676056338)
(1.0, 0.5492957746478874)
(1.0, 0.5633802816901409)
(1.0, 0.5774647887323944)
(1.0003559288443167, 0.5915492957746479)
(1.000503173483579, 0.6056338028169014)
(1.0005473656697053, 0.6197183098591549)
(1.000595319809059, 0.6338028169014085)
(1.0042167321148676, 0.647887323943662)
(1.0073650112075117, 0.6619718309859155)
(1.007377122987191, 0.676056338028169)
(1.0092359853210342, 0.6901408450704225)
(1.0103915648969126, 0.704225352112676)
(1.0146230822267066, 0.7183098591549296)
(1.0156853701389759, 0.7323943661971831)
(1.0182800166279244, 0.7464788732394366)
(1.0312495115993041, 0.7605633802816901)
(1.049054615194743, 0.7746478873239436)
(1.0558554001824125, 0.7887323943661971)
(1.0704480210314224, 0.8028169014084507)
(1.0720379284979813, 0.8169014084507042)
(1.1391792399026126, 0.8309859154929577)
(1.4488972933731088, 0.8450704225352113)
(2.8977945867462176, 0.8591549295774648)
(2.8977945867462176, 0.8732394366197183)
(2.8977945867462176, 0.8873239436619719)
(2.8977945867462176, 0.9014084507042254)
(2.8977945867462176, 0.9154929577464789)
(2.8977945867462176, 0.9295774647887324)
(2.8977945867462176, 0.9436619718309859)
(2.8977945867462176, 0.9577464788732394)
(2.8977945867462176, 0.971830985915493)
(2.8977945867462176, 0.9859154929577465)
(2.8977945867462176, 1.0)
(2.8977945867462176, 1.0140845070422535)
} node[below right] at (1.05, 0.81) {off};
\addplot+ coordinates {
(1.0, 0.014084507042253521)
(1.0, 0.028169014084507043)
(1.0, 0.04225352112676056)
(1.0, 0.056338028169014086)
(1.0, 0.07042253521126761)
(1.0, 0.08450704225352113)
(1.0, 0.09859154929577464)
(1.0, 0.11267605633802817)
(1.0, 0.1267605633802817)
(1.0, 0.14084507042253522)
(1.0, 0.15492957746478872)
(1.0, 0.16901408450704225)
(1.0, 0.18309859154929578)
(1.0, 0.19718309859154928)
(1.0, 0.2112676056338028)
(1.0, 0.22535211267605634)
(1.0, 0.23943661971830985)
(1.0, 0.2535211267605634)
(1.0, 0.2676056338028169)
(1.0, 0.28169014084507044)
(1.0, 0.29577464788732394)
(1.0, 0.30985915492957744)
(1.0, 0.323943661971831)
(1.0, 0.3380281690140845)
(1.0, 0.352112676056338)
(1.0, 0.36619718309859156)
(1.0, 0.38028169014084506)
(1.0, 0.39436619718309857)
(1.0, 0.4084507042253521)
(1.0, 0.4225352112676056)
(1.0000083804766802, 0.43661971830985913)
(1.0000707616359776, 0.4507042253521127)
(1.0004266959494528, 0.4647887323943662)
(1.0009879243450825, 0.4788732394366197)
(1.0019645590299984, 0.49295774647887325)
(1.003186088717428, 0.5070422535211268)
(1.004191076128712, 0.5211267605633803)
(1.0053916016295315, 0.5352112676056338)
(1.0060819677724873, 0.5492957746478874)
(1.0066957296815793, 0.5633802816901409)
(1.0071858931158273, 0.5774647887323944)
(1.0074255546958233, 0.5915492957746479)
(1.007626669229803, 0.6056338028169014)
(1.0077700222894073, 0.6197183098591549)
(1.0086884748886853, 0.6338028169014085)
(1.009893367507655, 0.647887323943662)
(1.0101657649389242, 0.6619718309859155)
(1.0112538989686954, 0.676056338028169)
(1.0122237740598854, 0.6901408450704225)
(1.0127629202494512, 0.704225352112676)
(1.0128297817426108, 0.7183098591549296)
(1.0165449771665438, 0.7323943661971831)
(1.0198571164685166, 0.7464788732394366)
(1.0203486783206899, 0.7605633802816901)
(1.0229614227672867, 0.7746478873239436)
(1.023833458671458, 0.7887323943661971)
(1.0298416471753349, 0.8028169014084507)
(1.0483637381171063, 0.8169014084507042)
(1.0542190397693367, 0.8309859154929577)
(1.054541911266375, 0.8450704225352113)
(1.066504780490537, 0.8591549295774648)
(1.080942397424563, 0.8732394366197183)
(1.1017160067429295, 0.8873239436619719)
(1.1158989603146101, 0.9014084507042254)
(1.1621345792947504, 0.9154929577464789)
(1.224115943668958, 0.9295774647887324)
(1.432142319973736, 0.9436619718309859)
(2.8977945867462176, 0.9577464788732394)
(2.8977945867462176, 0.971830985915493)
(2.8977945867462176, 0.9859154929577465)
(2.8977945867462176, 1.0)
(2.8977945867462176, 1.0140845070422535)
} node[above right] at (1.05, 0.885) {on};
\end{axis}
\end{tikzpicture}
\caption{MSD - execution time.}
\label{fig:scale:msdtime}
\end{subfigure}
\\[0.7em]
\begin{subfigure}{0.48\textwidth}
\centering\begin{tikzpicture}[font = \footnotesize]
\begin{axis}[xmax=2, xtick={2}]
\addplot+ coordinates {
(1.0, 0.014084507042253521)
(1.0, 0.028169014084507043)
(1.0, 0.04225352112676056)
(1.0, 0.056338028169014086)
(1.0, 0.07042253521126761)
(1.0, 0.08450704225352113)
(1.0, 0.09859154929577464)
(1.0, 0.11267605633802817)
(1.0, 0.1267605633802817)
(1.0, 0.14084507042253522)
(1.0, 0.15492957746478872)
(1.0, 0.16901408450704225)
(1.0, 0.18309859154929578)
(1.0, 0.19718309859154928)
(1.0, 0.2112676056338028)
(1.0, 0.22535211267605634)
(1.0, 0.23943661971830985)
(1.0, 0.2535211267605634)
(1.0, 0.2676056338028169)
(1.0, 0.28169014084507044)
(1.0, 0.29577464788732394)
(1.0, 0.30985915492957744)
(1.0, 0.323943661971831)
(1.0, 0.3380281690140845)
(1.0, 0.352112676056338)
(1.0, 0.36619718309859156)
(1.0, 0.38028169014084506)
(1.0, 0.39436619718309857)
(1.0, 0.4084507042253521)
(1.0, 0.4225352112676056)
(1.0, 0.43661971830985913)
(1.0, 0.4507042253521127)
(1.0, 0.4647887323943662)
(1.0, 0.4788732394366197)
(1.0, 0.49295774647887325)
(1.0, 0.5070422535211268)
(1.0, 0.5211267605633803)
(1.0, 0.5352112676056338)
(1.0, 0.5492957746478874)
(1.0, 0.5633802816901409)
(1.0, 0.5774647887323944)
(1.0, 0.5915492957746479)
(1.0, 0.6056338028169014)
(1.0, 0.6197183098591549)
(1.0, 0.6338028169014085)
(1.0, 0.647887323943662)
(1.0, 0.6619718309859155)
(1.0, 0.676056338028169)
(1.0, 0.6901408450704225)
(1.0, 0.704225352112676)
(1.0, 0.7183098591549296)
(1.0, 0.7323943661971831)
(1.0, 0.7464788732394366)
(1.0, 0.7605633802816901)
(1.0769230769230769, 0.7746478873239436)
(1.12, 0.7887323943661971)
(1.125, 0.8028169014084507)
(1.2, 0.8169014084507042)
(1.3571428571428572, 0.8309859154929577)
(1.3846153846153846, 0.8450704225352113)
(1.4, 0.8591549295774648)
(1.5, 0.8732394366197183)
(1.5, 0.8873239436619719)
(3.0, 0.9014084507042254)
(3.0, 0.9154929577464789)
(3.0, 0.9295774647887324)
(3.0, 0.9436619718309859)
(3.0, 0.9577464788732394)
(3.0, 0.971830985915493)
(3.0, 0.9859154929577465)
(3.0, 1.0)
(3.0, 1.0140845070422535)
} node[below right] at (1.1, 0.795) {off};
\addplot+ coordinates {
(1.0, 0.014084507042253521)
(1.0, 0.028169014084507043)
(1.0, 0.04225352112676056)
(1.0, 0.056338028169014086)
(1.0, 0.07042253521126761)
(1.0, 0.08450704225352113)
(1.0, 0.09859154929577464)
(1.0, 0.11267605633802817)
(1.0, 0.1267605633802817)
(1.0, 0.14084507042253522)
(1.0, 0.15492957746478872)
(1.0, 0.16901408450704225)
(1.0, 0.18309859154929578)
(1.0, 0.19718309859154928)
(1.0, 0.2112676056338028)
(1.0, 0.22535211267605634)
(1.0, 0.23943661971830985)
(1.0, 0.2535211267605634)
(1.0, 0.2676056338028169)
(1.0, 0.28169014084507044)
(1.0, 0.29577464788732394)
(1.0, 0.30985915492957744)
(1.0, 0.323943661971831)
(1.0, 0.3380281690140845)
(1.0, 0.352112676056338)
(1.0, 0.36619718309859156)
(1.0, 0.38028169014084506)
(1.0, 0.39436619718309857)
(1.0, 0.4084507042253521)
(1.0, 0.4225352112676056)
(1.0, 0.43661971830985913)
(1.0, 0.4507042253521127)
(1.0, 0.4647887323943662)
(1.0, 0.4788732394366197)
(1.0, 0.49295774647887325)
(1.0, 0.5070422535211268)
(1.0, 0.5211267605633803)
(1.0, 0.5352112676056338)
(1.0, 0.5492957746478874)
(1.0, 0.5633802816901409)
(1.0, 0.5774647887323944)
(1.0, 0.5915492957746479)
(1.0, 0.6056338028169014)
(1.0, 0.6197183098591549)
(1.0, 0.6338028169014085)
(1.0, 0.647887323943662)
(1.0, 0.6619718309859155)
(1.0, 0.676056338028169)
(1.0, 0.6901408450704225)
(1.0, 0.704225352112676)
(1.0, 0.7183098591549296)
(1.0, 0.7323943661971831)
(1.0, 0.7464788732394366)
(1.0, 0.7605633802816901)
(1.0, 0.7746478873239436)
(1.0, 0.7887323943661971)
(1.0, 0.8028169014084507)
(1.0, 0.8169014084507042)
(1.0, 0.8309859154929577)
(1.0, 0.8450704225352113)
(1.0, 0.8591549295774648)
(1.0, 0.8732394366197183)
(1.0, 0.8873239436619719)
(1.0, 0.9014084507042254)
(1.0, 0.9154929577464789)
(1.0, 0.9295774647887324)
(1.0, 0.9436619718309859)
(1.1666666666666667, 0.9577464788732394)
(1.25, 0.971830985915493)
(3.0, 0.9859154929577465)
(3.0, 1.0)
(3.0, 1.0140845070422535)
} node[below right] at (1.1, 0.957) {on};
\end{axis}
\end{tikzpicture}
\caption{Iter - iteration count.}
\label{fig:scale:iteriters}
\end{subfigure}%
\hspace{0.04\textwidth}%
\begin{subfigure}{0.48\textwidth}
\centering\begin{tikzpicture}[font = \footnotesize]
\begin{axis}[xmax=6, xtick={6}]
\addplot+ coordinates {
(1.0, 0.014084507042253521)
(1.0, 0.028169014084507043)
(1.0, 0.04225352112676056)
(1.0, 0.056338028169014086)
(1.0, 0.07042253521126761)
(1.0, 0.08450704225352113)
(1.0, 0.09859154929577464)
(1.0, 0.11267605633802817)
(1.0, 0.1267605633802817)
(1.0, 0.14084507042253522)
(1.0, 0.15492957746478872)
(1.0, 0.16901408450704225)
(1.0, 0.18309859154929578)
(1.0, 0.19718309859154928)
(1.0, 0.2112676056338028)
(1.0, 0.22535211267605634)
(1.0, 0.23943661971830985)
(1.0, 0.2535211267605634)
(1.0, 0.2676056338028169)
(1.0, 0.28169014084507044)
(1.0, 0.29577464788732394)
(1.0, 0.30985915492957744)
(1.0, 0.323943661971831)
(1.0, 0.3380281690140845)
(1.0, 0.352112676056338)
(1.0, 0.36619718309859156)
(1.0, 0.38028169014084506)
(1.0, 0.39436619718309857)
(1.0, 0.4084507042253521)
(1.0, 0.4225352112676056)
(1.0, 0.43661971830985913)
(1.0083333333333333, 0.4507042253521127)
(1.046456223942823, 0.4647887323943662)
(1.0606060606060606, 0.4788732394366197)
(1.0735294117647058, 0.49295774647887325)
(1.0754716981132075, 0.5070422535211268)
(1.1053702196908055, 0.5211267605633803)
(1.129740093736685, 0.5352112676056338)
(1.131578947368421, 0.5492957746478874)
(1.1400222257501191, 0.5633802816901409)
(1.1452513966480447, 0.5774647887323944)
(1.1538461538461537, 0.5915492957746479)
(1.1551339285714286, 0.6056338028169014)
(1.1718170580964153, 0.6197183098591549)
(1.1738035264483626, 0.6338028169014085)
(1.1848837209302325, 0.647887323943662)
(1.1913646617333924, 0.6619718309859155)
(1.1996587030716723, 0.676056338028169)
(1.2052173913043478, 0.6901408450704225)
(1.2195121951219512, 0.704225352112676)
(1.346935448815475, 0.7183098591549296)
(1.348977135980746, 0.7323943661971831)
(1.3827624309392266, 0.7464788732394366)
(1.3953319407864861, 0.7605633802816901)
(1.4285714285714286, 0.7746478873239436)
(1.571124763705104, 0.7887323943661971)
(1.5799434413748097, 0.8028169014084507)
(1.6515751931840696, 0.8169014084507042)
(3.0, 0.8309859154929577)
(3.1666666666666665, 0.8450704225352113)
(10.5, 0.8591549295774648)
(10.5, 0.8732394366197183)
(10.5, 0.8873239436619719)
(10.5, 0.9014084507042254)
(10.5, 0.9154929577464789)
(10.5, 0.9295774647887324)
(10.5, 0.9436619718309859)
(10.5, 0.9577464788732394)
(10.5, 0.971830985915493)
(10.5, 0.9859154929577465)
(10.5, 1.0)
(10.5, 1.0140845070422535)
} node[below right] at (1.6, 0.82) {off};
\addplot+ coordinates {
(1.0, 0.014084507042253521)
(1.0, 0.028169014084507043)
(1.0, 0.04225352112676056)
(1.0, 0.056338028169014086)
(1.0, 0.07042253521126761)
(1.0, 0.08450704225352113)
(1.0, 0.09859154929577464)
(1.0, 0.11267605633802817)
(1.0, 0.1267605633802817)
(1.0, 0.14084507042253522)
(1.0, 0.15492957746478872)
(1.0, 0.16901408450704225)
(1.0, 0.18309859154929578)
(1.0, 0.19718309859154928)
(1.0, 0.2112676056338028)
(1.0, 0.22535211267605634)
(1.0, 0.23943661971830985)
(1.0, 0.2535211267605634)
(1.0, 0.2676056338028169)
(1.0, 0.28169014084507044)
(1.0, 0.29577464788732394)
(1.0, 0.30985915492957744)
(1.0, 0.323943661971831)
(1.0, 0.3380281690140845)
(1.0, 0.352112676056338)
(1.0, 0.36619718309859156)
(1.0, 0.38028169014084506)
(1.0, 0.39436619718309857)
(1.0, 0.4084507042253521)
(1.0, 0.4225352112676056)
(1.0, 0.43661971830985913)
(1.0, 0.4507042253521127)
(1.0, 0.4647887323943662)
(1.0, 0.4788732394366197)
(1.0, 0.49295774647887325)
(1.0, 0.5070422535211268)
(1.0, 0.5211267605633803)
(1.0, 0.5352112676056338)
(1.0, 0.5492957746478874)
(1.0, 0.5633802816901409)
(1.0, 0.5774647887323944)
(1.0, 0.5915492957746479)
(1.0, 0.6056338028169014)
(1.0, 0.6197183098591549)
(1.0, 0.6338028169014085)
(1.0, 0.647887323943662)
(1.0, 0.6619718309859155)
(1.0, 0.676056338028169)
(1.0, 0.6901408450704225)
(1.0, 0.704225352112676)
(1.0, 0.7183098591549296)
(1.0, 0.7323943661971831)
(1.0, 0.7464788732394366)
(1.0, 0.7605633802816901)
(1.0799640610961365, 0.7746478873239436)
(1.0914060791602886, 0.7887323943661971)
(1.1181024215958715, 0.8028169014084507)
(1.124755700325733, 0.8169014084507042)
(1.1484375, 0.8309859154929577)
(1.17524115755627, 0.8450704225352113)
(1.2, 0.8591549295774648)
(1.2789590606156775, 0.8732394366197183)
(1.3157894736842106, 0.8873239436619719)
(2.088235294117647, 0.9014084507042254)
(2.2763948497854076, 0.9154929577464789)
(2.9, 0.9295774647887324)
(5.25, 0.9436619718309859)
(10.5, 0.9577464788732394)
(10.5, 0.971830985915493)
(10.5, 0.9859154929577465)
(10.5, 1.0)
(10.5, 1.0140845070422535)
} node[above right] at (1.6, 0.875) {on};
\end{axis}
\end{tikzpicture}
\caption{MSD - node count.}
\label{fig:scale:msdnodes}
\end{subfigure}
\caption{$\K^*$ certificate cut scaling performance profiles.}
\end{figure}
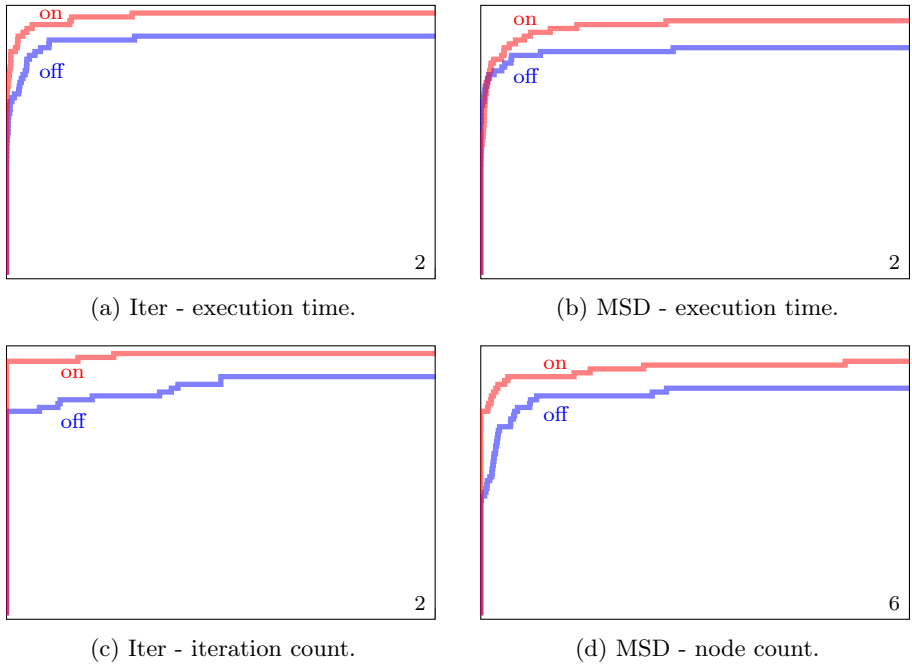

\Cref{tab:scale:summ} summarizes the status counts and shifted geomeans of performance metrics on instance subsets, and the performance profiles \crefrange{fig:scale:itertime}{fig:scale:msdnodes} compare the execution times or iteration/node counts. 
For both the iterative and MSD methods, using scaling improves the robustness of the pure-certificate-based variant, allowing us to converge on $6$ or $7$ additional instances. 
On the subset of instances solved by all four solvers (the `aco' columns), scaling slightly reduces conic subproblem counts and iteration or node counts, but has small and ambiguous effects on the execution times. 

\begin{table}[!htb]
\footnotesize
\begin{tabular}{l c *{4}{S[table-format=2]} *{3}{S[table-format=1.2]} *{3}{S[table-format=2.2]} *{3}{S[table-format=1.2]}}
\toprule
& & \multicolumn{4}{l}{statuses} & \multicolumn{3}{l}{time (s)} & \multicolumn{3}{l}{subproblems} & \multicolumn{3}{l}{iters or nodes} 
\\
\cmidrule(lr){3-6} \cmidrule(lr){7-9} \cmidrule(lr){10-12} \cmidrule(lr){13-15}
& \rotatebox{90}{\rlap{scale}} & {co} & {li} & {er} & {ex} & {aco} & {tco} & {all} & {aco} & {tco} & {all} & {aco} & {tco} & {all} \\
\midrule
\multirow{2}{*}{\rotatebox[origin=c]{90}{Iter}} 
& off & 
63 & 1 & 28 & 3 & 4.54 & 4.41 & 6.59 & 5.15 & 5.03 & 4.23 & 5.18 & 5.06 & 4.40 \\
& on & 
69 & 1 & 22 & 3 & 4.35 & 5.20 & 6.73 & 4.90 & 4.92 & 3.88 & 4.99 & 4.99 & 4.04 \\
\addlinespace
\multirow{2}{*}{\rotatebox[origin=c]{90}{MSD}} 
& off & 
60 & 0 & 30 & 5 & 2.68 & 2.77 & 3.15 & 12.44 & 14.48 & 12.78 & {$193$} & {$240$} & {$366$} \\
& on & 
67 & 0 & 26 & 2 & 2.92 & 4.02 & 3.86 & 11.88 & 15.77 & 12.07 & {$188$} & {$392$} & {$393$} \\
\bottomrule
\end{tabular}
\caption{$\K^*$ certificate cut scaling performance summary (larger $\delta$).}
\label{tab:scale:summ}
\end{table}

\appendix
\section{The Standard Primitive Nonpolyhedral Cones}
\label{sec:spec}

As we discussed in \cref{sec:soft:access}, Pajarito recognizes three standard primitive nonpolyhedral cones defined by MathProgBase: exponential, second-order, and positive semidefinite cones. 
Here, we tailor the general techniques for tightening OAs from \cref{sec:tight} to a primitive cone constraint involving one of these three cones. 
In particular, we describe the initial fixed OAs (see \cref{sec:tight:init}), extreme ray disaggregations (see \cref{sec:tight:extr}), and separation procedures (see \cref{sec:tight:sep}) implemented in Pajarito. 

These ideas could be adapted to other primitive nonpolyhedral cones if the user desires.
For example, consider a convex constraint (in NLP form) $f(\bt) \leq r$, where $f : \bbR^n \to \bbR$ is a non-homogeneous convex function. 
If we define a closed convex cone using the closure of the epigraph of the perspective of $f$:
\begin{equation}
\K_f = \cl \{ (r, s, \bt) \in \bbR^{2+n} : s > 0, r \geq s f ( \sfrac{\bt}{s} ) \}
\label{eq:kf}
,
\end{equation}
then the equivalent conic constraint is $(r, 1, \bt) \in \K_f$. 
From \citet{Zhang2014}, the dual cone of $\K_f$ is the closure of the perspective of the epigraph of the convex conjugate of $f$, $f^*(\bw) = \sup \{ -\bw^T \bt - f(\bt) : \bt \in \dom f \}$:
\begin{equation}
\K_f^* = \cl \{ (u, v, \bw) \in \bbR^{2+n} : u > 0, v \geq u f^* ( \sfrac{\bw}{u} ) \}
\label{eq:ksf}
.
\end{equation}
Note $\K_f^* \neq \K_{f^*}$ because of the permuting of the first two indices.

\subsection{Exponential Cone}
\label{sec:spec:exp}

The exponential cone $\KE$ is defined from the convex univariate exponential function $f(t) = \exp(t)$ in equation \cref{eq:kf}:
\begin{align}
\KE & = \cl \{ (r, s, t) \in \bbR^3 : s > 0, r \geq s \exp ( \sfrac{t}{s} ) \}
\\
& = \{ (r, 0, t) : r \geq 0, t \leq 0 \} \cup \{ (r, s, t) : s > 0, r \geq s \exp ( \sfrac{t}{s} ) \}
\notag
.
\end{align}
The convex conjugate is $f^*(w) = w - w \log(-w)$, so by equation \cref{eq:ksf}, the dual cone of the exponential cone is:
\begin{align}
\KE^* & = \cl \{ (u, v, w) \in \bbR^3 : u > 0, w < 0, v \geq w - w \log ( \sfrac{{-}w}{u} ) \}
\\
& = \{ (u, v, 0) : u, v \geq 0 \} \cup \{ (u, v, w) : u > 0, w < 0, v \geq w - w \log ( \sfrac{{-}w}{u} ) \}
\notag
.
\end{align}

\subsubsection{Initial Fixed Polyhedral Relaxation}
\label{sec:spec:exp:init}

Suppose we have a primitive cone constraint $(r, s, t) \in \KE$. 
We use the two $\KE^*$ extreme rays $(1, 0, 0), (0, 1, 0)$ to impose the simple bound constraints $r, s \geq 0$. 
We use more $\KE^*$ extreme rays of the form $(1, w - w \log (-w), w)$ by picking several different values $w < 0$. 
Note the corresponding cuts separate any point $(0, 0, t)$ satisfying $t > 0$.

\subsubsection{Extreme Ray Disaggregation}
\label{sec:spec:exp:extr}

Suppose we have the $\KE^*$ point $(u, v, w)$.
If $w = 0$, the point is already a nonnegative combination of the initial fixed $\KE^*$ points from \cref{sec:spec:exp:init}, so we discard it. 
If $w < 0$, we use the $\KE^*$ extreme ray $(u, w - w \log (\sfrac{{-}w}{u}), w)$, which when added to some nonnegative multiple of $(0, 1, 0)$, gives $(u, v, w)$.%
\footnote{This also projects $(u, v, w) \notin \KE^*$ with $u > 0, w < 0$ onto $\KE^*$.}

\subsubsection{Separation Of An Infeasible Point}
\label{sec:spec:exp:sep}

Suppose we want to separate a point $(r, s, t) \notin \KE$ that satisfies the initial fixed cuts. 
Then $r, s \geq 0$ and if $r = s = 0$ then $t \leq 0$. 
If $s = 0$, then $t > 0$ and $r > 0$, and we use the $\KE^*$ extreme ray $(\sfrac{t}{r}, -2 + 2 \log (\sfrac{2r}{t}), -2)$.
If $s > 0$, then $r < s \exp (\sfrac{t}{s})$, and we use the $\KE^*$ extreme ray $(1, (\sfrac{t}{s} - 1) \exp (\sfrac{t}{s}), -\exp (\sfrac{t}{s}))$.

\subsection{Second-Order Cone}
\label{sec:spec:soc}

For $n \geq 2$, the second-order cone is the epigraph of the $\ell_2$-norm, which is convex and homogeneous: 
\begin{equation}
\KL^{1+n} = \{ (r, \bt) \in \bbR^{1+n} : r \geq \lVert \bt \rVert_2 \}
.
\end{equation} 
We sometimes drop the dimension $1 + n$ when implied by context. This cone is self-dual ($\KL^* = \KL$). 
We also define the (self-dual) rotated second-order cone:
\begin{equation}
\KR^{2+n} = \{ (r, s, \bt) \in \bbR^{2+n} : r, s \geq 0, 2rs \geq \lVert \bt \rVert_2^2 \}
.
\end{equation} 
Note that $\KR$ is an invertible linear transformation of $\KL$, since $(r, s, \bt) \in \KR^{2+n}$ if and only if $(r + s, r - s, \sqrt{2} t_1, \ldots, \sqrt{2} t_n ) \in \KL^{2+n}$, so for simplicity we restrict attention to $\KL$.%
\footnote{As noted in \cref{sec:soft:access}, Pajarito transforms any $\KR$ constraints to equivalent $\KL$ constraints during preprocessing.}

\subsubsection{Initial Fixed Polyhedral Relaxation}
\label{sec:spec:soc:init}

Suppose we have a primitive cone constraint $(r, \bt) \in \KL^{1+n}$.
First, we note that the $\ell_{\infty}$-norm lower-bounds the $\ell_2$-norm, since for any $\bt \in \bbR^n$ we have:
\begin{equation}
\lVert \bt \rVert_{\infty} = \max_{\mathclap{i \in \iNs}} \, \lvert t_i \rvert \leq \lVert \bt \rVert_{2}
.
\end{equation}
Let $\be(i) \in \bbR^n$ be the $i$th unit vector in $n$ dimensions.
We use the $2n$ $\KL^*$ extreme rays $(1, \pm \be(i)), \forall i \in \iNs$, which imply the conditions $r \geq \lvert t_i \rvert, \forall i \in \iNs$, equivalent to the homogenized box relaxation $r \geq \lVert \bt \rVert_{\infty}$.
Second, we note that the $\ell_1$-norm also provides a lower bound for the $\ell_2$-norm, since for any $\bt \in \bbR^n$ we have:
\begin{equation}
\lVert \bt \rVert_1 = \sum_{\mathclap{i \in \iNs}} \, \lvert t_i \rvert \leq \sqrt{n} \lVert \bt \rVert_2
.
\end{equation}
We use the $2^n$ $\KL^*$ extreme rays $(1, \sfrac{\bm{\sigma}}{\sqrt{n}}), \forall \bm{\sigma} \in \{-1, 1\}^n$, which imply the homogenized diamond relaxation $r \geq \sfrac{\lVert \bt \rVert_1}{\sqrt{n}}$.
Although the number of initial fixed cuts is exponential in the dimension $n$, in \cref{sec:socef} we describe how to use an extended formulation introduced by \citet{Vielma2016} with $n$ auxiliary variables to imply an initial fixed OA that is no weaker but uses only a polynomial number of cuts.
Note that the $\KL^*$ point $(1, \bm{0})$, which corresponds to the simple variable bound $r \geq 0$, is a nontrivial conic combination of these initial fixed $\KL^*$ extreme rays.

\subsubsection{Extreme Ray Disaggregation}
\label{sec:spec:soc:extr}

Suppose we have the $\KL^*$ point $(u, \bw)$.
If $\bw = \bm{0}$, the point is already a nonnegative multiple of the $\KL^*$ point $(1, \bm{0})$, so we discard it.
Otherwise, we use the $\KL^*$ extreme ray $(\lVert \bw \rVert_2, \bw)$, which when added to some nonnegative multiple of $(1, \bm{0})$, gives the original point $(u, \bw)$.%
\footnote{This also projects $(u, \bw) \notin \KL^*$ with $u < \lVert \bw \rVert_2$ onto $\KL^*$.}

\subsubsection{Separation Of An Infeasible Point}
\label{sec:spec:soc:sep}

Suppose we want to separate a point $(r, \bt) \notin \KL$ that satisfies the initial fixed cuts.
Then $r \geq 0$ and so $\bt \neq \bm{0}$, and we use the $\KL^*$ extreme ray $(1, \sfrac{{-}\bt}{\lVert \bt \rVert_2})$.

\subsection{Positive Semidefinite Cone}
\label{sec:spec:psd}

For $n \geq 2$, we define the $n \times n$-dimensional positive semidefinite (PSD) matrix cone $\bbS^n_+$ as a subset of the symmetric matrices $\bbS^n = \{ \bT \in \bbR^{n \times n} : \bT = \bT^T \}$, avoiding the need to enforce symmetry constraints.
From the minimum eigenvalue function $\lambda_{\min} : \bbS^n \to \bbR$, we have:
\begin{equation}
\bbS^n_+ = \{ \bT \in \bbS^n : \lambda_{\min} (\bT) \geq 0 \}
.
\end{equation}
For $\bW, \bT \in \bbS^n$, we use the trace inner product $\langle \bW, \bT \rangle = \sum_{i,j \in \iNs} W_{i,j} T_{i,j}$.
$\bbS^n_+$ is self-dual and its extreme rays are the rank-$1$ PSD matrices \citep{Ben-Tal2001}, i.e. any $\bom \bom^T$ for $\bom \in \bbR^n$. 
An extreme ray $(\bbS^n_+)^*$ cut has the form:
\begin{equation}
\langle \bom \bom^T, \bT \rangle = \bom^T \bT \bom \geq 0
\label{cut:psdext}
.
\end{equation}
In \cref{sec:psdsoc}, we describe how to strengthen extreme ray $\bbS_+^*$ cuts to rotated second-order cone constraints (for MISOCP OA; see \cref{sec:soft:exten}).

Recall that our MI-conic form \ref{mod:micp} uses vector cone definitions, as does MathProgBase. 
\citet{MOSEK2016} refers to the matrix cone $\bbS^n_+ \subset \bbS^n$ as the \textit{smat PSD cone}, and to its equivalent vectorized definition $\KP^{\sfrac{n (n {+} 1)}{2}} \subset \bbR^{\sfrac{n (n {+} 1)}{2}}$ as the \textit{svec PSD cone}.
In svec space, we use the usual vector inner product, and $\KP$ is also self-dual. 
The invertible linear transformations for an smat-space point $\bT \in \bbS^n$ and an svec-space point $\bt \in \bbR^{\sfrac{n (n {+} 1)}{2}}$ are:
\begin{subequations}
\begin{align}
\svec (\bT) 
& = ( T_{1,1}, \sqrt{2} T_{2,1}, \ldots, \sqrt{2} T_{n,1}, T_{2,2}, \sqrt{2} T_{3,2}, \ldots, T_{n,n} ) 
,
\\
\smat (\bt) 
& = 
\begin{bmatrix}
t_1 & \sfrac{t_1}{\sqrt{2}} & \cdots & \sfrac{t_n}{\sqrt{2}}
\\
\sfrac{t_2}{\sqrt{2}} & t_{n+1} & \cdots & \sfrac{t_{2n-1}}{\sqrt{2}}
\\
\vdots & \vdots & \ddots & \vdots
\\
\sfrac{t_n}{\sqrt{2}} & \sfrac{t_{n-1}}{\sqrt{2}} & \cdots & t_{\sfrac{n (n {+} 1)}{2}}
\end{bmatrix}
.
\end{align}
\end{subequations}

\subsubsection{Initial Fixed Polyhedral Relaxation}
\label{sec:spec:psd:init}

Suppose we have a primitive cone constraint $\bT \in \bbS_+^n$.
Let $\be(i) \in \bbR^n$ be the $i$th unit vector in $n$ dimensions.
For each $i \in \iNs$, we let $\bom = \be(i)$ in the extreme ray $\bbS^*_+$ cut \cref{cut:psdext}, which imposes the diagonal nonnegativity condition $T_{i,i} \geq 0$ necessary for PSDness.
For each $i,j \in \iNs : i > j$, we let $\bom = \be(i) \pm \be(j)$ in \cref{cut:psdext}, which enforces the condition $T_{i,i} + T_{j,j} \geq 2 \lvert T_{i,j} \rvert$ necessary for PSDness.
\citet{Ahmadi2015} discuss an LP inner approximation of the PSD cone called the cone of diagonally dominant (DD) matrices. Our initial fixed $\bbS^*_+$ points $\bom \bom^T$ are exactly the extreme rays of the DD cone, so our initial fixed OA is the dual cone of the DD cone.

\subsubsection{Extreme Ray Disaggregation}
\label{sec:spec:psd:extr}

Suppose we have the $\bbS_+^*$ point $\bW$, not necessarily and extreme ray of $\bbS_+^*$. 
We perform an eigendecomposition $\bW = \sum_{i \in \iNs} \lambda_i \tilde{\bom}_i \tilde{\bom}_i^T$, where for all $i \in \iNs$, $\lambda_i$ is the $i$th eigenvalue and $\tilde{\bom}_i$ is its corresponding eigenvector.%
\footnote{Note that for real symmetric matrices, all eigenvalues are real. We select the eigenvectors to be orthonormal.}
Since $\bW$ is PSD, every eigenvalue is nonnegative, and there are $\rank (\bW) \leq n$ positive eigenvalues.
For each $i \in \iNs : \lambda_i > 0$, we let $\bom = \sqrt{\lambda_i} \tilde{\bom}_i$ in \cref{cut:psdext}.%
\footnote{By dropping any $i \in \iNs : \lambda_i < 0$, this projects $\bW \in \bbS \backslash \bbS_+^*$ onto $\bbS_+^*$.}
These extreme ray $\bbS^*_+$ cuts aggregate to imply the original $\bbS_+^*$ cut $\langle \bW, \bT \rangle \geq 0$.

\subsubsection{Separation Of An Infeasible Point}
\label{sec:spec:psd:sep}

Suppose we want to separate a point $\bT \in \bbS^n \backslash \bbS^n_+$. 
We perform an eigendecomposition $\bT = \sum_{i \in \iNs} \lambda_i \btau_i \btau_i^T$, for which at least one eigenvalue is negative.
For each $i \in \iNs : \lambda_i < 0$, we let $\bom = \btau_i$ in \cref{cut:psdext} (note $\langle \btau_i \btau_i^T, \bT \rangle = \btau_i^T \bT \btau_i = \lambda_i < 0$).

\section{The Second-Order Cone Extended Formulation}
\label{sec:socef}

Recall the definitions of the second-order cone $\KL$ and the rotated-second-order cone $\KR$ in \cref{sec:spec:soc}. Both $\KL$ and $\KR$ are self-dual. 
As discussed in \cref{sec:soft:exten}, Pajarito can optionally use an \textit{extended formulation} (EF) for second-order cone constraints, leading to tighter polyhedral relaxations.
\citet{Vielma2016} show that the constraint $(r, \bt) \in \KL^{1+n}$ is equivalent to the following $1 + n$ constraints on $r$, $\bt$, and the auxiliary variables $\bpi \in \bbR^n$:
\begin{subequations}
\begin{alignat}{3}
\sum_{\mathclap{i \in \iNs}} 2 \pi_i & \leq r 
\label{con:dis:lin}
\\
(r, \pi_i, t_i) & \in \KR^3 && \qquad & \forall i & \in \iNs 
\label{con:dis:kr}
.
\end{alignat}
\end{subequations}
By projecting out the $\bpi$ variables, the equivalence is obvious. 
Constraints \cref{con:dis:kr} imply $r \geq 0$ and $\pi_i \geq 0$ and $2 r \pi_i \geq t_i^2$ for all $i \in \iNs$. 
Aggregating the latter conditions and using the linear inequality \cref{con:dis:lin}, we see $r^2 \geq \sum_{i \in \iNs} 2 r \pi_i \geq \sum_{i \in \iNs} t_i^2$, which is equivalent to the original constraint (for $r \geq 0$).
We only use $\KL^{1+n}$ and $\KR^3$ here, so for convenience we drop the dimensions.%
\footnote{\citet{Ben-Tal2001a} introduced an alternative extended formulation for the $\KL^{1+n}$. See \citet{Vielma2016} for a discussion and computational comparison of various $\KL^{1+n}$ extended formulations in the context of a separation-based B\&B-OA algorithm for MISOCP.}

Suppose $(u, \bw)$ is a $\KL^*$ extreme ray, so from \cref{sec:spec:soc:extr}, we have $\bw \neq \bm{0}$ and $u = \lVert \bw \rVert_2 > 0$. 
Then $u r + \bw^T \bt \geq 0$ is a $\KL^*$ cut. 
Note that the linear constraint \cref{con:dis:lin} in the EF implies:
\begin{equation}
u r + \bw^T \bt \geq \frac{u r}{2} + u \sum_{\mathclap{i \in \iNs}} \pi_i + \bw^T \bt = \sum_{\mathclap{i \in \iNs}} \left( \frac{w_i^2 r}{2u} + u \pi_i + w_i t_i \right) 
\label{eq:socef:cuts}
.
\end{equation}
For each $i \in \iNs$, consider the $\KR^*$ extreme ray $(\sfrac{w_i^2}{2u}, u, w_i)$, which implies a $\KR^*$ cut for the $i$th constraint \cref{con:dis:kr} in the EF. 
The RHS of \cref{eq:socef:cuts} is an aggregation of these $n$ $\KR^*$ cuts, which means the $\KR^*$ cuts imply the $\KL^*$ cut condition $u r + \bw^T \bt \geq 0$. 
Therefore, there is no loss of strength in the polyhedral relaxations, and we maintain the certificate $\K^*$ cut guarantees from \cref{sec:guar:exact}.%
\footnote{Without the ability to rescale the linear constraint \cref{con:dis:lin}, we cannot recover the guarantees under an LP solver with a feasibility tolerance from \cref{sec:guar:tol}. However, Pajarito heuristically scales up each $\KR^*$ point by a factor $n$.}

We now apply this lifting procedure to the initial fixed $\KL^*$ points described in \cref{sec:spec:soc:init}.
The $\KL^*$ points for the $\ell_{\infty}$-norm relaxation are $(1, \pm \be(i)), \forall i \in \iNs$; for each $i \in \iNs$, we get three unique $\KR^*$ extreme rays $(0, 1, 0)$ (for $w_i = 0$) and $(\sfrac{1}{2}, 1, \pm 1)$ (for $w_i = \pm 1$).
The $\KL^*$ points for the $\ell_1$-norm relaxation are $(1, \sfrac{\bm{\sigma}}{\sqrt{n}}), \forall \bm{\sigma} \in \{-1, 1\}^n$; for each $i \in \iNs$, we get two unique $\KR^*$ extreme rays $(\sfrac{1}{2n}, 1, \sfrac{{\pm} 1}{\sqrt{n}})$ (for $w_i = \sfrac{{\pm} 1}{\sqrt{n}}$).
The polyhedral relaxation implied by these $5n$ $\KR^*$ points in the EF \cref{con:dis:kr,con:dis:lin} is at least as strong as that implied by the $2n + 2^n$ $\KL^*$ points from \cref{sec:spec:soc:init}, so our initial fixed OA can be imposed much more economically with the EF.

\section{SOCP Outer Approximation For PSD Cones}
\label{sec:psdsoc}

Recall the definitions of the self-dual smat-space PSD cone $\bbS_+$ in \cref{sec:spec:psd} and the self-dual rotated-second-order cone $\KR$ in \cref{sec:spec:soc}.%
\footnote{$\KR^3$ is in fact a simple linear transformation of $\bbS_+^2$.}
For a primitive cone constraint $\bT \in \bbS_+^n$, we demonstrate how to strengthen an $(\bbS^n_+)^*$ extreme ray cut $\langle \bom \bom^T, \bT \rangle \geq 0$ to up to $n$ different $\KR^3$ constraints.
As discussed in \cref{sec:soft:exten}, Pajarito can optionally solve an MISOCP OA model including these $\KR^3$ constraints, leading to tighter relaxations of a challenging $\bbS_+^n$ constraint.

Fix the index $i \in \iNs$. Let $\ubar{\omega} = \omega_i$ be the $i$th element of $\bom$, and $\ubar{\bom} = (\omega_j)_{j \in \iNs \backslash \{i\}} \in \bbR^{n-1}$ be the (column) subvector of $\bom$ with the $i$th element removed. 
Similarly, let $\ubar{t} = T_{i,i}$ and $\ubar{\bt} = (T_{i,j})_{j \in \iNs \backslash \{i\}} \in \bbR^{n-1}$, and let $\ubar{\bT} = (T_{k,j})_{k,j \in \iNs \backslash \{i\}} \in \bbS^{n-1}$ be the submatrix of $\bT$ with the $i$th column and row removed.
\citet{Kim2003} prove a variant of the standard Schur-complement result that $\bT \in \bbS^n_+$ if and only if $\ubar{t} \geq 0$ and:
\begin{subequations}
\begin{align}
\ubar{\bT} & \in \bbS^{n-1}_+
\label{con:psd:kim:2}
\\
\ubar{t} \ubar{\bT} - \ubar{\bt} \ubar{\bt}^T & \in \bbS^{n-1}_+
\label{con:psd:kim:3}
.
\end{align}
\end{subequations}
Consider the $3$-dimensional rotated second-order cone constraint:
\begin{equation}
( \ubar{t}, \ubar{\bom}^T \ubar{\bT} \ubar{\bom}, \sqrt{2} \ubar{\bom}^T \ubar{\bt} ) \in \KR^3
\label{con:psd:kr}
.
\end{equation}
By the definition of $\KR$, constraint \cref{con:psd:kr} is equivalent to the conditions $\ubar{t} \geq 0$ and:
\begin{subequations}
\begin{align}
\ubar{\bom}^T \ubar{\bT} \ubar{\bom} & \geq 0
\label{con:psd:kr:2}
\\
\ubar{\bom}^T (\ubar{t} \ubar{\bT}) \ubar{\bom} & \geq (\ubar{\bom}^T \ubar{\bt})^2 
\label{con:psd:kr:3}
.
\end{align}
\end{subequations}
Condition \cref{con:psd:kim:2} implies condition \cref{con:psd:kr:2}, by the dual cone definition \cref{eq:kd}. 
Since $(\ubar{\bom}^T \ubar{\bt})^2 = \ubar{\bom}^T \ubar{\bt} \ubar{\bt}^T \ubar{\bom}$, condition \cref{con:psd:kr:3} is equivalent to $\ubar{\bom}^T (\ubar{t} \ubar{\bT} - \ubar{\bt} \ubar{\bt}^T) \ubar{\bom} \geq 0$. 
Thus, by the dual cone definition again, condition \cref{con:psd:kim:3} implies condition \cref{con:psd:kr:3}. 
Therefore, constraint \cref{con:psd:kr} is a valid relaxation of the PSD constraint $\bT \in \bbS^n_+$. 
Furthermore, from Theorem 3.3 of \citet{Kim2003}, constraint \cref{con:psd:kr} holds if and only if:
\begin{equation}
\langle \bW, \bT \rangle \geq 0 \qquad \forall \bW \in \bbS^n_+ : ( W_{k,j} = \omega_k \omega_j, \forall k,j \in \iNs \backslash \{i\} )
.
\end{equation}
Thus constraint \cref{con:psd:kr} potentially implies an infinite family of $(\bbS^n_+)^*$ cuts, including the original $(\bbS^n_+)^*$ cut $\langle \bom \bom^T, \bT \rangle \geq 0$.%
\footnote{We do not explore how to scale these $\KR$ constraints to recover the guarantees from \cref{sec:guar:tol} for an SOCP solver with an absolute feasibility tolerance.}
Note that the choice of $i \in \iNs$ is arbitrary, so we can derive $n$ different $\KR$ constraints of the form \cref{con:psd:kr}.%
\footnote{For strengthening separation or certificate $\K^*$ cuts, Pajarito heuristically picks one of the $n$ possible $\KR$ constraints by choosing $i$ as the coordinate of the largest absolute value in $\bom$.}

We now apply this strengthening procedure to the initial fixed $(\bbS^n_+)^*$ extreme rays described in \cref{sec:spec:psd:init}. 
Letting $\bom = \be(i) \pm \be(j)$ for each $i,j \in \iNs : j > i$ in constraint \cref{con:psd:kr}, we get the $\sfrac{n (n {-} 1)}{2}$ initial fixed $\KR$ constraints:
\begin{equation}
( T_{i,i}, T_{j,j}, \sqrt{2} T_{i,j} ) \in \KR^3 \qquad \forall i,j \in \iNs : j > i
\label{con:psd:init}
.
\end{equation}
These constraints enforce that every $2 \times 2$ principal matrix of $\bT$ is PSD, a necessary but insufficient condition for $\bT \in \bbS^n_+$.
\citet{Ahmadi2015} discuss an SOCP inner approximation of the PSD cone called the cone of scaled diagonally dominant (SDD) matrices. 
Our initial fixed SOCP OA is the dual SDD matrix cone, a strict subset of the polyhedral dual DD matrix cone that corresponds to our initial fixed LP OA from \cref{sec:spec:psd:init}.


\newpage
\bibliographystyle{abbrvnat}
\bibliography{refs} 


\end{document}